\documentclass{article}
\usepackage[a4paper, total={6in, 10in}]{geometry}
\usepackage[utf8]{inputenc}
\usepackage{amsthm,mathtools,amsmath}
\usepackage{enumitem}
\usepackage{hyperref}
\usepackage[charter]{mathdesign}
\usepackage{color}
\usepackage{titlesec}
\usepackage{subcaption}

 \usepackage{ulem}

\titleformat{\paragraph}
{\normalfont\normalsize\bfseries}{\theparagraph}{1em}{}
\titlespacing*{\paragraph}
{0pt}{3.25ex plus 1ex minus .2ex}{1.5ex plus .2ex}

\DeclareFontFamily{U}{BOONDOX-calo}{\skewchar\font=45 }
\DeclareFontShape{U}{BOONDOX-calo}{m}{n}{
  <-> s*[1.05] BOONDOX-r-calo}{}
\DeclareFontShape{U}{BOONDOX-calo}{b}{n}{
  <-> s*[1.05] BOONDOX-b-calo}{}
\DeclareMathAlphabet{\mathcalboondox}{U}{BOONDOX-calo}{m}{n}
\SetMathAlphabet{\mathcalboondox}{bold}{U}{BOONDOX-calo}{b}{n}
\DeclareMathAlphabet{\mathbcalboondox}{U}{BOONDOX-calo}{b}{n}
\newcommand{\mcb}[1]{{\mathcalboondox #1}}

\makeatletter \@addtoreset{equation}{section}
\@addtoreset{figure}{section}
 \makeatother

\newtheorem{theorem}{Theorem}[section]
\newtheorem{corollary}{Corollary}[theorem]
\newtheorem{lemma}[theorem]{Lemma}
\newtheorem{proposition}[theorem]{Proposition}
\newtheorem{definition}{Definition}[section]
\newtheorem{remark}{Remark}[]

\DeclareMathOperator*{\maxx}{{\vphantom{p}\mathrm{max}}}

\newcounter{hypeq}

\usepackage{tikz}
\usetikzlibrary{shapes.misc}
\usetikzlibrary{arrows.meta}
\tikzset{cross/.style={cross out, draw=black, minimum size=2*(#1-\pgflinewidth), inner sep=0pt, outer sep=0pt},
cross/.default={0.2cm}}
\definecolor{blue-violet}{rgb}{0.54, 0.17, 0.89}
\usetikzlibrary{decorations.pathreplacing}

\newcommand{\lsim}{\lesssim}

\newenvironment{acknowledgements}{%
  \begin{abstract}
}{%
  \end{abstract}
}

\title{Non-equilibrium fluctuations for \\ SEP($\alpha$) with open boundary} 
\author{C. Franceschini \thanks{University of Modena and Reggio Emilia, FIM, Via G. Campi 213/B 41125, Modena, Italy. E-mail: {\tt chiara.fraceschini@unimore.it}}
\and P. Gon\c calves\thanks{Instituto Superior T\'ecnico, Department of Mathematics, Av. Rovisco Pais 1, 1049-001, Lisbon. E-mail: {\tt pgoncalves@tecnico.ulisboa.pt}}\and M. Jara\thanks{Instituto de Matemática Pura e Aplicada, Estrada Dona Castorina, no. 110, Rio de Janeiro, Brazil. E-mail: {\tt mdjara@gmail.com}}\and B. Salvador\thanks{Instituto Superior T\'ecnico, Department of Mathematics, Av. Rovisco Pais 1, 1049-001, Lisbon. E-mail: {\tt beatriz.salvador@tecnico.ulisboa.pt}}}

\date{ }

\begin{document}

\maketitle

\begin{abstract}{We analyze the non-equilibrium fluctuations of the partial symmetric simple exclusion process, SEP($\alpha$), which allows at most $\alpha \in \mathbb{N}$ particles per site, and we put it in contact with stochastic reservoirs whose strength is regulated by a parameter $\theta \in \mathbb{R}$. Setting $\alpha = 1$, we find the results of \cite{LMO,FGN19,GJMN} and extend the known results to cover all range of $\theta$.}
\end{abstract}



\section{Introduction}

Interacting particle systems are stochastic systems on which individual units (the so-called \textit{particles}) perform Markovian evolutions influenced by the presence of other particles. The objective is to study the emergence of collective behavior out of simple interaction rules for the individual units of the system.
Among the most studied interacting particle systems \cite{Liggett} is the so-called \textit{exclusion process}, on which the interaction between particles is reduced to a simple \textit{exclusion rule}, under which particles evolving on a graph can never share the same position. The exclusion model has been used as a landmark for a myriad of collective behavior, among which mass transport, interface growth and motion by mean curvature. The success of the exclusion process as an interacting particle system comes from one side from its striking combinatorial and algebraic properties, which makes the analysis of the collective behavior of particles a mathematically tractable problem, and from the other side from the fact that it is rich enough to allow modelling a great variety of collective behaviors. 
A generalization of the exclusion process that shares many of its algebraic properties is the so-called \textit{partial exclusion process}: in this model, the exclusion rule is relaxed to allow at most $\alpha$ particles per site, where $\alpha \in \mathbb N$ is a fixed parameter.

The partial exclusion process that we investigate here, the SEP($\alpha$), was first introduced in Section B of \cite{gunter}. We restrict ourselves to the choice of a simple symmetric dynamics on a one-dimensional lattice, i.e.~nearest-neighbor jumps with $p(1)=p(-1)=1/2$. For $N\in\mathbb N$, we consider the finite lattice $\Lambda_N =\{1,\dots,N-1\}$  which we call bulk. For a site $x\in\Lambda_N$, we fix the rate at which a particle jumps from  $x$ to $x+1$ (resp. from $x+1$ to $x$) to be equal to $\eta(x)(\alpha-\eta(x+1))$  (resp. $\eta(x+1)(\alpha-\eta(x))$), where $\eta(x)$ denotes the quantity of particles at site $x$ on the configuration $\eta$. If $\alpha=1$, the model coincides with the so-called symmetric simple exclusion process (SSEP).  This specific choice of the rates was introduced in \cite{gunter}, see equation (2.30) in that article. The SEP($\alpha$) has been further studied in other settings, such as in \cite{carinci2013duality} and \cite{duality_neg_corr} where the system is put in contact with stochastic reservoirs, in \cite{FRShydrodynamics2021} under a random enviroment and also in \cite{CGRexact2020}, \cite{ChenSau2021}, always from a duality point of view. We note that for the choice of rates given above this model is what is called a \textit{gradient model}, since the instantaneous current of the system at the bond $\{x,x+1\}$, i.e.~the difference between the jump rate from $x $ to $x+1$  and the jump rate from $x+1$ to $x$ can be written as the gradient of a local function. More precisely, that current is equal to $\alpha (\eta(x)-\eta(x+1))$. We also observe that the number of particles is conserved by the dynamics of the SEP($\alpha$) and that the symmetry of the jump rates of the individual particles makes the system reversible with respect to Binomial measures of product form. 

Non-equilibrium phenomena have become increasingly relevant in recent years, and the study of how collective behavior is modified by breaking reversibility is an active research subject. A natural way to modify the SEP($\alpha$) in order to make it non-reversible, is to attach to the lattice \textit{density reservoirs} with at least two different densities. This creates currents through the system, which \textit{drive} the system out of equilibrium. In this article, this will be the setting we will be working on, i.e. we will attach a stochastic reservoir to each boundary point of  $\Lambda_N$. These reservoirs will break the conservation of the total number of particles, since they can inject and remove particles, even-though the individual units of the system will still be conserved \textit{locally}. With the aim of exploring various possible answers to the question whether the limiting collective behavior of particles retains the non-reversible behavior, we will choose the particles injection and removal rates to scale with the size $N$ of the system, through a parameter $\theta \in \mathbb{R}$, and to be such that the system is no longer in equilibrium. When $\theta<0$, the reservoirs are fast and when $\theta\geq 0$, the reservoirs are slow. 

The main question here is whether this non-reversible behavior is observed at the level of the scaling limits of the model. The hydrodynamic limit of the SEP($\alpha$) turns out to be  a non-reversible PDE, which answers this question at the level of the law of large numbers. The next question is whether the non-reversible behavior has a stochastic component, which motivates the analysis of the fluctuations of the density around its hydrodynamic limit. The question can thus be restated as whether a non-reversible behavior is observed in the limiting SPDE. The Macroscopic Fluctuation Theory (MFT), as formulated in \cite{MFT1, MFT2} can be used to predict the behavior of large scale limits of driven-diffusive systems. This description depends on two macroscopic quantities, the \textit{diffusivity} and the \textit{mobility} of the system. One assumes that these quantities are local functions of the thermodynamic variables. In the case of the SEP($\alpha$), the density of particles $\rho \in [0,\alpha]$ is the only thermodynamic variable. The diffusivity is constant and equal to $\alpha$, while the mobility is quadratic and equal to $ \rho(\alpha - \rho)$. Our main result confirms the predictions of MFT for the Central Limit Theorem (CLT) fluctuations of the density of particles.

In this article, we will be interested on the analysis of the fluctuations of the density of particles around its hydrodynamic limit. This corresponds to the derivation of the CLT associated to the hydrodynamic limit of the system. The limiting equation is no longer a PDE, but a linear SPDE on which the time evolution is given by the hydrodynamic equation, plus a stochastic conservative noise with a covariance structure given in terms of solutions of the hydrodynamic equation. More precisely, in this paper we will analyse the non-equilibrium time dependent fluctuations for SEP($\alpha$) for all $\theta \in \mathbb{R}$ and $\alpha \in \mathbb{N}$. We remark that the equilibrium case can also be easily proved by the same type of arguments as in the case $\alpha=1$, obtained in \cite{FGN16}. For that reason, we omit the proof of this case here and we refer the reader to that article for a proof.

Now we recall the state-of-the-art of some of the scaling limits for this model. For the case of the exclusion process with open boundary and $\alpha=1$, the hydrodynamic limit was derived in \cite{baldasso2017exclusion} for slow reservoirs and in \cite{byronMPRF} for fast reservoirs. In \cite{FGS2022}, the derivation of the hydrodynamic limit was extended to $\alpha \in \mathbb{N}$ in both the slow and fast regimes, with a proof that relies on the entropy method introduced in \cite{GPV}. An extension of these hydrodynamic limits to general domains based on duality can be found in \cite{schiavo2021scaling}. The hydrodynamic equation of the  SEP($\alpha$) is the heat equation given by $\partial_{t} \rho_t(u) = \alpha\Delta \rho_t(u),$ that needs to be complemented with suitable boundary conditions. Depending on the choice of the parameter $\theta$, the boundary conditions are of Dirichlet type (for $\theta<1$), Robin type (for $\theta=1$) or Neumann type (for $\theta>1$). The non-equilibrium fluctuations for the case $\alpha=1$ were analysed in several works, namely in: \cite{LMO} when $\theta = 0$, where the non-equilibrium stationary fluctuations were derived as a consequence of its non-equilibrium fluctuations; \cite{FGN19} when $\theta = 1$ and \cite{GJMN} when $\theta\in[0,\infty)$. The equilibrium fluctuations, also for the case $\alpha = 1$, were analysed in \cite{FGN16} for $\theta \geq 0$. Nevertheless, the case $\theta < 0$ was an open problem up to now, apart in the equilibrium setting, which was derived in \cite{BGJS21}.

The main difficulty on the rigorous
mathematical derivation of the non-equilibrium fluctuations relies on the fact that the systems typically exhibit long-range space–time correlations. For that reason, one has to face the problem of obtaining good estimates of the two-point centered correlation function, that we denote by $\varphi^N_t$. This is one of the main topics discussed in this article and we consider that it is here that relies the major contribution of our work. For the case $\alpha = 1$, by writing down the Chapman-Kolmogorov equations directly for $\varphi^N_t$, one gets 
\begin{equation} \label{eq_1_varphi_intro}
\partial_t \varphi_t^N(x, y) 
		=  N^2 \Delta^{i}_N \varphi_t^N(x, y) 
			+ g^{N}_t(x,y) \mathbb {1} ((x,y) \in \mathcal{D}_N^{\pm}),
\end{equation} where $\Delta^i_N$ is the infinitesimal generator of a certain bi-dimensional random walk, $\mathcal{D}_N^{\pm}$ is a certain finite set that we will define later and $g^{N}_t$ is a non-positive function that only has support on $\mathcal{D}_N^{\pm}$. From last identity, one can use Duhamel's formula to obtain an expression for such function. From that, we reduce the problem to estimating three simple quantities: the initial correlations  $\varphi^N_0$, the term $g^{N}_t$ and the occupation time on  $\mathcal{D}_N^{\pm}$ of the bi-dimensional random walk with infinitesimal generator $\Delta^i_N$. Unfortunately, for $\alpha \geq 2$, if one tries to write down the Chapman-Kolmogorov equations directly for $\varphi^N_t$ defined as in the case $\alpha=1$, an additional interaction term appears at the diagonal $\{x=y\}$, which breaks down the previous approach. To overcome such issue, we construct an extension of $\varphi^N_t$ to the diagonal $\{x=y\}$, to which a similar approach as the one previously described can be applied to obtain the decay in $N$ of $\varphi^N_t$. By analyzing this extension function, we are able to obtain a generalization of the results  for $\alpha = 1$ that were derived in  \cite{LMO,GJMN,FGN19}. The novelty of our approach to obtain the decay in $N$ of $\varphi^N_t$ is the construction and use of such a well chosen extension function that can be compared with $\varphi^N_t$ and also the use of some discrete versions of the maximum principle (see Appendix \ref{appendix_maximum_principles}) to, after applying Duhamel's formula, compare occupation times for different values of $\theta$. After some trial and error, we discovered that the right choice of the extension function is related to the \textit{duality function} of the SEP($\alpha$), see \cite{carinci2013duality} and Remark \ref{remark_duality}.  Nevertheless, we observe that there are other ways on which one can arrive to the right extension function for the correlation function $\varphi^N_t$.  In order to follow a fully analytical method, for example, one can introduce a boundary layer at the diagonal to discover the best approximation of the heat equation with sources at the diagonal.

To determine the non-equilibrium fluctuations of the system we follow the same strategy outlined in \cite{LMO, FGN19, GJMN} (with similar ideas to the ones described in Chapter 11 of \cite{KL}), and, for that reason, some details in the proofs are omitted here. The idea of the argument is the classical probabilistic approach to functional convergence of stochastic processes, namely, to prove tightness of the sequence of density fluctuation fields and then characterize all limit points. If on top of the conditions that we will need to ask in order to prove tightness, we also ask that, at the initial time, the sequence of density fields converges to a mean-zero Gaussian process, then the convergence takes place for any time $t$ and the unique limiting process is a generalized Ornstein–Uhlenbeck process which is a solution of \eqref{OU_equation}. 

Now we comment on the main tools and difficulties of our approach. We first observe that depending on the range of $\theta$, the density fluctuation fields have to be defined on proper spaces of test functions, which typically are quite regular and satisfy the boundary conditions of the hydrodynamic equation but with an appropriate choice of parameters. Second,  in order to prove tightness, we use both Aldous and Kolmogorov-Centsov criteria (as in \cite{GJMN}), where this last one is mainly applied to the boundary integral terms of the Dynkin's martingales. Recall that on the proof of tightness at the level of the hydrodynamic limit, i.e.~of the sequence of empirical measures associated with the density profile, the quadratic variation of the Dynkin's martingale $\{M^N_t(\phi)\}_{N \in \mathbb{N}}$ converges to zero. Now, in the case of fluctuations, the corresponding Dynkin's martingale converges, as $N$ goes to infinity, in the $J_1$-Skorohod space $\mathcal{D}_N([0, T]; \mathbb{R})$ of c\`adl\`ag functions from $[0,T]$ to $\mathbb{R}$, to a mean-zero Gaussian process which is a
martingale with continuous trajectories and with a deterministic, non-degenerated quadratic variation. We also note that from our results we can  obtain the non-equilibrium fluctuations starting the process from a product measure with slowly varying parameter or even a constant one. In particular, if we fix a profile $\rho: [0,1] \to [0,1]$  and consider $\mu^N$ as the product measure whose marginals are given by the Binomial($\alpha$, $\rho(\tfrac{x}{N})$) distribution, the result also holds, leading to an Ornstein–Uhlenbeck process in the limit.

In our work, we also consider the case $\theta<0$ for $\alpha \in \mathbb{N}$ in the non-equilibrium scenario, extending therefore the results of \cite{BGJS21}. This case is more demanding than the others since the boundary terms are of order $O(N^{-\theta})$ and therefore, they blow up when taking $N\to+\infty$. To overcome this difficulty, we take a space of test functions that have all derivatives equal to zero at the boundary. Since this space of test functions is too little we supplement the characterization of limit points by showing that the limit field when integrated in time satisfies the Dirichlet conditions as in the case $\theta\in[0,1)$. This is reminiscent of item 2 (ii) of Theorem 2.13 of \cite{BGJS21}, where it was proved that when the system is in its equilibrium state, this extra  condition gives in fact the uniqueness of the limit. Here we extended that result to the non-equilibrium setting, though we lack a proof of uniqueness in that general case. 

Here is a summary of our contributions in this article.
First we provide a natural extension of the two-point correlation function to the diagonal in such a way that it satisfies a consistent set of equations that allows  estimating the non-stationary two-point correlations of the SEP($\alpha$) for any value of $\alpha \in \mathbb{N}$ and $\theta \in \mathbb{R}$. As a consequence, we characterize the non-equilibrium fluctuations of SEP($\alpha$) for any value of $\alpha \geq 2$ and $\theta \in \mathbb{R}$. Moreover, our approach also allows characterizing the non-equilibrium fluctuations of SEP($1$), for $\theta <0$.

To conclude we comment on the fluctuations starting from non-equilibrium stationary state (NESS). Observe that the Ornstein-Uhlenbeck equation \eqref{OU_equation} has a unique invariant measure, which is given by a Gaussian spatial process on the interval $[0,1]$. Observe as well that the SEP($\alpha$) as defined here is irreducible, and in particular has a unique invariant measure. A relevant question is the derivation of a fluctuation result for the empirical density of particles of the SEP($\alpha$) starting from its NESS. This question has been solved for the SSEP in \cite{LMO, GJMN}, and more recently in \cite{GJMM} for reaction-diffusion models. Unfortunately, our estimates are  not sharp enough to allow for the limit exchange which needed to derive such a result. Recall that, for $\alpha=1$, the matrix ansatz method (MPA) developed by \cite{derrida} provides detailed information about the NESS of SEP($1$) and recently \cite{FloreaniCasanova} found a characterization of such measure. For SEP($1$), the MPA enables one to obtain explicitly the n-point correlation function of the system for any value of $\theta\in\mathbb{R}$, see, for example, Section 2.2 of \cite{GJMN} and references therein. Knowing the decay in $N$ of such objects is one of the main ingredients to analyze both its stationary fluctuations as its hydrostatic limit. We observe that, when $\alpha \neq 1$, the model we consider has no matrix ansatz formulation available. As a consequence, there is not much information about its non-equilibrium stationary measure. Even-though it is known that the two-point stationary correlations of SEP($\alpha$) are negative (see Theorem 3.4 of \cite{duality_neg_corr}), nevertheless, its decay with $N$ is still an open problem. In this paper, we will not treat the case of the  fluctuations from the NESS since our method depends on having such bounds on correlations. From our results, we can not just simply take $t \to \infty$ to obtain the stationary fluctuations of SEP($\alpha$) because some of the estimates we use here depend on time and would blow up as $t$ goes to infinity. This is left as future work. Nonetheless, for the case $\theta = 0$ and any $\alpha \in \mathbb{N}$, since we can find explicit expressions for the two-point correlations for certain choices of the  parameters at the boundary rates (see for example in \cite{carinci2013duality} equation (6.8)), one can follow the same strategy of the proof developed  here and easily obtain the non-equilibrium stationary fluctuations of the system when $\theta=0$, we leave this to the reader.

Now we provide an outline of this article. In Section \ref{sec_2} we introduce the  SEP($\alpha$); we recall some known facts regarding its equilibrium measure  (see Section  \ref{sub_model_def}) and its hydrodynamic behavior (see Section \ref{sec_hydrodynamics}); and we  introduce  the setting for the analysis of the non-equilibrium fluctuations (see Section \ref{sub_non_eq}) and state our main results, namely, Proposition \ref{proposition_corr_decay} and Theorems \ref{th_flutuations} and \ref{th_conv_OU}. In Section \ref{proof_th_1_fluct} we provide the proof of Theorem \ref{th_flutuations}, which relies on  showing tightness and  characterizing the limit points; and we also  prove Theorem \ref{th_conv_OU} by spotting the main differences with respect to the  results known in the literature - in particular Proposition 2.5 of \cite{BGJS21} and Theorem 2.13 of \cite{BGJS21}. In Sections \ref{section_after_proof_main_theorems} and \ref{sec_results_occupation_times}  we obtain a collection of auxiliary results that we use in our proofs mainly related to estimating the two-point correlation function. In Appendix \ref{appendix_maximum_principles} we state and provide the proofs of various versions of the maximum principle.  In Appendix \ref{equation_correlations_computation} we provide some details on the Chapman-Kolmogorov equation for $\varphi^N_t$, when $\alpha \geq 2$, with the aim of facilitating the reading of the article. In Appendix \ref{remark_construction_G_N} we show two different arguments for the construction of the extension function that we use to bound $\varphi^N_t$: the first one via stochastic duality and the second one by analytic methods.Finally, Appendix \ref{sec_RL} is devoted to the proof of  a replacement lemma.

\section{The model and statement of results} \label{sec_2}

\subsection{The model: the SEP($\alpha$)}
\label{sub_model_def}
Fix $\alpha\in\mathbb N$ and for each $N \in \mathbb N$ let $\Lambda_N:=\{1,\dots, N-1\}$ be the one-dimensional, discrete interval and let $\overline{\Lambda}_N:=\Lambda_N\cup{\{0,N\}}$.  We will call $\Lambda_N$ the \textit{bulk}. We say that $x,y \in \Lambda_N$ are \textit{nearest neighbors} if $|y-x|=1$, and we denote it by $x \sim y$. We consider a Markov chain with state space $\Omega_N:=\{0,\dots, \alpha\}^{\Lambda_N}$. We call the elements of $\Omega_N$ \textit{configurations} and we denote them by $\eta = (\eta(x); x \in \Lambda_N)$. We interpret $\eta(x)$ as the number of particles at site $x \in \Lambda_N$ and we call the functions $(\eta(x); x \in \Lambda_N)$ the \textit{occupation variables}.
For each $x \in \Lambda_N$, let us denote by $\delta_x$ the configuration in $\Omega_N$ with exactly one particle, located at $x$, that is, 
\[
\delta_x(y) :=
\left\{
\begin{array}{c@{\;;\;}l}
1 & y =x,\\
0 & y \neq  x.
\end{array}
\right.
\]
For each $f: \Omega_N \to \mathbb R$, let $\mcb{L}_{\mathrm{bulk}}f = \mcb{L}_{\mathrm{bulk},N} f: \Omega_N \to \mathbb R$ be given by
\[
\begin{split} \label{generator_bulk_SEP}
\mcb{L}_{\mathrm{bulk}} f(\eta) 
		&:= \sum_{x=1}^{N-2} \eta(x)(\alpha-\eta(x+1)) \big\{ f(\eta+\delta_{x+1}-\delta_x) - f(\eta)\big\}\\
		&\quad + \sum_{x=1}^{N-2} \eta(x+1)(\alpha-\eta(x)) \big\{ f(\eta+\delta_x-\delta_{x+1}) - f(\eta)\big\}
\end{split}
\]
for every $\eta \in \Omega_N$. In this expression, we adopt the convention that $0 \cdot f(\eta + \delta_y-\delta_x) =0$ whenever $f(\eta +\delta_y -\delta_x)$ is not well defined. The linear operator $\mcb{L}_{\mathrm{bulk}}$ defined in this way is a Markov generator, which describes the \textit{bulk} dynamics.

For every $j \in\{ \ell, r\}$, let $0 < \lambda^j \leq 1$  and $\rho^j \in (0,\alpha)$ be fixed, and let $\theta \in \mathbb R$ be fixed. Define $x^\ell = 1$ and $x^r=N-1$. For $f: \Omega_N \to \mathbb R$, let $\mcb{L}_{j} f = \mcb{L}_{j,N} f: \Omega_N \to \mathbb R$ be given by
\[
\mcb{L}_{j} f(\eta) := \lambda^j\rho^j (\alpha-\eta(x^j)) \big\{ f(\eta +\delta_{x^j})- f(\eta)\big\} + \lambda^j(\alpha-\rho^j\big) \eta(x^j)\big\{ f(\eta - \delta_{x^j}) -f(\eta) \big\}
\]
for every $\eta \in \Omega_N$. The SEP($\alpha$) with \textit{slow{/fast} reservoirs} at $0$ and $N$ is the Markov chain $(\eta_t; t \geq 0)$ in $\Omega_N$ generated by the operator
\[
\mcb{L}_N := \mcb{L}_{\mathrm{bulk}} + \frac{1}{N^\theta} \big( \mcb L_\ell + \mcb L_r\big).
\]
Observe that the operator $\mcb{L}_N$ depends on the parameters $\alpha, \lambda^\ell, \lambda^r, \rho^\ell, \rho^r,\theta$. Sometimes it will be useful to state this dependence explicitly on the notation. Whenever we need to do this, we will use the generic index $i$ to denote the vector of parameters $(\alpha,\lambda^\ell, \lambda^r, \rho^\ell, \rho^r,\theta)$.

The dynamics of the  SEP($\alpha$) with parameters $(\lambda^\ell, \lambda^r, \rho^\ell, \rho^r,\theta)$ is described in the figure below.

\textit{
\begin{figure}[h]
\centering
\begin{tikzpicture}[scale=0.90]
\node at (5,0.7) (a) [fill,circle,inner sep=0.15cm,color=blue] { };
\node at (6,0.7) (b) [fill,circle,inner sep=0.15cm,color=blue] { };
\node at (6,1.3) (b') [fill,circle,inner sep=0.15cm,color=blue] { };
\node at (6,1.8) (b*) [fill,circle,inner sep=0.05cm,color=black] { };
\node at (6,2) (b'') [fill,circle,inner sep=0.05cm,color=black] { };
\node at (6,2.2) (b**) [fill,circle,inner sep=0.05cm,color=black] { };
\node at (6,2.7) (o) [fill,circle,inner sep=0.15cm,color=blue] { };
\draw (0,0) -- (14,0);
\foreach \x in {0,...,2} \draw (\x, 0.3) -- (\x, -0.3) node (\x) [above=1cm] { } node[below] {\x}; 
\node at (0,0) [below = 0.8cm] {left};
\node at (0,0) [below = 1.2cm] {reservoir};
\node at (0,0) (A) [below = 1.8cm] { };
\node at (1,0) (A') [below = 1.2cm] { };
\node at (3,-0.5) {\ldots};
\draw[-latex] (1) to[out=60,in=120] node[midway,font=\scriptsize,above] {$\hspace{0.7cm}\eta(1)(\alpha-\eta(2))$} (2);
\draw[-latex,color = red] (1) to[out=130,in=80] node[midway,font=\scriptsize,above,color=red] {$\hspace{-0.5cm}\frac{\lambda^\ell (\alpha-\rho^\ell)}{N^\theta} \eta(1)$} (0);
\draw[-latex, color = blue] (A) to [bend right] node[midway,font=\scriptsize,below=0.1cm,color=blue] {$\frac{\lambda^\ell \rho^\ell}{ N^\theta} (\alpha-\eta(1))$} (A');
\foreach \x\y in {4/x - 1,6/x + 1} \draw (\x, 0.3) -- (\x, -0.3) node (\x) [above=1.4cm] { } node [below] {\y};
\foreach \x\y in {5/x} \draw (\x, 0.3) -- (\x, -0.3) node (\x) [above=1.4cm] { } node [below=0.065cm] {\y};
\node at (7,-0.5) {\ldots};
\draw[-latex] (5) to[out=100,in=120] node[midway,font=\scriptsize,above] { } (o);
\draw (5.1,2.9) node[cross,red] { };
\draw[-latex] (5) to[out=130,in=80] node[midway,font=\scriptsize,above] {$\hspace{-1.5cm}\eta(x)(\alpha-\eta(x-1))$} (4);
\foreach \z\w in {8/y-1,10/y+1} \draw (\z, 0.3) -- (\z, -0.3) node (\z) [above=1.4cm] { } node [below] {\w};
\foreach \z\w in {9/y} \draw (\z, 0.3) -- (\z, -0.3) node (\z) [above=1.4cm] { } node [below=0.065cm] {\w};
\node at (8,0.7)  [fill,circle,inner sep=0.15cm,color=blue] { };
\node at (8,1.2) [fill,circle,inner sep=0.05cm,color=black] { };
\node at (8,1.4) [fill,circle,inner sep=0.05cm,color=black] { };
\node at (8,1.6) [fill,circle,inner sep=0.05cm,color=black] { };
\node at (8,2.1) [fill,circle,inner sep=0.15cm,color=blue] { };
\node (c) at (9,0.7) [fill,circle,inner sep=0.15cm,color=blue] { };
\node (c) at (9,1.3) [fill,circle,inner sep=0.15cm,color=blue] { };
\node (c) at (9,1.8) [fill,circle,inner sep=0.05cm,color=black] { };
\node (c) at (9,2) [fill,circle,inner sep=0.05cm,color=black] { };
\node (c) at (9,2.2) [fill,circle,inner sep=0.05cm,color=black] { };
\node (c1) at (9,2.7) [fill,circle,inner sep=0.15cm,color=blue] { };
\node at (9,2.8) (Z) [above] { };
\node at (8,2.8) (W) [above] { };
\node at (10,2.8) (Y) [above] { };
\draw[-latex] (Z) to[out=130,in=80] node[midway,font=\scriptsize,above]  { } (W);
\draw [decorate,decoration={brace,amplitude=0.1cm, mirror}] (6.5,0.4) -- (6.5,3) node [black,midway,xshift=0.3cm] {\scriptsize$\alpha$};
\node at ([shift={(40:-2)}]9.4,4.8) {\scriptsize $\hspace{-0.2cm}\eta(y)(\alpha - \eta(y-1))$};
\draw[-latex] (Z) to[out=60,in=120] node[midway,font=\scriptsize,above] {\scriptsize $\hspace{1.3cm}\eta(y)(\alpha-\eta(y+1))$} (Y);
\foreach \x\y in {12/N-2,13/N-1,14/N} \draw (\x, 0.3) -- (\x, -0.3) node (\x) [above=1cm] { } node [below] {\y};
\node at (11,-0.5) {\ldots};
\node at (13,0) (B') [below = 1.2cm] { };
\node at (14,0) (B) [below = 1.8cm] { };
\node at (14,0) [below = 0.8cm] {right};
\node at (14,0) [below = 1.2cm] {reservoir};
\draw[-latex] (13) to[out=130,in=60] node[midway,font=\scriptsize,above] {$\hspace{-1.8cm}\eta(N-1)(\alpha-\eta(N-2))$} (12);
\draw[-latex,color = red] (13) to[out=80,in=130] node[midway,font=\scriptsize,above,color=red] {$\hspace{1cm}\frac{\lambda^r(\alpha-\rho^r)}{N^\theta} \eta(N-1)$} (14);
\draw[-latex, color = blue] (B) to [bend left] node[midway,font=\scriptsize,below=0.1cm,color=blue] {$\frac{\lambda^r \rho^r}{ N^\theta} (\alpha-\eta(N-1))$} (B');
\end{tikzpicture}
\vspace*{-1cm}
\caption{Dynamics of SEP($\alpha$).}
\end{figure}
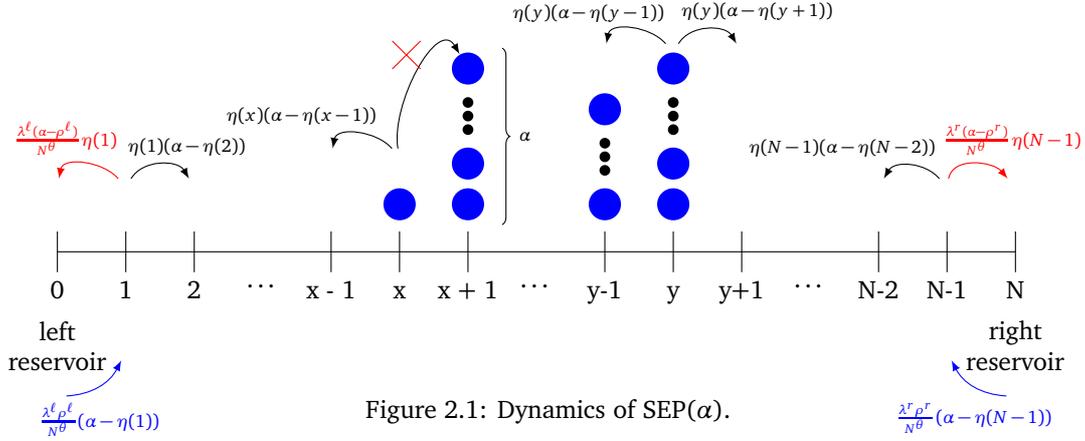}

The choice of such parametrization allows to interpret the reservoirs' dynamics in a similar way to the bulk dynamics.
More precisely, let us define 
\begin{equation}\label{eq:choiceofparameters}
\epsilon = \lambda^{\ell} \rho^{\ell}, \qquad \delta = \lambda^{r} \rho^{r}, \qquad \gamma = \lambda^{\ell} (\alpha - \rho^{\ell}), \qquad \beta = \lambda^{r} ( \alpha - \rho^{r}).
\end{equation}
Interpreting $\lambda^j \rho^j$ for $j=\ell, r$ as the corresponding particle densities at the two reservoirs,
then the jump rates of the reservoirs' dynamics corresponds to the jump rates of the bulk dynamics on which the occupation variables of sites outside the interval $\Lambda_N$ are replaced by their corresponding densities.

Hereafter we fix $T>0$ and we consider a finite time horizon $[0,T]$.  For each $N \geq 1$, we denote by $\mcb{D}_N([0,T],\Omega_N)$ the space of càdlàg trajectories endowed with the $J_1$-Skorohod topology. We fix a sequence of probability measures {$(\mu^N)_{N \geq 1}$} on $\Omega_N$. In order to see a non-trivial evolution of macroscopic quantities we need to speed up the process in the diffusive time scale $tN^2$, and in that case $\eta_{tN^2}$ has generator $N^2\mcb L_N$.  Let $\mathbb{P}_{\mu^N}$ be the probability measure on $\mcb{D}_N([0,T],\Omega_N)$ induced by the Markov process $(\eta_{tN^2};{t\geq 0})$ and by the initial measure $\mu^N$. We denote the expectation with respect to $\mathbb{P}_{\mu^N}$ by $\mathbb{E}_{\mu^N}$.
\subsection{Stationary measures}

Since the SEP($\alpha$) is an irreducible continuous time Markov chain with a finite state space, then it admits a unique  stationary measure. In fact this stationary measure can be identified for a certain choice of the parameters of the model. 
\begin{proposition}
If 
$\rho^\ell = \rho^r=:\rho$, then the stationary (equilibrium) measure is given by an homogeneous product measure with Binomial marginal distributions with parameters $\alpha \in \mathbb{N}$ and $\frac{\rho}{\alpha} \in (0,1)$:
\begin{equation} \label{reversible_measure_SEP}
\nu(\eta) = \prod_{x \in \Lambda_N} {\alpha \choose \eta(x)} \Big(\frac{\rho}{\alpha}\Big)^{\eta(x)} \Big(1-\frac{\rho}{\alpha}\Big)^{\alpha - \eta(x)}.
\end{equation}
\end{proposition}
See \cite{carinci2013duality} for a proof when $\theta = 0$, for $\theta \neq 0$ the proof is identical.

We note that for $\rho^\ell \neq \rho^r$ we do not have any information about this measure.

\subsection{Hydrodynamic limit} \label{sec_hydrodynamics}

Here we recall the hydrodynamic limit  for the SEP($\alpha$) which was obtained in \cite{FGS2022}.
For $\eta \in \Omega_N$, we define the empirical measure $\pi^{N}(\eta,du)$ by 
\begin{equation*}\label{MedEmp}
\pi^{N}(\eta, du):=\dfrac{1}{N}\sum _{x\in\Lambda_N }\eta(x)\delta_{\frac{x}{N}}\left( du\right),
 \end{equation*}
where $\delta_{b}(du)$ is a Dirac measure at $b \in [0,1]$. For every $G: [0,1] \rightarrow \mathbb{R}$ continuous, we denote the integral of $G$ with respect to $\pi^N$ by $\langle \pi^N\!, G \rangle$
and we observe that 
$$ \langle \pi^N\!, G \rangle = \frac{1}{N} \sum_{x \in \Lambda_N} \eta(x) G \left( \tfrac{x}{N}  \right).$$ 
We denote by $\mcb{M}$ the space of non-negative Radon measures on $[0,1]$ with total mass bounded by $\alpha$ and equipped with the weak topology. Also, we denote by $\mcb{D}_N([0,T],\mcb{M})$ the space of càdlàg trajectories in $\mcb{M}$ endowed with the Skorohod topology. We define $\pi^{N}_{t}(\eta, du):=\pi^{N}(\eta_{tN^2}, du)$.

\begin{definition} 
Let $\gamma: [0,1]\rightarrow [0,\alpha]$ be a measurable function. We say that a sequence  of probability measures {$(\nu^N)_{N \geq 1}$} on $\Omega_N$ is associated to the profile $\gamma$ if for every continuous function $G:[0,1] \to \mathbb{R}$ and for every $\delta >0$, it holds
\begin{equation} \label{def_seq_associated_to_prof}
\lim_{N \rightarrow \infty} \nu^N \Big( \eta \in \Omega_N: \big| \langle \pi^N\!, G \rangle - \int_0^1 G(u) \gamma (u) du \big| > \delta \Big) =0. 
\end{equation}
\end{definition}
From now on we make the following  assumption on the  sequence of probability measures:
\begin{equation}
\label{H0}
\refstepcounter{hypeq}
\tag{H\thehypeq}
\text{ {$(\mu^N)_{N \geq 1}$} is associated to a {measurable} function $\gamma: [0,1]\rightarrow [0,\alpha]$.}
\end{equation}

In order to properly state the hydrodynamic limit, i.e.~Theorem \ref{th:hyd_ssep}, we  need to recall the notion of weak solutions stated in \cite{FGS2022}.
To this end, we need to consider a proper space of test functions. We denote by $C^{1,\infty}([0, T] \times [0,1])$ the space of continuous functions defined on $[0, T] \times [0,1] $ that are continuously differentiable on the first variable and infinitely differentiable on the second variable. We also denote by $C^{1,\infty}_c ([0,T] \times [0,1])$ the space of functions $G \in C^{1,\infty}([0, T] \times[0, 1])$ such that for each time $t$, the support of $G_t$ is contained in $(0,1)$. We denote by $C^\infty([0,1])$ the space of infinitely differentiable functions defined in $[0,1]$ and we denote by $C_{c}^{m}([0,1])$ (resp.~$C_c^\infty ([0,1])$) the space of $m$-continuously differentiable (resp.~infinitely differentiable) real-valued functions defined on $[0,1]$ with support contained in $(0,1)$. We denote by $\langle \cdot, \cdot \rangle$ the inner product in ${L}^2([0,1])$ and we denote by $\|\cdot\|_{{L}^2}$ the corresponding ${L}^2$-norm. Now we define the Sobolev space $\mcb{H}^1$ on $[0,1]$. For that purpose, we define the semi inner-product $\langle \cdot, \cdot \rangle_{1}$ on the set $C^{\infty} ([0,1])$ by 
$\langle G, H \rangle_{1} := \langle\partial_u G,\partial_u H\rangle
$
for $G,H\in C^{\infty} ([0,1])$
and we denote the corresponding semi-norm by $\| \cdot \|_{1}$.
\begin{definition}
\label{Def. Sobolev space}
The Sobolev space $\mcb{H}^{1}$ on $[0,1]$ is the Hilbert space defined as the completion of $C^\infty ([0,1])$ with respect to the norm 
$\| \cdot\|_{{\mcb{H}}^1}^2 :=  \| \cdot \|_{L^2}^2  +  \| \cdot \|^2_{1}$ and  its elements coincide a.e.~with continuous functions. The space $L^{2}(0,T;\mcb{H}^{1})$ is the set of measurable functions $f:[0,T]\rightarrow  \mcb{H}^{1}$ such that 
$\int^{T}_{0} \Vert f_{t} \Vert^{2}_{\mcb{H}^{1}}dt< \infty.$
\end{definition}

We remark that in $\mcb{H}^{1}$ we can define the trace operator, and so it makes sense to talk about boundary values of functions in this space when interpreted in the trace sense.

\begin{definition}
\label{Def. Dirichlet source Condition-g_ssep}
Let $\gamma_0:[0,1] \to [0,\alpha]$ be a measurable function. We say that  $\rho:[0,T]\times[0,1] \to [0,\alpha]$ is a weak solution of the heat equation
\begin{align} \label{hydrodynamic_equation}
\begin{cases}\partial_{t}\rho_{t}(u)= {\alpha} \Delta\, {\rho} _{t}(u), \quad (t,u) \in (0,T]\times(0,1)\\
\rho_0(u) = \gamma_0(u) , \quad u \in [0,1].
\end{cases}
\end{align}
with initial condition $\gamma_0(\cdot)$ and:
\begin{enumerate}
    \item Dirichlet boundary conditions given by
 \begin{equation}
 \label{Dirichlet Equation-g_ssep}
 \rho_{t}(0)= {\rho}^\ell\quad \textrm{and}\quad\rho_{t}(1)={\rho}^r,  \quad t \in (0,T],
 \end{equation}
 if $\rho \in L^{2}(0,T;\mcb{H}^{1})$, ${ \rho} _{t}(0)={\rho}^\ell$ and ${ \rho}_{t}(1)= {\rho}^r$ for a.e.~$t \in (0,T]$, and  for all $t\in [0,T]$ and all $G \in C_c^{1,\infty} ([0,T]\times[0,1])$ it holds
\[
\langle \rho_{t},  G_{t}\rangle  -\langle \gamma_0,   G_{0} \rangle
- \int_0^t\langle \rho_{s},\Big( {\alpha}\Delta + \partial_s\Big) G_{s}\rangle ds=0.
\]

    \item Robin  boundary conditions given by
 \begin{equation}\label{Robin Equation-g_ssep}
 \partial_{u}\rho _{t}(0)=\lambda^\ell \big(\rho_{t}(0) -\rho^\ell\big),\quad \partial_{u} \rho_{t}(1)= \lambda^r \big(\rho^r - \rho_t(1)\big),\quad t \in (0,T],
 \end{equation}
if $\rho \in L^{2}(0,T;\mcb{H}^{1})$ and for all $t\in [0,T]$ and all $G \in C^{1,\infty} ([0,T]\times[0,1])$ it holds
\begin{equation*}\label{eq:Robin integral-g_ssep}
\begin{split}
\langle \rho_{t},  G_{t} \rangle -\langle \gamma_0,   G_{0} &\rangle   - \int_0^t\langle\rho_{s}, \Big({\alpha}\Delta + \partial_s\Big) G_{s} \rangle ds + {\alpha} \int^{t}_{0} \left[   \rho_{s}(1)\partial_uG_{s}(1)-\rho_{s}(0) \partial_u G_{s}(0)  \right] \, ds\\
& \qquad-\alpha \int^{t}_{0}  \left[ G_{s}(0)\lambda^\ell \big(\rho_{s}(0) -{\rho}^\ell\big) +G_{s}(1) \lambda^r \big({\rho}^r - \rho_s(1)\big)  \right] \,  ds=0.
\end{split}   
\end{equation*}
\item Neumann boundary conditions given by
\begin{equation}
\label{Neumann Equation-g_ssep}
 \partial_{u}\rho _{t}(0)= \partial_{u} \rho_{t}(1)= 0,
 \end{equation}
 if $\rho \in L^{2}(0,T;\mcb{H}^{1})$ and for all $t\in [0,T]$ and any $G \in C^{1,\infty} ([0,T]\times[0,1])$ it holds
\begin{equation*}\label{eq:Neumann integral-g_ssep}
\begin{split}
\langle \rho_{t},  G_{t} \rangle -\langle \gamma_0,   G_{0} &\rangle   - \int_0^t\langle\rho_{s}, \Big({\alpha}\Delta + \partial_s\Big) G_{s} \rangle ds + {\alpha} \int^{t}_{0} \left[   \rho_{s}(1)\partial_uG_{s}(1)-\rho_{s}(0) \partial_u G_{s}(0)  \right] \, ds = 0.
\end{split}   
\end{equation*}
\end{enumerate}
\end{definition}

We observe that there exists one and only one weak solution of the heat equation with any of the previous boundary conditions, see \cite{baldasso2017exclusion}.
We are now ready to state the hydrodynamic limit of  \cite{FGS2022}.

\begin{theorem}
\label{th:hyd_ssep}
Let $\gamma:[0,1]\rightarrow[0,\alpha]$ be a measurable function and $\lbrace\mu ^{N}\rbrace_{N\geq 1}$ a sequence of probability measures associated to $\gamma(\cdot)$, i.e.~satisfying \eqref{H0}. For any $t \in [0,T]$, any continuous function $G: [0,1] \to \mathbb R$ and any $\delta >0$, it holds
\begin{equation*}\label{limHidreform}
 \lim _{N\to\infty } \mathbb P_{\mu^{N}}\big( \eta_{\cdot} : \Big|\dfrac{1}{N}\sum_{x \in \Lambda_{N} }G\left(\tfrac{x}{N} \right)\eta_{tN^2}(x) - \langle G,\rho_{t}\rangle\Big|    > \delta \Big)= 0,
\end{equation*}
where  $\rho_{t}(\cdot)$ is the unique weak solution of the heat equation with initial condition $\gamma$ and for:

a)  $\theta<1$,  Dirichlet boundary conditions \eqref{Dirichlet Equation-g_ssep};

b)   $\theta=1$, Robin boundary conditions (\ref{Robin Equation-g_ssep}); 

c)  $\theta>1$, Neumann boundary conditions (\ref{Neumann Equation-g_ssep}).
\end{theorem}

Our focus on this article is to describe the fluctuations of the system around the hydrodynamical profile. And this is what we discuss in the next subsection.

\subsection{Non-equilibrium fluctuations}
\label{sub_non_eq}

\subsubsection{The space of test functions}
\label{set_of_test_function_invariant_remark}
As we did before stating Theorem \ref{th:hyd_ssep}, in order to show the non-equilibrium fluctuations of the SEP($\alpha$), we need to introduce a proper space of test functions. Observe that realizations of white noises are not well defined as measures, but only as distributions. Therefore, we need to introduce Schwarz-like spaces of test functions. Recall that a subscript or superscript $i$ represents dependence on the parameters $i= (\alpha,\lambda^\ell, \lambda^r, \rho^\ell,\rho^r, \theta)$ of the model.

\begin{definition} \label{def_S_set_test_functions}
We define $\mathcal{S}_i$ as the set of functions $\phi$ in $C^{\infty}([0,1])$ that satisfy, for all $k \in \mathbb{N} \cup \{0\}$,
\begin{enumerate}
    \item if $\theta < 0$: $\partial^{k}_u \phi (0) = \partial^{k}_u \phi (1) = 0$;

    \item if $0 \leq \theta < 1$: $\partial^{2k}_u \phi (0) = \partial^{2k}_u \phi (1) = 0$;
    
    \item if $\theta = 1$: $\partial^{2k+1}_u \phi (0) = \lambda^\ell \partial^{2k}_u \phi (0), \ \partial^{2k+1}_u \phi (1) = - \lambda^r \partial^{2k}_u \phi(1)$;
    
    \item if $\theta > 1$: $\partial^{2k+1}_u \phi (0) = \partial^{2k+1}_u \phi (1) = 0$.
\end{enumerate}
\end{definition}

As in \cite{FGN19, GJMN}, the previous choice is to  make $\mathcal{S}_i$ invariant under taking second derivatives, which in turn implies  that the Markov semigroup associated to the operator $\alpha \Delta$ with the corresponding boundary conditions, which we denote by $S^{i}_{t}$, is such that, if $\phi \in \mathcal{S}_i$, then $S^{i}_{t}\phi \in \mathcal{S}_i$. This property will be useful later on.  Indeed, as in the proof of Proposition 3.1 of \cite{FGN19}, for the case $\theta = 1$, and for the other values of $\theta$ as in Remark 2.5. of \cite{GJMN}, given $\phi \in \mathcal{S}_i$, $S^{i}_{t}\phi$ is solution to
\begin{align*}
\begin{cases}\partial_{t}S^{i}_{t}\phi(u)= {\alpha} \Delta\, S^{i}_{t}\phi(u), \quad (t,u) \in [0,T]\times(0,1)\\
S^i_0 \phi(u) = \phi(u) , \quad u \in [0,1].
\end{cases}
\end{align*} with boundary conditions:
\begin{enumerate}
\item if $\theta > 1$
\begin{equation} \label{semigroup_eq_theta_bigger_1}
\partial_u S^{i}_{t}\phi(0) = \partial_u S^{i}_{t}\phi(1) = 0;
\end{equation}

\item if $\theta =1$
\begin{equation} \label{semigroup_eq_theta_equal_1}
\partial_u S^{i}_{t}\phi(0) = \lambda^\ell S^{i}_{t}\phi(0)  \quad \textrm{ and } \quad \partial_u S^{i}_{t}\phi(1) = - \lambda^r S^{i}_{t}\phi(1);
\end{equation}

\item if $\theta < 1$
\begin{equation} \label{semigroup_eq_theta_less_1}
S^{i}_{t}\phi(0) = S^{i}_{t}\phi(1) = 0.
\end{equation}
\end{enumerate}

Let us compute $S_t^i$ by the separation of variables method. The aim is to look for solutions of the form
\begin{equation} \label{formula_semigroup_to_susbtitute}
S^{i}_{t}\phi (u) = g(t) f(u),
\end{equation} 
with $g$ a function of $t$ and $f$ a function of $x$ to be computed. This leaves us with $g(t) = C e^{\mu \alpha t}$, where $C, \mu \in \mathbb{R}$ to be computed, and the Sturm-Liouville problem
$
f''(u) - c f(u) = 0,$ for $u \in (0,1)
$ with boundary conditions
\begin{enumerate}
\item if $\theta > 1$, $f'(0) = f'(1) = 0$;

\item if $\theta = 1$, $f'(0) = \lambda^\ell f(0)$ and $f'(1) = -\lambda^r f(1)$;

\item if $\theta < 1$, $f(0) = f(1) = 0$.
\end{enumerate}

The previous problems have a solution of the form
$f(u) = A \sin(\omega_1 u) + B \cos(\omega_2 u)$, where $A, B, \omega_1, \omega_2$ have to be computed. A simple but long computation shows that 
\begin{enumerate}
\item if $\theta > 1$,
$
f(u) = B(k) \cos(\pi k u), \textrm{ for some } k \in \mathbb{Z},
$ where $B(k)$ has to be computed.
Thus,
\begin{equation} \label{semigroup_formula_theta_1}
S^i_t \phi (u) = \sum_{k \in \mathbb{Z}} e^{-\pi^2 k^2 \alpha t} \langle \phi, 2 \cos (\pi k \cdot)\rangle \cos(\pi k u).
\end{equation}

\item if $\theta = 1$,
$
f(u) = B(k)\left[ \frac{\lambda^\ell}{\beta_k} \sin(\beta_k u) + \cos(\beta_k u)\right], \textrm{ for some } k \in \mathbb{Z},
$ where $B(k)$ has to be computed and $\beta_k$ are the solutions of
$
\frac{(\lambda^\ell + \lambda^r)x}{x^2 + \lambda^\ell \lambda^r} = \tan(x).
$
Thus,
\begin{equation} \label{semigroup_formula_theta_lower_1}
S^i_t \phi (u) = \sum_{k \in \mathbb{Z}} e^{- \beta_k^2 \alpha t} B(k)\left[ \frac{\lambda^\ell}{\beta_k} \sin(\beta_k u) + \cos(\beta_k u)\right],
\end{equation} with $B(k)$ such that $\sum_{k \in \mathbb{Z}} B(k)\left[ \frac{\lambda^\ell}{\beta_k} \sin(\beta_k u) + \cos(\beta_k u)\right] = \phi(u)$.

\item if $\theta < 1$,
$
f(u) = A(k) \sin(\pi k u), \textrm{ for some } k \in \mathbb{Z},
$ where $A(k)$ has to be computed.
Thus,
\begin{equation} \label{semigroup_formula_theta_bigger_1}
S^i_t \phi (u) = \sum_{k \in \mathbb{Z}} e^{-\pi^2 k^2 \alpha t} \langle \phi, 2 \sin (\pi k \cdot)\rangle \sin(\pi k u).
\end{equation}
\end{enumerate}

For every $\theta \in \mathbb{R}$, we showed that $S^i_t \phi$ can be written in terms of the eigenvalues and eigenfunctions of the Laplace operator with different boundary conditions. From here we easily conclude that, for every $\phi \in \mathcal{S}_i$, $S^i_t \phi \in \mathcal{S}_i$.


We equip $\mathcal{S}_i$ with the topology induced by the family of seminorms $\left\{ |||\cdot|||_{j} \right\}_{j \in \mathbb{N} \cup \{0\}}$ where for $\phi\in\mathcal S_i$ 
\begin{align} \label{seq_seminorms_def}
|||\phi|||_{j} :=  \sup_{u \in [0,1]} | \phi^{(j)} (u)|.
\end{align}
The space $\mathcal S_i$ endowed with this topology turns out to be a nuclear Fréchet space, i.e.~a complete Hausdorff space whose topology is induced
by a countable family of semi-norms and such that all summable sequences in $\mathcal S_i$ are absolutely summable. We will denote by $\mathcal{S}'_i$ the topological dual of $\mathcal{S}_i$, i.e.~the set of linear bounded functionals over $\mathcal{S}_i$ and we  equip it with the weak topology. Let $\mcb{D}_N([0,T],\mathcal{S}'_i)$ denote the set of c\`adl\`ag time trajectories of linear functionals acting on $\mathcal{S}_i$.

\subsubsection{The discrete profile and the density fluctuation field}

Observe that Theorem \ref{th:hyd_ssep} can be understood as a law of large numbers for the random trajectories $(\langle\pi^N_t\!,G\rangle; t \geq 0)$. Therefore, it is natural to study the corresponding central limit theorem. In order to do that, one needs to specify how to center and how to rescale the random variables $\langle \pi_t^N\!, G\rangle$. Whenever possible, the most natural way to do this is to consider the quantity
\[
\sqrt{N} \big( \langle \pi_t^N\!, G\rangle - \mathbb E_{\mu^N} [ \langle \pi_t^N\!, G\rangle ] \big).
\]
Thanks to the duality properties of the SEP($\alpha$), the expectation $\mathbb E_{\mu^N} [ \langle \pi_t^N\!, G\rangle ]$ can be computed in a fairly explicit way.
Let us define the \textit{expected density of particles} $\rho_t^N(x)$ for all $t \geq 0$ and $x \in \overline{\Lambda}_N$ as
\begin{equation*}
\rho_t^N(x) := \mathbb E_{\mu^N} [ \eta_{tN^2}(x)] \text{ for } x \in \Lambda_N \text{ and } \rho_t^N(0) := \rho^\ell, \quad \rho_t^N(N) := \rho^r \;.
\end{equation*}

This last definition serves as a boundary condition for the expected density of particles. Using that the monomials $\left( \frac{\eta_x}{\alpha}; x \in \Lambda_N \right)$ are self-duality functions for the SEP($\alpha$), one can show that $(\rho_t^N(x); t \geq 0, x \in \overline{\Lambda}_N)$ is the unique solution of the discrete heat equation
\begin{align} \label{system_of_eq_for_discrete_rho}
\left\{
\begin{array}{r@{\;=\;}l}
\partial_t \rho_t^N(x) & N^2 \Delta^i_N \rho_t^N(x), x \in \Lambda_N, t \geq 0,\\
\rho_t^N(0) & \rho^\ell, t \geq 0,\\
\rho_t^N(N) & \rho^r, t \geq 0,
\end{array}
\right.
\end{align}
 with initial condition $\rho_0^N(x) := \mathbb E_{\mu^N}[\eta_0^N(x)]$. Here the operator $\Delta^i_N$ is a discrete Laplacian with modified rates at the boundary depending on $i$. More precisely, let us define the jump rate 
 \[
 c^i\!: \{ (x,y) \in \Lambda_N \times \overline{\Lambda}_N ; x \sim y\}
\]
as
\begin{equation} \label{ratesss}
c^i_{x,y} \!:=
\left\{
\begin{array}{c@{\;;\;}l}
\alpha & x,y \in \Lambda_N\\
\displaystyle{\frac{\alpha \lambda^\ell}{N^\theta}} & x=1, y=0\\
\displaystyle{\frac{\alpha \lambda^r}{N^\theta}} & x=N-1, y=N.\\
\end{array}
\right.
\end{equation}
Then the operator $\Delta^i_N$ acts on functions $f: \overline{\Lambda}_N \to \mathbb R$ as
\begin{equation} \label{laplaciannn}
\Delta^i_N f(x) = c_{x,x-1}^i \big( f(x-1)-f(x) \big) + c^i_{x,x+1} \big( f(x+1)-f(x) \big),
\end{equation}
for every $x \in \Lambda_N$.

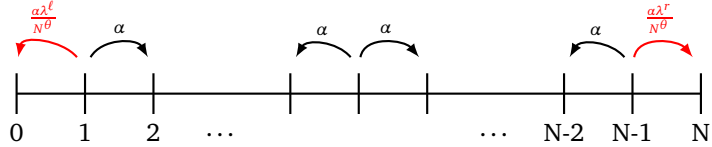
\begin{figure}[h!]
    \centering
\begin{tikzpicture}[thick, scale=0.9]
\draw (0,0) -- (10,0);
\foreach \x in {0,...,2} \draw (\x, 0.3) -- (\x, -0.3) node (\x) [above=0.5cm] { } node[below] {\x}; 
\node at (0,0) (A) [below = 1.cm] { };
\node at (1,0) (A') [below = 0.4cm] { };
\node at (3,-0.6) {\ldots};
\draw[-latex] (1) to[out=60,in=120] node[midway,font=\scriptsize,above] {$\alpha$} (2);
\draw[-latex,color = red] (1) to[out=130,in=80] node[midway,font=\scriptsize,above,color=red] {$\frac{\alpha \lambda^\ell}{N^\theta}$} (0);
\foreach \x in {4,5,6} \draw (\x, 0.3) -- (\x, -0.3) node (\x) [above=0.5cm] { } node [below] { };
\node at (7,-0.6) {\ldots};
\draw[-latex] (5) to[out=80,in=130] node[midway,font=\scriptsize,above] {$\alpha$} (6);
\draw[-latex] (5) to [out=130,in=60] node[midway,font=\scriptsize,above] {$\alpha$} (4);
\foreach \x\y in {8/N-2,9/N-1,10/N} \draw (\x, 0.3) -- (\x, -0.3) node (\x) [above=0.5cm] { } node [below] {\y};
\node at (8,0) (B') [below = 1.2cm] { };
\node at (9,0) (B) [below = 1.8cm] { };
\draw[-latex] (9) to[out=130,in=60] node[midway,font=\scriptsize,above] {$\alpha$} (8);
\draw[-latex,color = red] (9) to[out=80,in=130] node[midway,font=\scriptsize,above,color=red] {$\frac{\alpha \lambda^r}{N^\theta}$} (10);
\end{tikzpicture}
\vspace*{-1cm}
\caption{Illustration through arrows of the jump rate $c^i$ defined above.}
\end{figure}

The stationary solution of \eqref{system_of_eq_for_discrete_rho}, that we denote by $\rho_{ss}^N(\cdot)$, is given, for every $x \in \Lambda_N$ by
\begin{equation}
\rho_{ss}^N(x) := a^i_N x + b^i_N,
\end{equation} where
\begin{equation} \label{a_Nplusb_N}
a^i_N = \frac{\lambda^\ell}{N^\theta - \lambda^\ell}(b_N - \rho^\ell) \quad \textrm{ and } \quad b^i_N = \frac{\lambda^r \rho^r(N^\theta - \lambda^\ell) + \lambda^\ell \rho^\ell(N^\theta + (N-1)\lambda^r)}{\lambda^\ell \lambda^r(N-1) + \lambda^\ell N^\theta + \lambda^r(N^\theta - \lambda^\ell)}.
\end{equation}
\begin{definition}
We define the density fluctuation field
$(Y_t^N; t \geq 0)$ associated to the SEP($\alpha$), $(\eta_{tN^2}; t \geq 0)$, with initial measure $(\mu^N)_{ N \in \mathbb N}$ as the time trajectory of linear functionals acting on functions $\phi \in \mathcal{S}_i$ as
\begin{equation} \label{fluctuation_field_def}
    Y_t^N(\phi) = \frac{1}{\sqrt{N}} \sum_{x\in \Lambda_N} \phi\left(\frac{x}{N}\right) \bar{\eta}_{tN^2}(x),
\end{equation}
where, for each $x \in \Lambda_N$,  we centered $\eta_{tN^2}(x)$ by taking $\bar{\eta}_{tN^2}(x):= \eta_{tN^2}(x) - \rho^N_t(x)$.
\end{definition}

For each $N \in \mathbb N$, let $\mathbb{Q}_N$ be the probability measure in $\mcb{D}_N([0,T],\mathcal{S}'_i)$, induced by the density fluctuation field $(Y_t^N)_{t \geq 0}$. Our goal is to prove, under suitable assumptions, that $(\mathbb{Q}_N)_{N \in \mathbb{N}}$ weakly converges to $\mathbb{Q}$, a probability measure on $\mcb{D}_N([0,T],\mathcal{S}'_i)$, that can be uniquely characterized. A limit theorem of this form is known in the literature as the derivation of the \textit{non-equilibrium fluctuations} of the SEP($\alpha$). To achieve our goal, it will be enough to: show that the sequence of measures $(\mathbb{Q}_N)_{N \in \mathbb{N}}$ is tight, guaranteeing the weak convergence up to a subsequence and then characterize (uniquely) the limit point. Roughly speaking, this is the content of Theorems \ref{th_flutuations} and \ref{th_conv_OU}. 


\subsubsection{Main results}

To properly state our results, we need to introduce some definitions and notations. A crucial estimate for the non-equilibrium fluctuations is a sharp estimate on the decay of both space and space-time correlation function of the SEP($\alpha$).
Define the two-dimensional set  $V_N := \{(x,y) \in (\Lambda_N)^2 \ | \ x \leq y\}$ and its boundary  by
\[
\partial {V}_N:= \{ (x,y)\,: \, x\in\{0,N\} \,\,\textrm {and} \,\, y\in\Bar\Lambda_N\}\cup \{ (x,y)\,: \, y\in\{0,N\} \,\,\textrm {and} \,\,x\in\Bar\Lambda_N\,\}.
\]
We denote its closure by $\overline{V}_N := V_N\cup \partial V_N$, and we denote  its upper diagonal and its diagonal, respectively,   by
\begin{align} \label{def_D_N_and_main_diag}
    \mathcal{D}_N^+ := \{(x,y)\in V_N \ | \ y=x + 1\} \textrm{ and } \mathcal{D}_N := \{(x,y)\in  V_N \ | \ y = x\}.
\end{align}

\begin{figure}[h!]
    \centering
\begin{tikzpicture}[thick, scale=0.9]
\draw[->] (-0.5,0)--(5.5,0) node[anchor=north]{$x$};
\draw[->] (0,-0.5)--(0,5.5) node[anchor=east]{$y$};
\begin{scope}[scale=0.75]
\foreach \x in {1,...,3} 
	\foreach \y in {\x,...,3,4}
		\shade[ball color=black](\x,1+\y) circle (0.15);
  \foreach \x in {1,...,3,4} 
	\foreach \y in {\x}
		\shade[ball color=blue](\x,1+\y) circle (0.15);
\foreach \x in {1,...,5} 
	\shade[ball color = red](\x,6) circle (0.15); 
\foreach \x in {1,...,5} 
	\shade[ball color = gray!50](\x,\x) circle (0.15); 
\foreach \y in {1,...,5} 
	\shade[ball color = red](0,\y) circle (0.15); 
\shade[ball color= red](6,6) circle (0.15);
\shade[ball color= red](0,6) circle (0.15);
\shade[ball color= red](0,0) circle (0.15);
\end{scope}	
\draw (0,0) node[anchor=north east] {$0$};
\draw (0.75,2pt)--(0.75,-2pt) node[anchor=north] {$1$};
\draw (1.5,2pt)--(1.5,-2pt) node[anchor=north]{$2$};
\draw (2.25,2pt)--(2.25,-2pt);
\draw (3,2pt)--(3,-2pt);
\draw (3.75,2pt)--(3.75,-2pt) node[anchor=north]{\small $N-1$};
\draw (4.5,2pt)--(4.5,-2pt) node[anchor=north]{\small $N$};
\draw (-0.05,.75) node[anchor=east] {$1$};
\draw (-0.05,1.5) node[anchor=east]{$2$};
\draw (-0.05,3.75) node[anchor= east]{$N-1$};
\draw (-0.05,4.5) node[anchor=east]{$N$};
\node at (7,4.25) {\color{black} \underline{Color interpretation:}};
\node at (7,2.75) {\color{red} $\partial V_N$};
\node at (7,2) {\color{blue} $\mathcal{D}_N^+$};
\node at (7,1.25) {\color{gray} $\mathcal{D}_N$};
\end{tikzpicture}
\caption{Illustration of the sets  $\partial V_N$ (in red), $\mathcal{D}_N^+$ (in blue) and $\mathcal{D}_N$ (in gray).}
    \label{fig:triangle}
\end{figure}
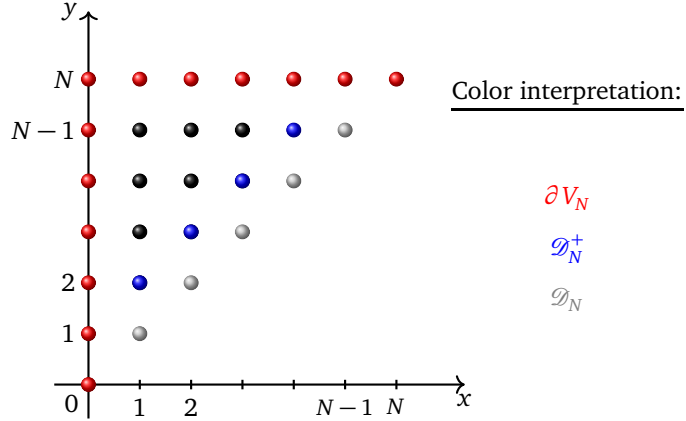

\begin{definition}
Let  $(\varphi_t^N; t \geq 0)$ be the time-dependent, two-point correlation function, defined on $(x,y) \in V_N$ with $x \neq y$ by
\begin{align} \label{time_dependent_correlation_def}
 \varphi^N_t(x,y)  &:= 
\begin{cases} 
\mathbb{E}_{\mu^N} [\bar{\eta}_{tN^2}(x)\bar{\eta}_{tN^2}(y)], \textrm{ if }  (x,y) \notin \partial V_N,\\
0, \textrm{ if } (x,y) \in \partial V_N,
\end{cases}
\end{align} and extended symmetrically to $(\overline{\Lambda}_N)^2 \setminus \mathcal{D}_N$.
\end{definition}

Now we make some extra assumptions on the initial measures, besides \eqref{H0}. We assume that there exists a continuous profile $\gamma:[0,1] \to [0,\alpha]$  such that
\begin{equation}\label{H1_3}
\refstepcounter{hypeq}
\tag{H\thehypeq}
\frac{1}{N} \sum_{x=1}^N \Big|\rho_0^N(x) - \gamma \left(\tfrac{x}{N}\right)\Big| \xrightarrow[\text{ }]{N \to \infty} 0.
\end{equation}
We also assume that there exists a sequence of profiles
$g_N(\cdot)$ of class $C^6$ that satisfy, for each $N \geq 1$ 
\begin{equation}\label{H1_6}
\refstepcounter{hypeq}
\tag{H\thehypeq}
\partial_u^j g_N = \partial_u( N a_N^i u + b_N^i),
\end{equation}
for $u \in [0,1]$ and $j=0,1,2,3$,  where $a_N^i$
and $b_N^i$ were defined in \eqref{a_Nplusb_N} and such that, for every $N \geq 1$,
\begin{equation}
\label{H2}
\refstepcounter{hypeq}
\tag{H\thehypeq}
\max_{x \in \Lambda_N} \big| \rho_0^N(x) - g_N \left(\tfrac{x}{N}\right) \big| \lsim \frac{1}{N}.
\end{equation}

We also assume that
\begin{equation}
\refstepcounter{hypeq}
\tag{H\thehypeq}
\label{decay_corr_time0_bulk}
    \max_{\substack{(x,y) \in V_N\\x\neq y}} |\varphi^N_0(x,y)| \lesssim 
    \frac{1}{N}, \quad \max_{x \in \Lambda_N \setminus \{1,N-1\}} \big|\mathbb E_{\mu^N} \big[ \alpha \eta_0(x)(\eta_0(x)-1) -(\alpha-1) \rho_0^N(x)^2\big] \big| \lesssim \frac{1}{N},
\end{equation}
\\
and that  for $x = 1$ and $x= N-1$, 
\begin{equation}
\refstepcounter{hypeq}
\tag{H\thehypeq}
\label{decay_corr_time0_bound}
  \max_{\substack{y \in \Lambda_N \\ x \neq y}} |\varphi^N_0(x,y)| \lsim 
\frac{1}{N} \min \{ 1, N^{\theta-1}\},
\;\; \max_{x =1,N-1} \big|\mathbb E_{\mu^N} \big[ \alpha \eta_0(x)(\eta_0(x)-1) -(\alpha-1) \rho_0^N(x)^2\big] \big| 
		\lesssim \frac{1}{N} \min \{ 1, N^{\theta-1}\}.
\end{equation}

\noindent
\textbf{\underline{Notation:}} Above and in {what follows}, we denote by $\lesssim$ an inequality that is correct up to a multiplicative constant independent of $N$.

Now we present the main results of this article.
\begin{theorem}[Non-Equilibrium Fluctuations] \label{th_flutuations} Let $\alpha \geq 1$ and $\theta \in \mathbb{R}$. Let $\gamma \in C^6([0,1])$ and  $(\mu^N)_{N \in \mathbb{N}}$ a sequence of probability measures satisfying {\eqref{H0} -\eqref{decay_corr_time0_bound}}. 
Then, the sequence of probability measures $\{\mathbb{Q}_N\}_{N \in \mathbb{N}}$ is tight with respect to the $J_1$-Skorohod topology of $\mcb{D}_N([0,T],\mathcal{S}'_i)$ and all limit points $\mathbb{Q}$ are probability measures concentrated on paths $Y$ satisfying
\begin{equation}\label{eq:charac_lim_point}
    Y_t(f) = Y_0(S^{i}_t f) + W^i_t(f),
\end{equation}
for any $f \in \mathcal{S}_i$ and any $t\in[0,T]$. Above $S^{i}_t: \mathcal{S}_i \to \mathcal{S}_i$ is the semigroup associated to the hydrodynamic equation \eqref{hydrodynamic_equation} with the respective boundary conditions, and $W_t^i$ is a mean-zero Gaussian random variable of variance
\begin{equation*}
    \int_0^t \|S^{i}_{t-s}f\|^2_{L^2(\rho_s)}ds,
\end{equation*}
where, for every $s\in[0,T]$ and $g, h \in L^2(\rho_s)$,
\begin{align*}
    \langle h, g \rangle_{L^2(\rho_s)} & := \int_0^1 2 \chi_\alpha(\rho_s(u)) h(u) g(u) du  \\
    & + \mathbb{1}({\theta = 1}) \left\lbrace \left[ \lambda^\ell(1-2\rho^\ell) \rho_{s}(0) + \lambda^\ell \rho^\ell \alpha \right]  h(0)g(0) + \left[ \lambda^r(1-2\rho^r) \rho_{s}(1) + \lambda^r \rho^r \alpha \right] h(1)g(1) \right\rbrace 
\end{align*} 
and $\rho_\cdot$ is the unique weak solution of the corresponding hydrodynamic equation \eqref{hydrodynamic_equation}. Above,  \begin{equation}\label{eq:chi_rho}
\chi_\alpha(\rho) = \rho(\alpha - \rho)
\end{equation} represents the mobility of our model. Moreover, $Y_0$ and $W_t^i$ are uncorrelated in the sense that for all $f, g \in \mathcal{S}_i$ it holds $\mathbb{E}_{\mathbb{Q}}[Y_0(f) W_t^i(g)] = 0$ .
\end{theorem}

In the last theorem, we do not guarantee the convergence of the hole sequence $\{\mathbb{Q}_N\}_{N \in \mathbb{N}}$ but only the convergence up to a subsequence, whose limit points we are not able to prove their uniqueness (i.e.~independence with respect to the convergent subsequence) with only the assumptions of the theorem. Nevertheless, when we also impose the convergence at the initial time $t=0$ of $Y^N_t$ to a Gaussian process, then, uniqueness holds and we prove that the hole sequence $\{\mathbb{Q}_N\}_{N \in \mathbb{N}}$ converges to a measure $\mathbb{Q}$ which is concentrated on the unique solution of the next martingale problem, which is an Ornstein-Uhlenbeck (O.U.) process. With this extra assumption at time $t=0$, we prove uniqueness of the limit point and convergence of $\{\mathbb{Q}_N\}_{N \in \mathbb{N}}$ follows.

\begin{definition}[Ornstein-Uhlenbeck - Definition 2.4 of \cite{BGJS21}] \label{def_ou_martingale_problem}
Fix some time horizon $T > 0$. Let $C$ be a topological vector space, $A : C \to C$ an operator letting $C$ invariant and $c : C \to [0, \infty)$  a continuous functional satisfying
$
    c(\lambda H) = |\lambda| c(H),
$
for all $\lambda \in \mathbb{R}$ and $H \in C$. Let $C'$ be the topological dual of $C$ equipped with the weak-$*$ topology. Denote by $\mcb C([0, T], C')$ the set of continuous trajectories in $[0,T]$ of functionals in $C'$. We say that the process $\{Y_t ; t \in [0, T]\} \in \mcb C([0, T], C')\}$ is a solution of the O.U. martingale problem $OU(C,A,c)$ on the time interval $[0, T]$ with initial (random) condition $y_0 \in C'$ if:
\begin{enumerate}
    \item for any $H \in C$ the two real-valued processes $M_\cdot(H)$ and $N_\cdot(H)$ defined by
\begin{align*}
    &M_t(H) = Y_t(H) - Y_0(H) - \int_0^t Y_s(AH)ds,\\
    &N_t(H) = ( M_t(H))^2 - tc^2(H),
\end{align*}
    are martingales with respect to the natural filtration of the process, that
is, $\{\mathcal{F}_t \ ; \ t \in [0,T]\} = \{ \sigma( Y_s(H) \ | \ s \leq t, H \in C) \ ; \ t \in [0,T]\}.$

    \item $Y_0 = y_0$ in law.
\end{enumerate}
\end{definition}

\begin{theorem}[Convergence to the Ornstein-Uhlenbeck process] \label{th_conv_OU}
Let $\alpha \in \mathbb{N}$ and $\theta \in \mathbb{R}$. Assume the conditions of Theorem \ref{th_flutuations} and also  that the sequence of initial density fluctuation field $\{Y^N_0\}_{N \in \mathbb{N}}$ converges, as $N \to +\infty$, to a mean-zero Gaussian field
$Y_0$ with covariance given, for $f,g \in \mathcal{S}_i$, by
\begin{align*}
    \sigma(f,g) := \mathbb{E}[Y_0(f)Y_0(g)] = \lim_{N \to +\infty} \mathbb{E}[Y^N_0(f)Y^N_0(g)].
\end{align*}
Then:
\begin{enumerate}
    \item if $\theta \geq 0$, the sequence $\{\mathbb{Q}_N\}_{N \in \mathbb{N}}$ converges, as $N \to +\infty$, to a measure $\mathbb{Q}$ which is concentrated on the unique solution $Y_t$ of the O.U. martingale problem $OU(\mathcal{S}_i,\alpha\Delta, ||\cdot||_{L^2(\rho_s)})$ on the time interval $[0, T]$ with the initial (random) condition equal to $Y_0$. Thus, $Y_t$ is a generalized O.U. process, which is the unique (in law) formal solution of the stochastic partial differential equation:
\begin{align} \label{OU_equation}
    \partial_t Y_t = \alpha\Delta Y_t dt + \sqrt{2 \chi_\alpha(\rho_t)} \nabla dW_t,
\end{align}
where $dW_t$ is a space-time white noise with unit variance and $\alpha\Delta$ is the same operator as in \eqref{hydrodynamic_equation} with the corresponding boundary conditions depending on the value of $\theta$. As a consequence, the covariance of the limit field $Y_t$ is given, for $f,g \in \mathcal{S}_i$, by
\begin{align*}
\mathbb{E}[Y_t(f)Y_s(g)] = \sigma(S^i_t f, S^i_s g) + \int_0^s \langle \partial_u S^i_{t-r} f, \partial_u S^i_{s-r} g \rangle_{L^2(\rho_r)} dr,
\end{align*} with $\partial_u h (0)$ (respectively, $\partial_u h (1)$) identified with $ \displaystyle\partial_u h (0^+) = \lim_{x \downarrow 0} \partial_u h (x)$ (respectively, $\displaystyle \partial_u h (1^-) = \lim_{x \uparrow 1} \partial_u h (x)$), for $h \in \mathcal{S}_i$.

    \item if $\theta < 0$, the sequence $\{\mathbb{Q}_N\}_{N \in \mathbb{N}}$ converges, as $N \to +\infty$, to a measure $\mathbb{Q}$ which is concentrated on the unique solution $Y_t$ of the Ornstein-Uhlenbeck martingale problem $OU(\mathcal{S}_i,\alpha\Delta, ||\cdot||_{L^2(\rho_s)})$ on the time interval $[0, T]$ with initial (random) condition equal $Y_0$, and whose uniqueness (in law) of solution is guaranteed when remarking that $Y_t$ satisfies the following two extra conditions:
\begin{enumerate}
    \item regularity condition: $\mathbb{E}[(Y_t(H))^2] \lsim ||H||_{L^2}$, for any $H \in \mathcal{S}_i$;

    \item boundary condition: For each $j \in \{0,1\}$, let $\iota^j_\epsilon$ be defined as, for $j=0$, $\iota^0_\epsilon(u):=\epsilon^{-1}\mathbb {1}_{(0,\epsilon]}(u)$ and, for $j=1$, $\iota^1_\epsilon(u):=\epsilon^{-1}\mathbb {1}_{[1-\epsilon,1)}(u)$ $u \in[0,1]$.
For any $t \in [0, T]$ and $j \in \{0, 1\}$, it holds that
    \begin{equation*}
        \lim_{\epsilon \to 0} \mathbb{E}\left[ \left(\int_0^t Y_s(\iota^j_\epsilon) ds \right)^2  \right] = 0.
    \end{equation*}
\end{enumerate}
\end{enumerate}
\end{theorem}

As a consequence of the previous result we obtain the
non-equilibrium fluctuations starting from a local Gibbs state. 

\begin{corollary}
Fix a measurable profile $\gamma_0:[0,1] \to [0,\alpha]$ satisfying \eqref{H1_6} and \eqref{H2}; and start the process SEP($\alpha$) from the Binomial product
measure with marginals given by $$\nu_{\gamma_0}^N \{\eta \ | \ \eta(x) = k\} = {\alpha \choose k} \left[\tfrac{\gamma_0 \left(\tfrac{x}{N} \right)}{\alpha} \right]^k \left[ 1 - \tfrac{\gamma_0 \left(\tfrac{x}{N} \right)}{\alpha}\right]^{\alpha - k},$$ for $k \in \{0,\dots,\alpha\}$. Let $f,g \in \mathcal{S}_i$. Then Theorem \ref{th_conv_OU} holds with
\begin{align*}
\sigma(S^i_t f, S^i_s g) = \int_0^1 \chi_\alpha(\gamma_0(u)) S^i_t f(u) S^i_s g(u) du.
\end{align*}
\end{corollary}

Observe that the remaining assumptions of
Theorems \ref{th_flutuations} and \ref{th_conv_OU} are satisfied by the starting measure $\nu_{\gamma_0}^N$, so that above, we only need to impose \eqref{H1_6} and \eqref{H2} from the initial profile.

To prove Theorem \ref{th_flutuations} and Theorem \ref{th_conv_OU}, we will need some auxiliary results that we will leave their proofs and details for after showing each of the above theorems.

\section{Proof of Theorems \ref{th_flutuations} and  \ref{th_conv_OU}} \label{proof_th_1_fluct}

The proof of both theorems follows by showing first the tightness of the sequence of probability measures $\{\mathbb{Q}_N\}_{N \in \mathbb{N}}$ with respect to the Skorohod topology of $\mcb{D}_N([0,T],\mathcal{S}'_i)$; and to show that all limit points $\mathbb{Q}$ are probability measures concentrated on paths $Y$ satisfying \eqref{eq:charac_lim_point}. 
We start now with the former.

\subsection{Tightness} \label{tightness_proof}

Recall that  the spaces $\mathcal{S}_i$ are nuclear Fréchet spaces when endowed with the seminorms defined in  \eqref{seq_seminorms_def}. Therefore, in order to prove tightness, we can use Mitoma’s criterium (that we recall below) and restrict ourselves to showing tightness of the sequence of real-valued processes $\{Y^N_t(\phi)\}_{N \in \mathbb{N}}$, for every $\phi \in \mathcal{S}_i$.

\begin{theorem}[Mitoma’s criterium - Theorem 4.1 of \cite{Mitoma}]
A sequence of processes $\{X^N_t; t \in [0,T]\}_{N \in \mathbb{N}}$ in $\mcb{D}([0,T], \mathcal{S}'_i)$ is tight with respect to the Skorohod topology if, and only if, for every $H \in \mathcal{S}_i$, the sequence of real-valued processes $\{X^N_t(H); t \in [0,T]\}_{N \in \mathbb{N}}$ is tight with respect to the Skorohod topology of $\mcb{D}([0,T], \mathbb{R})$.
\end{theorem}

Recall that, from Lemma 5.1 of Appendix 1 of \cite{KL}, 
\begin{align} \label{def_dynkin_martingale}
    M^N_t(\phi) := Y^N_t(\phi) - Y^N_0(\phi) -  \int_0^t (N^2 \mcb{L}_N + \partial_s)Y^N_s(\phi) ds, 
\end{align}
is a martingale for every $\phi \in \mathcal{S}_i$.
Therefore, 
in order to show that $\{Y^N_t(\phi_t)\}_{N \in \mathbb{N}}$ is tight, it is enough to show that  \begin{equation*}
\{Y^N_0(\phi)\}_{N \in \mathbb{N}}
 	\ , \ \{[M^N_t(\phi)]_{t \geq 0}\}_{N \in \mathbb{N}} 
		 \textrm{ and } \left\{\int_0^t (N^2 \mcb{L}_N+ \partial_s) Y^N_s(\phi) ds \right\}_{N \in \mathbb{N}}
\end{equation*} are tight. We start by showing that $\{Y^N_0(\phi)\}_{N \in \mathbb{N}}$ is tight.

\subsubsection{Initial time}
By Helly-Bray theorem, it is enough to show that
\begin{equation*}
    \lim_{A \to \infty} \limsup_{N \to +\infty} \mathbb{P}_{\mu^N} [|Y_0^N(\phi)| > A] = 0.
\end{equation*}
By Markov's inequality, for every $A > 0$ and for every $N \in \mathbb{N}$, 
\begin{align*}
\mathbb{P}_{\mu^N} [|Y_0^N(\phi)| > A] &\leq \frac{1}{A^2} \mathbb{E}_{\mu^N}[|Y_0^N(\phi)|^2] \\
&= \frac{1}{A^2} \frac{1}{N} \Big( \sum_{x\in \Lambda_N} [\phi(\tfrac{x}{N}) ]^2 \mathbb{E}_{\mu^N} [\bar{\eta}_{0}(x)^2] + \sum_{\substack{x, y\in \Lambda_N\\ y \neq x}} \phi\left(\tfrac{x}{N}\right) \phi\left(\tfrac{y}{N}\right) \varphi^N_{0}(x,y) \Big).
\end{align*}
Using \eqref{decay_corr_time0_bulk} and the fact that the occupation variables are bounded by $\alpha$, we can bound the last display from above by a constant independent of $A$ and $N$ times $$\frac{1}{A^2 N}\left(\alpha^2N + N\right) \lesssim \frac{ 1}{A^2}.$$ Therefore by taking  $A \to \infty$ the result follows.

\subsubsection{The sequence of martingales}

For the martingales $\{M^N_t(\phi) \ ; \ t \in [0,T]\}_{N \in \mathbb{N}}$, tightness is just a consequence of the fact that $\{M^N_t(\phi) \ ; \ t \in [0,T]\}_{N \in \mathbb{N}}$ converges in law with respect to the Skorohod topology of $\mcb{D}([0,T], \mathbb{R})$ (see the next lemma) and therefore it has to be tight.

\begin{lemma} \label{lemma_conv_martingale}
For $\phi \in \mathcal{S}_i$, the sequence of martingales $\{M^N_t(\phi) \ ; \ t \in [0,T]\}_{N \in \mathbb{N}}$
converges in law with respect to
the topology of $\mcb{D}([0,T]; \mathbb{R})$, as $N \to +\infty$, towards a mean-zero Gaussian process $W_t^i(\phi)$
with quadratic variation given by
\[
\begin{split}
    \int_0^t \|\nabla \phi\|^2_{L^2(\rho_s)} ds 
    	&:= \int_0^t \int_0^1 2 \chi_\alpha(\rho_s(u))  \nabla \phi(u)^2 du ds \\
	&+\mathbb{1}(\theta = 1)  \int_0^t \Big\{\big(\lambda^\ell(\alpha-2\rho^\ell) \rho_{s}(0) +  \alpha\lambda^\ell\rho^\ell\big) \nabla \phi(0)^2 \\
    	&\quad \quad \quad \quad \quad + \big(\lambda^r(\alpha-2\rho^r) \rho_{s}(1) + \alpha\lambda^r\rho^r\big) \nabla \phi(1)^2 \big\} ds.
\end{split}
\]
\end{lemma}

\begin{proof}

Let us fix $\phi \in \mathcal{S}_i$. To prove that $\{M^N_t(\phi) \ ; \ t \in [0,T]\}_{N \in \mathbb{N}}$
converges in law with respect to the topology of $\mcb{D}([0,T]; \mathbb{R})$, as $N \to +\infty$, it is enough to verify conditions (1)-(3) of Theorem 3.2 of \cite{BGJS21}. 

Let us verify condition (1), that is, that
\begin{equation}
\label{maipu}
\text{for any $N > 1$, the quadratic variation of $M^N_t(\phi)$ has continuous trajectories almost surely.} 
\end{equation}
The quadratic variation of $M^N_t(\phi)$ is given by
\[
\langle M^N(\phi) \rangle_t
	: = \int_0^t \Gamma^N_s(\phi) ds,
\]
where $\Gamma^N_s(\phi) := N^2 \mcb{L}_N Y^N_s(\phi)^2 - 2N^2 Y^N_s(\phi) \mcb{L}_N Y^N_s(\phi)$. A long, but simple computation shows that this quadratic variation is given by
\begin{equation}
\label{renca}
\begin{split}	
\langle M^N(\phi) \rangle_t 
		&= \frac{N}{N^\theta}\int_0^t \Big(  \phi\Big(\frac{1}{N}\Big)^2 \big( \lambda^\ell(\alpha-2\rho^\ell)\eta_{sN^2}(1) + \alpha\lambda^\ell\rho^\ell\big)  \\
		&\quad \quad \quad \quad \quad + \phi\Big(\frac{N-1}{N}\Big)^2 \big( \lambda^r(\alpha-2\rho^r)\eta_{sN^2}(N-1) +  \alpha\lambda^r\rho^r\big) \Big) ds\\ 
		&+ \int_0^t \frac{1}{N}\sum_{x=1}^{N-2} \nabla_N \phi\Big(\frac xN\Big)^2 \Big(\eta_{sN^2}(x)\big(\alpha - \eta_{sN^2}(x+1)\big) + \eta_{sN^2}(x+1)\big(\alpha - \eta_{sN^2}(x)\big) \Big) ds,
\end{split}
\end{equation}
where
\begin{equation}
 \label{gradient_in_continuous}
\nabla_N\phi\Big(\tfrac{x}{N}\Big) := N\Big(\phi\Big(\frac{x+1}{N}\Big) - \phi\Big(\frac{x}{N}\Big)\Big)
\end{equation}
is the discrete gradient of $\phi$.
Therefore \eqref{maipu} follows from the fact that the number of particles is bounded by $\alpha$ and from the observation that the integral in time of a bounded function is a continuous function of time. 
    
Let us verify condition (2) in Theorem  3.2 of  \cite{BGJS21}, that is, that
 \begin{equation*}
 \lim_{N \to +\infty} \mathbb{E}_{\mu^N} \left[\sup_{0\leq s\leq T} |M^N_s(\phi) - M^N_{s^-}(\phi)| \right] = 0.
\end{equation*}

Observe that the integral term in \eqref{def_dynkin_martingale} is continuous, by exactly the same reason as in \eqref{maipu}. Therefore, in order to prove  the last limit, it is enough to show that
\begin{align*}
\lim_{N \to +\infty} \mathbb{E}_{\mu^N} \left[\sup_{0\leq s\leq T} |Y^N_s(\phi) - Y^N_{s^-}(\phi)|\right]=0.
\end{align*}
Since  a jump only changes a configuration in (at most) two sites and the occupation variables are bounded, we can bound the last expectation from above by $\frac{2}{\sqrt{N}} \|\phi'\|_\infty$, from where the result follows. 
 
We are left to verify condition (3) in Theorem 3.2 of \cite{BGJS21}, that is, that
\[
\text{for any $t \in [0, T]$, $\langle M^N(\phi)\rangle_t$ converges, as $N\to+\infty$, and in probability to $\int_0^t \|\nabla \phi\|^2_{L^2(\rho_s)} ds$.}
\]

Recall \eqref{renca}. We now argue that  $\int_0^t \Gamma^N_s (\phi)ds$  is an additive functional of the empirical measure plus some error that vanishes in the limit. To this end, we split the terms defining  $\Gamma^N_s (\phi)$ into bulk terms (the third line of \eqref{renca}) and boundary terms (the first two lines of \eqref{renca}). We present the argument for the leftmost term appearing in the bulk term, namely,
\begin{equation}\label{eq:bulk_left}
\begin{split}	
 \int_0^t \frac{1}{N}\sum_{x=1}^{N-2} \nabla_N \phi\Big(\frac xN\Big)^2 \eta_{sN^2}(x)\big(\alpha - \eta_{sN^2}(x+1)\big) ds,
\end{split}
\end{equation}
 but for the remaining one, it is completely analogous. The argument also extends to the boundary terms. We leave all this to the reader.   Let $0<\epsilon<1/2$ and 
 \begin{equation}\label{eq:sets}
 \Lambda^{\epsilon,\ell}_N := \{ 1, \dots, \epsilon (N-1)\}\quad \textrm{and}\quad\Lambda^{\epsilon,r}_N := \{ N-1-\epsilon (N-1), \dots, N-1\}\end{equation}
and  we consider
the sum  divided into $x\notin\Lambda_N^{\epsilon,r}\cup\Lambda_N^{\epsilon,\ell}$ and its complementary. Note that the terms in the complementary sets are uniformly (in $N$) bounded by $\epsilon$. Now, using twice   the replacement lemma (see Lemma 4.3 of \cite{FGS2022}, which we recall in  Lemma \ref{RL_lemma_bulk}) with proper choices of the function $\varphi$ appearing in the statement  of Lemma \ref{RL_lemma_bulk}, 
we can rewrite the terms in  \eqref{eq:bulk_left} for $x\notin \Lambda_N^{\epsilon,r}\cup\Lambda_N^{\epsilon,\ell}$ as
\begin{equation}
\begin{split}	
 \int_0^t \frac{1}{N}\sum_{x\notin \Lambda_N^{\epsilon,r}\cup\Lambda_N^{\epsilon,\ell}}  \nabla_N \phi\Big(\frac xN\Big)^2 \overleftarrow{\eta}^{\lfloor \epsilon N\rfloor}_{sN^2}(x)\Big(\alpha - \overrightarrow{\eta}^{\lfloor \epsilon N\rfloor}_{sN^2}(x+1)\Big) ds.
\end{split}
\end{equation}
Above, for $L\in\mathbb N$, \begin{equation}\label{eq:averages}
 \overrightarrow{\eta}^L(z) :=\frac{1}{L}\sum_{y=z+1}^{z+L}\eta(y) \quad \textrm{and}\quad \overleftarrow{\eta}^L(z) :=\frac{1}{L}\sum_{y=z-L}^{z-1}\eta(y).
 \end{equation}
Now it is enough to note that 
$
 \overrightarrow{\eta}^{\epsilon N}(x) =\langle\pi^N, \iota_{\epsilon}^{x/N}\rangle$ and similarly for the left average. {Above 
 $\iota_{\epsilon}^{x/N}(u):=\frac{1}{\epsilon}\mathbb {1}_{(x/N,x/N+\epsilon]}(u)$.}
From the fact that $\phi \in \mathcal{S}_i$ and the hydrodynamic limit namely Theorem \ref{th:hyd_ssep}, it follows the convergence in distribution, as $N\to+\infty$ and then $\epsilon\to 0$, to $\int_0^t \|\nabla \phi\|^2_{L^2(\rho_s)} ds$. Since the limit is deterministic, the convergence in probability also holds. 
\end{proof}

\subsubsection{The integral term}

Observe that, for every $\phi \in \mathcal{S}_i$,
\begin{align} \label{first_term_lambda_N_s_phi}
    \int_0^t (N^2 \mcb{L}_N + \partial_s) Y^N_s(\phi) ds 
    	&= \int_0^t Y^N_s( \alpha \Delta_N \phi) ds \\ 
\label{second_term_lambda_N_s_phi}
	&\quad- \int_0^t \frac{ \alpha N^{3/2}}{N^\theta}\left[\lambda^\ell\phi \left( \frac{1}{N}\right) \bar{\eta}_{sN^2}(1) + \lambda^r \phi \left( \frac{N-1}{N}\right) \bar{\eta}_{sN^2}(N-1) \right] ds\\ \label{third_term_lambda_N_s_phi}
&\quad \quad- \int_0^t \alpha \sqrt{N} \left[ \nabla_N \phi \left( \frac{N-1}{N} \right) \bar{\eta}_{sN^2}(N-1) - \nabla_N \phi \left(0 \right) \bar{\eta}_{sN^2}(1) \right]ds,
\end{align}
where, for every $x \in \Lambda_N$, 
\[
\Delta_N \phi\left(\frac{x}{N}\right) := N^2 \left[\phi\left(\frac{x+1}{N}\right) + \phi\left(\frac{x-1}{N}\right) - 2 \phi\left(\frac{x}{N}\right)\right]\footnote{By abuse of notation, we understand $Y_s^N(\alpha \Delta_N \phi)$ as the field $Y_s^N$ acting on any smooth function that coincides with $\alpha \Delta_N \phi$ in $\frac{1}{N} \Lambda_N$.}
\]
is the discrete Laplacian of $\phi$ evaluated at $\frac{x}{N}$. 
We will treat each of the integral terms  \eqref{first_term_lambda_N_s_phi}, \eqref{second_term_lambda_N_s_phi}, and \eqref{third_term_lambda_N_s_phi}, separately. We will rely on the \textit{Kolmogorov-Centsov's criterion}:

\begin{proposition}[Kolmogorov-Centsov criterion - Problem 2.4.11 of \cite{Karatzas}]
\label{KolCen}
A sequence $\{X^N_t ; t \in [0, T]\}_{ N \in \mathbb{N}}$ of
continuous, real-valued, stochastic processes is tight with respect to the uniform topology of $\mcb C([0, T];\mathbb{R})$ if the sequence of real-valued random variables $\{X^N_0\}_{N \in \mathbb{N}}$ is tight and there are constants $K, \gamma_1,\gamma_2 > 0$ such that, for any $t,s \in [0, T]$ and any $N \in \mathbb{N}$, it holds that
\begin{equation*}
    \mathbb{E}[|X^N_t - X^N_s|^{\gamma_1}] \leq K |t-s|^{1+\gamma_2}.
\end{equation*}
\end{proposition}

We start proving the tightness of \eqref{first_term_lambda_N_s_phi}: by the Cauchy–Schwarz inequality and Fubini's theorem, we have for every $t_1,t_2 \in [0,T]$ such that $t_1 < t_2$, that
\[
\begin{split}
\mathbb E_{\mu^N} \Bigg[ \Bigg( \int_{t_1}^{t_2} Y_s^N( \alpha \Delta_N \phi) ds\Bigg)^2 \Bigg] 
		&\leq (t_2-t_1) \int_{t_1}^{t_2} \mathbb E_{\mu^N} \big[ Y^N_s(\alpha \Delta_N \phi)^2 \big] ds\\
		&\lesssim \frac{t_2-t_1}{N}\int_{t_1}^{t_2} \sum_{x,y \in \Lambda_N} \mathbb E_{\mu^N}\big[\bar\eta_{sN^2}(x)\bar\eta_{sN^2}(y)\big]\Delta_N \phi \left(\frac{x}{N} \right) \Delta_N \phi \left(\frac{y}{N} \right) ds.
\end{split}
\]
Using the fact that the  occupation variables are bounded by $\alpha$ and {from} Proposition \ref{proposition_corr_decay}, last display is bounded from above by 
\begin{align} \label{toBound}
  C(t_2-t_1)^2\Big[ \alpha^2 \sup_{x \in \Lambda_N} \Delta_N \phi \left(\frac{x}{N} \right)^2 + \sup_{\substack{x,y \in \Lambda_N\\ y \neq x}} \Big|\Delta_N \phi \left(\frac{x}{N} \right) \Delta_N \phi \left(\frac{y}{N} \right) \Big| \Big],
\end{align}
for some constant $C$ independent of $N$.
Now, since $\phi \in \mathcal{S}_i \subseteq C^\infty([0,1])$, \eqref{toBound} is bounded {from} above by another constant times
\[
( \|\phi\|_\infty^2+ \|\phi''\|_\infty^2) (t_2-t_1)^2,
\]
which, by Proposition \ref{KolCen}, shows the tightness of \eqref{first_term_lambda_N_s_phi}. 

Let us now prove the tightness of the remaining terms, i.e.~\eqref{second_term_lambda_N_s_phi} and \eqref{third_term_lambda_N_s_phi}. 
We present the proof for the terms related to the left boundary of \eqref{second_term_lambda_N_s_phi} and \eqref{third_term_lambda_N_s_phi}; for the right boundary it is completely analogous. We start with the case $\theta = 1$. In this case we note that the terms related to the left boundary in \eqref{second_term_lambda_N_s_phi} and \eqref{third_term_lambda_N_s_phi} are equal to
     \begin{align} \nonumber
    &\int_0^t \alpha \sqrt{N}\left[ \lambda^\ell\phi \left( \frac{1}{N}\right) - \nabla_N \phi \left(0 \right) \right] \bar{\eta}_{sN^2}(1) ds.
    \end{align}
    Doing a Taylor expansion on $\phi$ at $x=0$ and noting that $\phi \in \mathcal S_i$, since the  occupation  variables  are bounded, we conclude that if  $X^N_t$ is defined as the integral term above, then
    \begin{equation}\label{criteria_Kolgomorov_Centsov_eq}
        \mathbb{E}[|X^N_t - X^N_s|^2]  \lesssim  |t-s|^2,
    \end{equation} 
    and tightness follows.

    Now we analyse the  case $\theta > 1$.  In this case it is enough to prove that $X^N_t$ defined as the  next integral term
\begin{equation*}
 \int_0^t \frac{N^{3/2}}{N^\theta} \alpha\lambda^\ell \phi \left( \frac{1}{N}\right) \bar{\eta}_{sN^2}(1)ds,
\end{equation*}  satisfies \eqref{criteria_Kolgomorov_Centsov_eq} with $\gamma_1 = 2$ and $\gamma_2 = \delta_\theta$ where $\delta_\theta$ is  defined in Lemma \ref{lemma_bound_corr_bound}. This result also implies that all the integral terms in \eqref{second_term_lambda_N_s_phi} and \eqref{third_term_lambda_N_s_phi} are  tight. But from Lemma \ref{lemma_bound_corr_bound}, we have that 
\begin{align*}
    \mathbb{E}\left[ \Big|\int_s^t \frac{N^{3/2}}{N^\theta} \alpha\lambda^\ell \phi \left( \frac{1}{N}\right) \bar{\eta}_{sN^2}(1)ds\Big| ^{2}\right]\lsim |t-s|^{1+\delta_\theta},
\end{align*} and we finish the proof for $\theta > 1$.

Now we go to the case $0 \leq \theta < 1$.  Note that since  $\phi \in \mathcal{S}_i$, then $\phi(0) = 0$. Thus
    \begin{align*}
    &\int_0^t \frac{N^{3/2}}{N^\theta} \alpha\lambda^\ell \phi \left( \frac{1}{N}\right) \bar{\eta}_{sN^2}(1)ds =  \int_0^t \frac{\sqrt N}{N^\theta}\alpha\lambda^\ell \nabla_N \phi \left( 0\right) \bar{\eta}_{sN^2}(1)ds.
    \end{align*}
Therefore, tightness in this case will follow if we show that 
\begin{equation*}
    \int_0^t \alpha \sqrt{N} \nabla_N \phi \left( 0\right)\bar{\eta}_{sN^2}(1)ds = \int_0^t \alpha \sqrt{N} \phi' \left( 0\right)\bar{\eta}_{sN^2}(1)ds + O\left(\frac{1}{\sqrt{N}}\right),
\end{equation*}
satisfies \eqref{criteria_Kolgomorov_Centsov_eq} with $\gamma_1 = 2$ and $\gamma_2 = \delta_\theta$ where $\delta_\theta$ is again defined as in Lemma \ref{lemma_bound_corr_bound}. This is a simple consequence of  Lemma \ref{lemma_bound_corr_bound}.
 
Finally, we treat the case $\theta < 0$. 
Note that now we need to prove tightness of
    \begin{align*}
    & \int_0^t \left[\frac{N^{3/2}}{N^\theta} \alpha\lambda^\ell \phi \Big( \frac{1}{N}\Big) \bar{\eta}_{sN^2}(1) - \alpha \sqrt{N} \nabla_N \phi \left(0 \right)\bar{\eta}_{sN^2}(1) \right]ds.
    \end{align*}
    From Lemma \ref{lemma_bound_corr_bound} the rightmost term in last display is tight. For the leftmost, we do a Taylor expansion of $\phi$ of order $\lfloor -\theta \rfloor +2$ around $x=0$, and we use that $\phi\in\mathcal S_i$, so that the leftmost term in last display writes as 
     \begin{align*}
    & \int_0^t \frac{N^{3/2}}{N^{\theta+\lfloor -\theta\rfloor +2}} \alpha\lambda^\ell \phi(t_N) \bar{\eta}_{sN^2}(1)  ds,
    \end{align*}
where $t_N$ is a point between $0$ and $1/N$. Since $3/2-\theta-\lfloor -\theta\rfloor -2<1/2$, then Lemma  \ref{lemma_bound_corr_bound} shows that the Kolmogorov–Centsov’s criteria is satisfied with $\gamma_1 = 2$ and $\gamma_2 = \min\{\delta_\theta,1\} > 0$  and tightness follows. 
This ends the proof of tightness.

\subsection{Characterization of the limit points}

Having proven tightness, we already know that there exists a subsequence {$\{\mathbb{Q}_{N_k}\}_{k \in \mathbb{N}}$ of $\{\mathbb{Q}_N\}_{N \in \mathbb{N}}$} which is convergent. Let us denote by $\mathbb{Q}$ its limit. We want now to characterize $\mathbb{Q}$. To do that, we will start by showing that $\mathbb{Q}$ gives probability one to all the paths of funcionals $\{Y_t \ | \ t \geq 0\}$ with a decomposition of the form \eqref{eq:charac_lim_point} - see Section \ref{decomposition_proof}. The strategy is to rewrite Dynkin's martingale $M^N_t$, see \eqref{def_dynkin_martingale},  applied to a particular test function $\phi$ defined in \eqref{choice_phi} and {to} prove that the integral term of $M^N_t$ goes to zero as $N \to +\infty$ in the $L^2(\mathbb{P}_{\mu^N})$-norm. This is what is done in the next {subsection}.

\subsubsection{Proof of the decomposition given in \eqref{eq:charac_lim_point}} \label{decomposition_proof}

Let $S^i_{t}$ be the semigroup associated to \eqref{hydrodynamic_equation}. We start by observing that, if $\lambda^\ell=\lambda^r=\alpha$, then $S^i_t = T^\theta_{\alpha t}$, where $T^\theta_{\alpha t}$ is the corresponding semigroup when taking in \eqref{hydrodynamic_equation} $\lambda^\ell=\lambda^r=1$  and that coincides with the semigroup taken in Definition 4 of \cite{GJMN}. In this case, due to the previous relation between semigroups, we can simply repeat the proof presented in case $\alpha = 1$ in \cite{FGN19} taking (for every fixed $t \in [0,T]$ and restricting the process to the time interval $[0,t]$) as test function 
\begin{equation} \label{choice_phi}
    \phi(u,s) := S^i_{t-s}f(u),
\end{equation}
where {$f \in \mathcal{S}_i$}, to obtain the decomposition of the limit point in the form
\begin{equation*}
    Y_t(f) = Y_0(S^i_t f) + W_t^i(f),
\end{equation*}
where $W_t^i(f)$ is the mean-zero Gaussian process characterized in Lemma \ref{lemma_conv_martingale}. For the previous choice of $\lambda^\ell=\lambda^r=\alpha$, this test function coincides with  $T^\theta_{\alpha(t-s)}f(u)$.

For completeness, we present here the  proof in the general case, which also follows the strategy of \cite{FGN19}. Taking $\phi_s(\cdot) = \phi(\cdot,s)$ defined in \eqref{choice_phi}, we have that 
\begin{align*}
    M^N_t(\phi_t) &= Y^N_t(\phi_t) - Y^N_0(\phi_0) -  \int_0^t [N^2 \mcb{L}_N Y^N_s(\phi_s) +  Y^N_s(\partial_s \phi_s)] ds
\end{align*} 
it is also a martingale. For every $s \in [0,T]$, if $f \in \mathcal{S}_i$, then $\phi_s \in \mathcal{S}_i$. Remarking that the proof of Lemma \ref{lemma_conv_martingale} still holds if the test function is time-dependent (and $C^1$ in time), we obtain that $\{M^N_t(\phi_t) \ ; \ t \in [0,T]\}_{N \in \mathbb{N}}$
converges in law with respect to
the topology of $\mcb{D}([0,T]; \mathbb{R})$, as $N \to +\infty$, towards a mean-zero Gaussian process $W_t^i(\phi_t) = W_t^i(f)$ with quadratic variation given by
\begin{align*}
    &\int_0^t ||\nabla \phi_s||^2_{L^2(\rho_s)} ds \\
    &:= \int_0^t \int_0^1 2 \chi_\alpha(\rho)(u) \left( \nabla \phi_s(u)\right)^2 du ds \\
    &+ \mathbb{1} (\theta = 1) \int_0^t \left\{\lambda^\ell[(\alpha-2\rho^\ell) \rho_{s}(0) + \rho^\ell \alpha] \left( \nabla \phi_s(0)\right)^2 + \lambda^r[(\alpha-2\rho^r) \rho_{s}(1) + \rho^r \alpha] \left( \nabla \phi_s(1)\right)^2 \right\} ds,
\end{align*}

Since, for every $N \in \mathbb{N}$,
\begin{equation} \label{To_converge}
     M^N_t(\phi_t) = Y^N_t(f) - Y^N_0(S^i_{t}f) -  \int_0^t [N^2 \mcb{L}_N Y^N_s(\phi_s)+ Y^N_s(\partial_s \phi_s)] ds,
\end{equation}  if we show that {the time integral in the  last display} goes to zero as $N \to +\infty$, then, using tightness and the previous reasoning about $\{M^N_t(\phi_t) \ ; \ t \in [0,T]\}_{N \in \mathbb{N}}$, taking the limit as $N \to +\infty$, we have, up to a subsequence, that \eqref{To_converge} converges in law with respect to the topology of $\mcb{D}([0,T]; \mathbb{R})$, to
\begin{equation*}
    W_t^i(f) = Y_t(f) - Y_0(S^i_{t}f),
\end{equation*} as we wanted. By the same computations done to obtain \eqref{first_term_lambda_N_s_phi}, \eqref{second_term_lambda_N_s_phi} and \eqref{third_term_lambda_N_s_phi}, we have
\begin{align} \nonumber
N ^2 \mcb{L}_N Y^N_s(\phi_s)+ Y^N_s(\partial_s \phi_s) &= \alpha 
     Y^N_s( \Delta_N S^i_{t-s}f - \Delta S^i_{t-s}f ) + Y^N_s( \alpha \Delta S^i_{t-s}f + \partial_s S^i_{t-s}f)\\ \label{eq_1_lambda_testfunction}
    &- \frac{ \alpha N^{3/2}}{N^\theta} \left[\lambda^\ell S^i_{t-s}f \left( \frac{1}{N}\right) \bar{\eta}_{sN^2}(1) + \lambda^r S^i_{t-s}f\left( \frac{N-1}{N}\right) \bar{\eta}_{sN^2}(N-1) \right]\\ \label{eq_2_lambda_testfunction}
    &- \alpha \sqrt{N} \left[ \nabla_N S^i_{t-s}f \left( \frac{N-1}{N} \right) \bar{\eta}_{sN^2}(N-1) - \nabla_N S^i_{t-s}f \left(0 \right) \bar{\eta}_{sN^2}(1) \right],
\end{align}
where $\Delta$ represents the continuous Laplacian operator. Since $S^i_{t-s}f$ is smooth (by the properties of the semigroup $S^i_{t-s}$), then $\Delta_N S^i_{t-s}f - \Delta S^i_{t-s}f$ is of order $O(N^{-2})$ and $\alpha \Delta S^i_{t-s}f + \partial_t S^i_{t-s}f$ is identically zero because $S^i_{t-s}f$ is solution to the heat equation with diffusion coefficient equal to $\alpha$ {with the corresponding boundary conditions depending on $\theta$} - recall \eqref{semigroup_eq_theta_bigger_1} for $\theta > 1$, \eqref{semigroup_eq_theta_equal_1} for $\theta = 1$, and \eqref{semigroup_eq_theta_less_1} for $\theta < 1$. It remains now to analyse the terms in \eqref{eq_1_lambda_testfunction} and  \eqref{eq_2_lambda_testfunction}. Here we treat the terms regarding the left boundary, since for the right boundary it is completely analogous.
\begin{enumerate}
    \item If $\theta = 1$, we have that
    \begin{align} \nonumber
        &- \frac{ \alpha N^{3/2}}{N^\theta}\lambda^\ell S^i_{t-s}f \left( \frac{1}{N}\right) \bar{\eta}_{sN^2}(1) + \alpha \sqrt{N} \nabla_N S^i_{t-s}f \left( \frac{1}{N} \right) \bar{\eta}_{sN^2}(1) \\ \nonumber
    &=\alpha \sqrt{N} \left[\nabla_N S^i_{t-s}f \left(0 \right) - \lambda^\ell S^i_{t-s}f \left( \frac{1}{N}\right) \right]\bar{\eta}_{sN^2}(1), \\ \label{order1_over2_2}
    &= \alpha \sqrt{N} \left[ \left(\nabla_N S^i_{t-s}f \left( \frac{1}{N}\right) - \partial_u S^i_{t-s}f (0) \right) - \lambda^\ell \left( S^i_{t-s}f\left( \frac{1}{N}\right) - S^i_{t-s}f\left( 0\right) \right) \right] \bar{\eta}_{sN^2}(1)  \\ \label{term0_2}
    &+ \alpha \sqrt{N} \left[ \partial_u S^i_{t-s}f \left( 0\right) - \lambda^\ell S^i_{t-s}f\left( 0\right) \right] \bar{\eta}_{sN^2}(1).
    \end{align}
    Since $S^i_{t-s}f$ is smooth, both terms in \eqref{order1_over2_2} are of order $O(N^{-1/2})$ and \eqref{term0_2} is identically zero because $S^i_{t-s}f$ satisfies the boundary conditions given in \eqref{semigroup_eq_theta_equal_1}. This immediately implies that, if $\theta = 1$, then $\displaystyle
    \int_0^t [N ^2 \mcb{L}_N Y^N_s(\phi_s)+ Y^N_s(\partial_s \phi_s)] ds$ goes to zero as $N \to +\infty$.

    \item If $\theta > 1$, since $f \in \mathcal{S}_i$ and so $S^i_t f \in \mathcal{S}_i$, we have that
    \begin{align} \nonumber
        &- \frac{N^{3/2}}{N^\theta}\alpha \lambda^\ell S^i_{t-s}f \left( \frac{1}{N}\right)\bar{\eta}_{sN^2}(1) + \alpha \sqrt{N} \nabla_N S^i_{t-s}f \left( \frac{1}{N} \right) \bar{\eta}_{sN^2}(1) \\\label{1}
		&= -\frac{\alpha \lambda^\ell}{N^{\theta-1/2}} \nabla_N S^i_{t-s}f \left(0 \right)\bar{\eta}_{sN^2}(1) + \alpha \sqrt{N} \left(\nabla_N S^i_{t-s}f \left(0 \right)  - \partial_u S^i_{t-s}f \left(0 \right) \right) \bar{\eta}_{sN^2}(1)\\  \label{2}
		&- N^{3/2-\theta}\alpha \lambda^\ell S^i_{t-s}f \left(0\right) \bar{\eta}_{sN^2}(1).
    \end{align}
    Since $S^i_{t-s}f$ is smooth and {the occupation variables are bounded,} then the first term of \eqref{1} is of order $O(N^{1/2-\theta})$ and the second is of order $O(N^{-1/2})$. Finally, integrating \eqref{2} between $0$ and $t$, and taking its $L^2(\mathbb{P}_{\mu_N})$-norm, by the Lemma \ref{lemma_bound_corr_bound} we conclude that the integral between $0$ and $t$ of this term goes to zero as $N \to +\infty$, and we are done.

    \item If $0 \leq \theta < 1$, by the invariance of the semigroup $S^i_t$ in $\mathcal{S}_i$, we have that
    \begin{align} \nonumber
        &- \frac{N^{3/2}}{N^\theta}\alpha \lambda^\ell S^i_{t-s}f \left( \frac{1}{N}\right)\bar{\eta}_{sN^2}(1) + \alpha \sqrt{N} \nabla_N S^i_{t-s}f \left( \frac{1}{N} \right) \bar{\eta}_{sN^2}(1) \\ \label{22}
        &= -\frac{\sqrt{N}}{N^{\theta}} \alpha \lambda^\ell \nabla_N S^i_{t-s}f \left(0 \right)\bar{\eta}_{sN^2}(1) + \alpha \sqrt{N} \nabla_N S^i_{t-s}f \left(0 \right) \bar{\eta}_{sN^2}(1).
    \end{align} Integrating both terms in \eqref{22} between $0$ and $t$, and taking the $L^2(\mathbb{P}_{\mu_N})$-norm of each term, by Lemma \ref{lemma_bound_corr_bound}, the integral between $0$ and $t$ of these terms go to zero as $N \to +\infty$. We can then conclude that, if $0 \leq \theta < 1$, then $\displaystyle
    \int_0^t [N ^2 \mcb{L}_N Y^N_s(\phi_s)+ Y^N_s(\partial_s \phi_s)] ds$ goes to zero as $N \to +\infty$.

    \item Finally, if $\theta < 0$, since $f \in \mathcal{S}_i$ implies that $S^i_{t-s}f \in \mathcal{S}_i$, then, writing the Taylor expansion of order $\lceil -\theta \rceil + 1$ of $S^i_{t-s}f$ around $0$ and substituting in \eqref{eq_1_lambda_testfunction} and \eqref{eq_2_lambda_testfunction}, we immediatly conclude that $\displaystyle
    \int_0^t [N ^2 \mcb{L}_N Y^N_s(\phi_s) + Y^N_s(\partial_s \phi_s)] ds$ goes to zero as $N \to +\infty$.
\end{enumerate} This completes the proof of the decomposition part of Theorem \ref{th_flutuations}.

Putting all the previous results together we end the proof of Theorem \ref{th_flutuations}. What distinguishes the two main theorems is the fact that in the first one, there is a  convergence result (only up to subsequences) since we are not able to show uniqueness of solution of the martingale problem. Nevertheless, since in Theorem \ref{th_conv_OU} we assume a convergence at the initial time, this gives the uniqueness of the limit point. In the next subsection we complete the proof of Theorem \ref{th_conv_OU} by showing uniqueness of the limit.
\subsection{Proof of Theorem \ref{th_conv_OU}} \label{proof_th_2_fluct}

The uniqueness of the O. U. process for $\theta \geq 0$ follows from Proposition 2.5 of \cite{BGJS21} once we show that $(S^i_t)_{t \geq 0}$, the semigroup associated to \eqref{hydrodynamic_equation}, satisfies
\begin{equation}
S^i_{t+\epsilon} H - S^i_t H = \epsilon \alpha \Delta S^i_t H + o(\epsilon, t),
\end{equation} for every $\epsilon > 0$, $t\geq 0$ and $H \in \mathcal{S}_i$, where $o(\epsilon, t)$ goes to 0, as $\epsilon$ goes to 0, in $\mathcal{S}_i$ uniformly on compact time intervals. But this is an immediate consequence of the explicit formulas given by \eqref{semigroup_formula_theta_bigger_1}, \eqref{semigroup_formula_theta_1} and \eqref{semigroup_formula_theta_lower_1}, if $\theta > 1$, $\theta =1$ or $\theta < 1$, respectively. Moreover, for $\theta < 0$, the uniqueness of solution of the O. U. martingale problem follows by repeating the arguments of Theorem 2.13. of \cite{BGJS21} and Proposition 2.5. of \cite{BGJS21}. Finally, to show that the two extra conditions, i.e.~\textit{regularity} and \textit{boundary conditions}, hold, we only have to observe that the first follows from the boundedness of the occupation variables jointly with Proposition \ref{proposition_corr_decay} and the second follows from Lemma \ref{statementlemma2}. This finishes the proof of Theorem \ref{th_conv_OU}.

\section{Auxiliary estimates} \label{section_after_proof_main_theorems}
This section is devoted to some estimates needed in order to proof our main results.
Let us denote by $\tilde{\nabla}^+_N$ the operator defined, for every $f:\Lambda_N \to \mathbb{R}$ and $x \in \Lambda_{N-1}$, by
\begin{align} \label{discrete_gradient}
    \tilde{\nabla}^+_N f(x) := N[f(x+1)-f(x)].
\end{align}

\begin{lemma} \label{lemma_bounded_discrete_grad_rho}
Assume that $\gamma\in C^6([0,1])$ satisfies \eqref{H1_3}, that there exists a sequence $(g_N)_{N \in \mathbb{N}}$ of functions of class $C^6([0,1])$ that satisfies \eqref{H1_6} and \eqref{H2} and that  $(\mu^N)_{N \in \mathbb{N}}$ is a sequence of probability measures satisfying \eqref{H0}. Then, there exists $C > 0$ such that
$$\max\limits_{x \in \Lambda_{N-1}} |\tilde{\nabla}_N^+ \rho^N_t(x)| \leq C,$$
for every $t \in [0,T]$.
\end{lemma}

The proof of the previous lemma can be found in Appendix \ref{edp_result_bound_discrete_derivative}.

One of the key ingredients to prove fluctuations is to obtain sharp estimates for the decay in $N$ of the time-dependent two-point correlation function, i.e.~on $\varphi^N_t$ defined in \eqref{time_dependent_correlation_def}, which we recall  that is not defined for $x=y$.

\begin{proposition} \label{proposition_corr_decay}
Under the assumption  \eqref{decay_corr_time0_bulk}, we have that
\begin{align} \label{results_bounds_corr_bulk}
    \sup_{t \in [0,T]} \max_{\substack{(x,y) \in V_N\\x\neq y }} |\varphi^N_t(x,y)| \lsim 
    \frac{1}{N},
\end{align}
and, under the assumption \eqref{decay_corr_time0_bound}, for $x = 1$ and for $x= N-1$,
\begin{align} \label{results_bounds_corr_bound}
    \sup_{t \in [0,T]} \max_{\substack{y \in \Lambda_N \\ y \neq x}} |\varphi^N_t(x,y)| \lsim R_N^\theta:= \begin{cases}
    \frac{1}{N}, \textrm{ if } \theta > 1,\\
    \frac{N^\theta}{N^2}, \textrm{ if } 0 \leq \theta \leq 1,\\
    \frac{N^\theta}{N}, \textrm{ if } -1 < \theta < 0,\\
    \frac{1}{N^2}, \textrm{ if } \theta \leq -1.
    \end{cases}
\end{align}

\end{proposition}

The proof of the previous proposition can be found in Section \ref{_correlation_decay_proof}.

\begin{lemma} \label{lemma_bound_corr_bound}
Recall that, {for $y \in \Lambda_N$, we denote by $\bar{\eta}(y)$} the centered variable.
Then, for every $\theta \in \mathbb{R}$, for $x \in \{1,N-1\}$ and $t,s \in [0,T]$, it holds 
\begin{equation} \label{bound_to_use_on_KC}
    \mathbb{E}_{\mu^N} \left[ \left( \int_s^t d_N^\theta \bar{\eta}_{sN^2}(x) dr\right)^2\right] \lsim |t-s|^{1+\delta_\theta} + |t-s|^2 (d^\theta_N)^2 R_N^\theta
\end{equation}
and
\begin{equation} \label{bound_to_use_on_CLM}
    \mathbb{E}_{\mu^N} \left[ \left( \int_s^t \bar{\eta}_{sN^2}(x) dr\right)^2\right] \lsim \frac{N^\theta}{N^2}|t-s| + |t-s|^2 R_N^\theta,
\end{equation}
where $d_N^\theta = \sqrt{N} \mathbb{1} (\theta \leq 1) + N^{3/2 - \theta} \mathbb{1} (\theta > 1)$, $\delta_\theta = \frac{|1-\theta|}{2} \mathbb{1} (\theta < 3) + \mathbb{1} (\theta \geq 3)$ and $R_N^\theta$ was introduced in the last proposition.
So, in particular, for $x \in \{1,N-1\}$, for every $t \in [0,T]$ and $\theta \in \mathbb{R}$,
\begin{align} \label{limit_to_use_on_CLM}
    \lim_{N \to +\infty} \mathbb{E}_{\mu^N} \left[ \left( \int_0^t d_N^\theta \bar{\eta}_{sN^2}(x) dr\right)^2\right] = 0.
\end{align}
\end{lemma}

The proof of the previous lemma is given in Section \ref{prob_estimates_Lemma1}.

For $\theta < 0$, for all $\alpha \in \mathbb{N}$, we will also need the following estimates.

\begin{proposition} \label{proposition_corr_decay_theta_less_1}
Let $\theta < 1$. Recall \eqref{eq:sets}.
  If \eqref{decay_corr_time0_bulk} holds, then, for every $\epsilon > 0$ and every $t \in (0,T]$, we have that
\begin{align} \label{results_bounds_Gt_bulk_close_bound}
    \max_{\substack{(x,y) \in \Lambda^{\epsilon,\ell}_N \times \Lambda_N \\ y \neq x} }|\varphi^N_t(x,y)| \lsim \left( 1 + \frac{1}{\sqrt{t}} \right)\frac{\epsilon}{N} + o\left(\frac{1}{N}\right),
\end{align}
and the same results holds for $(x,y) \in  \Lambda_N \times \Lambda^{\epsilon,r}_N$. 
\end{proposition}

The proof of the previous result can be found in Section \ref{proof_proposition_corr_decay_theta_less_1}.

\begin{lemma} \label{statementlemma2}
Let $\theta < 1$. Then, the following limit holds, for every $t \in [0,T]$ and $j \in \{0,1\}$.
\begin{align} \label{convergence_in_l2_flu_i_epsilon}
    \lim_{\epsilon \to 0} \limsup_{N \to +\infty} \mathbb{E}_{\mu_N}\left[ \left(\int_0^t Y^N_s(\iota^j_\epsilon) ds \right)^2  \right] = 0,
\end{align}
where $\iota^j_\epsilon$ was defined in item 2. (b) of Theorem \ref{th_conv_OU}.
\end{lemma}

The proof of the previous result is given in Section \ref{prob_estimates_Lemma2}.

\subsection{Proof of Proposition \ref{proposition_corr_decay}} \label{_correlation_decay_proof}

In this proof we will use some random walks that, for simplicity of the presentation, we define now:
\begin{enumerate}
\item $\{\mathcal  X_{t}^i; t \geq 0\}$ is the random walk evolving on the set of points $V^\alpha_N$ where \begin{equation}\label{new_V}
V^\alpha_N:= V_N\setminus \mathcal{D}_N
 \quad \textrm{  for } \alpha=1 \quad \textrm{and}\quad  V^\alpha_N:= V_N\quad \textrm{for }\quad \alpha\geq 2,
 \end{equation}
 that moves to nearest-neighbors at rate $\alpha$, except at the line $\mathcal D_N^+$ that moves right/up at rate $\alpha$ and left/down at rate $\alpha-1$ and that is reflected at the line $\mathcal D_N^+$ if $\alpha=1$, and at the line $\mathcal{D}_N$ if $\alpha\geq 2$. Moreover, it is absorbed at $\partial V_N$: with rate $\alpha\lambda^\ell/N^\theta$ at the set of points $\{(0,y): y\in{\overline{\Lambda}}_N\}$ and with rate $\alpha\lambda^r/N^\theta$ at the set of points $\{(x,N): x\in \overline{\Lambda}_N\}$. This random walk has generator $\Delta^{i}_N$ which is the operator that acts on functions $f:\overline{V}_N \to\mathbb R$ such that $f(x,y) =0$ for every $(x,y) \in \partial V_N$ as
 \begin{equation}\label{op_A}
 \Delta^{i}_N f(u)=\sum_{\substack{v \in \overline{V}_N \\v \sim u}} c^{i}_{u,v}\big[ f(v)-f(u)\big],
 \end{equation}
for every $u \in V_N$, with $c^i_{x,y}$ defined, for $\alpha=1$ by
\[
c^i\!: \big\{((x,y),(x',y')) \in V_N \times \overline{V}_N;  |x-x'|+|y-y'| =1\big\} \to [0,\infty) 
\]
as
\[
\left\{
\begin{array}{r@{\;\text{if}\;}l}
c^i_{(x,y),(x',y)} := c^i_{x,x'} \mathbb{1}(x' \neq y) & |x-x'|=1,\\
c^i_{(x,y),(x,y')} := c^i_{y,y'} \mathbb{1}(x \neq y') & |y-y'|=1,
\end{array}
\right.
\] and, for $\alpha\geq 2$,
\begin{equation} \label{c_i_2}
\left\{
\begin{array}{r@{\;\text{if}\;}l}
c^i_{(x,y),(x',y)} := c^i_{x,x'} - \mathbb{1}(x' = y) & |x-x'|=1,\\
c^i_{(x,y),(x,y')} := c^i_{y,y'} - \mathbb{1}(x = y') & |y-y'|=1,
\end{array}
\right. \end{equation} 
with $c^i_{x,y}$ as defined in equation \eqref{ratesss}.

\begin{figure}[h!] \label{fig:triangle}
\captionsetup[subfigure]{font=footnotesize}
\centering
\subcaptionbox{Illustration of the jump rates $c^i$ of the random walk $\{\mathcal  X_{t}^i; t \geq 0\}$ when $\alpha=1$.}[.45\textwidth]{%
\begin{tikzpicture}[thick, scale=1]
\draw[->] (-0.5,0)--(5.5,0) node[anchor=north]{$x$};
\draw[->] (0,-0.5)--(0,5.5) node[anchor=east]{$y$};
\begin{scope}[scale=0.75]
\foreach \x in {1,...,3} 
	\foreach \y in {\x,...,3,4}
		\shade[ball color=black](\x,1+\y) circle (0.15);
  \foreach \x in {1,...,3,4} 
	\foreach \y in {\x}
		\shade[ball color=black](\x,1+\y) circle (0.15);
\foreach \x in {1,...,5} 
	\shade[ball color = red](\x,6) circle (0.15); 
\foreach \y in {1,...,5} 
	\shade[ball color = red](0,\y) circle (0.15); 
\shade[ball color= red](6,6) circle (0.15);
\shade[ball color= red](0,6) circle (0.15);
\shade[ball color= red](0,0) circle (0.15);
\end{scope}	
\draw (0,0) node[anchor=north east] {$0$};
\draw (0.75,2pt)--(0.75,-2pt) node[anchor=north] {$1$};
\draw (1.5,2pt)--(1.5,-2pt) node[anchor=north]{$2$};
\draw (2.25,2pt)--(2.25,-2pt);
\draw (3,2pt)--(3,-2pt);
\draw (3.75,2pt)--(3.75,-2pt) node[anchor=north]{\small $N-1$};
\draw (4.5,2pt)--(4.5,-2pt) node[anchor=north]{\small $N$};
\draw (-0.05,.75) node[anchor=east] {$1$};
\draw (-0.05,1.5) node[anchor=east]{$2$};
\draw (-0.05,3.75) node[anchor= east]{$N-1$};
\draw (-0.05,4.5) node[anchor=east]{$N$};
\node at (1.5,3) (C) { };
\node at (0.75,3) (Cl) { };
\node at (2.25,3) (Cr) { };
\node at (1.5,3.75) (Cup) { };
\node at (1.5,2.25) (Cdow) { };
\node at (0.75,3.75) (B) { };
\node at (0,3.75) (Bl) { };
\node at (0.75,4.5) (Bup) { };
\node at (3,3.75) (D) { };
\node at (2.25,3.75) (Dl) { };
\node at (3,4.5) (Dup) { };
\node at (0.75,2.25) (El) { };
\draw[-latex] (C) to[out=130,in=60] node[midway,font=\scriptsize,above] {$1$} (Cl);
\draw[-latex] (C) to[out=-60,in=-120] node[midway,font=\scriptsize,below=0.2cm] {$1$} (Cr);
\draw[-latex] (C) to[out=60,in=-60] node[midway,font=\scriptsize,right] {$1$} (Cup);
\draw[-latex] (C) to[out=-130,in=130] node[midway,font=\scriptsize,left] {$1$} (Cdow);
\draw[-latex] (B) to[out=130,in=60] node[midway,font=\scriptsize,above] {$\frac{\lambda^\ell}{N^\theta}$} (Bl);
\draw[-latex] (D) to[out=130,in=60] node[midway,font=\scriptsize,above] {$1$} (Dl);
\draw[-latex] (D) to[out=60,in=-60] node[midway,font=\scriptsize,right] {$\frac{\lambda^r}{N^\theta}$} (Dup);
\draw[-latex] (Cdow) to[out=-120,in=-60] node[midway,font=\scriptsize,below] {$1$} (El);
\draw[-latex] (Cdow) to[out=60,in=-60] node[midway,font=\scriptsize,below=0.1cm,right=0.1cm] {$1$} (C);
\end{tikzpicture}}%
\hspace*{0.1cm}
\subcaptionbox{Illustration of the jump rates $c^i$ of the random walk $\{\mathcal  X_{t}^i; t \geq 0\}$ when $\alpha\geq 2$.}[.45\textwidth]{
\begin{tikzpicture}[thick, scale=1]
\draw[->] (-0.5,0)--(5.5,0) node[anchor=north]{$x$};
\draw[->] (0,-0.5)--(0,5.5) node[anchor=east]{$y$};
\begin{scope}[scale=0.75]
\foreach \x in {1,...,3} 
	\foreach \y in {\x,...,3,4}
		\shade[ball color=black](\x,1+\y) circle (0.15);
  \foreach \x in {1,...,3,4} 
	\foreach \y in {\x}
		\shade[ball color=black](\x,1+\y) circle (0.15);
\foreach \x in {1,...,5} 
	\shade[ball color = red](\x,6) circle (0.15); 
\foreach \x in {1,...,5} 
	\shade[ball color = gray!50](\x,\x) circle (0.15); 
\foreach \y in {1,...,5} 
	\shade[ball color = red](0,\y) circle (0.15); 
\shade[ball color= red](6,6) circle (0.15);
\shade[ball color= red](0,6) circle (0.15);
\shade[ball color= red](0,0) circle (0.15);
\end{scope}	
\draw (0,0) node[anchor=north east] {$0$};
\draw (0.75,2pt)--(0.75,-2pt) node[anchor=north] {$1$};
\draw (1.5,2pt)--(1.5,-2pt) node[anchor=north]{$2$};
\draw (2.25,2pt)--(2.25,-2pt);
\draw (3,2pt)--(3,-2pt);
\draw (3.75,2pt)--(3.75,-2pt) node[anchor=north]{\small $N-1$};
\draw (4.5,2pt)--(4.5,-2pt) node[anchor=north]{\small $N$};
\draw (-0.05,.75) node[anchor=east] {$1$};
\draw (-0.05,1.5) node[anchor=east]{$2$};
\draw (-0.05,3.75) node[anchor= east]{$N-1$};
\draw (-0.05,4.5) node[anchor=east]{$N$};
\node at (1.5,3) (C) { };
\node at (0.75,3) (Cl) { };
\node at (2.25,3) (Cr) { };
\node at (1.5,3.75) (Cup) { };
\node at (1.5,2.25) (Cdow) { };
\node at (0.75,3.75) (B) { };
\node at (0,3.75) (Bl) { };
\node at (0.75,4.5) (Bup) { };
\node at (3,3.75) (D) { };
\node at (2.25,3.75) (Dl) { };
\node at (3,4.5) (Dup) { };
\node at (3,3.75) (Di2) { };
\node at (3,3) (Dow2) { };
\node at (1.5,2.25) (Di) { };
\node at (0.75,2.25) (Dil) { };
\node at (2.25,2.25) (Dir) { };
\node at (1.5,3) (Diup) { };
\node at (1.5,1.5) (Dow) { };
\draw[-latex] (Dow2) to[out=60,in=-60] node[midway,font=\scriptsize,right] {$2\alpha$} (Di2);
\draw[-latex] (Dow2) to[out=-120,in=-60] node[midway,font=\scriptsize,below] {$2\alpha$} (Cr);
\draw[-latex] (Di) to[out=-60,in=-120] node[midway,font=\scriptsize,below,color=blue] {$\alpha-1$} (Dir);
\draw[-latex] (Di) to[out=-130,in=130] node[midway,font=\scriptsize,left,,color=blue] {$\alpha-1$} (Dow);
\draw[-latex] (Di2) to[out=-130,in=130] node[midway,font=\scriptsize,left,color=blue] {$\alpha-1$} (Dow2);
\draw[-latex] (C) to[out=130,in=60] node[midway,font=\scriptsize,above] {$\alpha$} (Cl);
\draw[-latex] (C) to[out=-60,in=-120] node[midway,font=\scriptsize,below] {$\alpha$} (Cr);
\draw[-latex] (C) to[out=60,in=-60] node[midway,font=\scriptsize,right] {$\alpha$} (Cup);
\draw[-latex] (C) to[out=-130,in=130] node[midway,font=\scriptsize,left] {$\alpha$} (Cdow);
\draw[-latex] (B) to[out=130,in=60] node[midway,font=\scriptsize,above] {$\frac{\alpha \lambda^\ell}{N^\theta}$} (Bl);
\draw[-latex] (D) to[out=130,in=60] node[midway,font=\scriptsize,above] {$\alpha$} (Dl);
\draw[-latex] (D) to[out=60,in=-60] node[midway,font=\scriptsize,right] {$\frac{\alpha \lambda^r}{N^\theta}$} (Dup);
\end{tikzpicture}}%
\end{figure}


\item $\{\widetilde{\mathcal  X}_{t}; t \geq 0\}$ is the random walk evolving  on the set of points $V_N$ that moves to nearest-neighbors at rate $\alpha$,  except at the line $\mathcal D_N^+$ that moves right/up at rate $\alpha$ and left/down at rate $\alpha-1$ and that is reflected at the  line $\mathcal D_N$  and at $\partial V_N$. We denote by $\mathfrak C^i_N$ the Markov generator of 
$\{\widetilde{\mathcal  X}_{t}^i; t \geq 0\}$ which is the operator that acts on functions $f: \overline{V}_N \to \mathbb{R}$ as, for every $u \in V_N$,
\begin{equation} \label{op_C}
 \mathfrak C^i_N f (u)=\sum_{\substack{v \in V_N \\v \sim u}} c^{i}_{u,v}\big[ f(v)-f(u)\big],
 \end{equation} where $c^{i}_{u,v}$ are the same as defined in \eqref{c_i_2}.

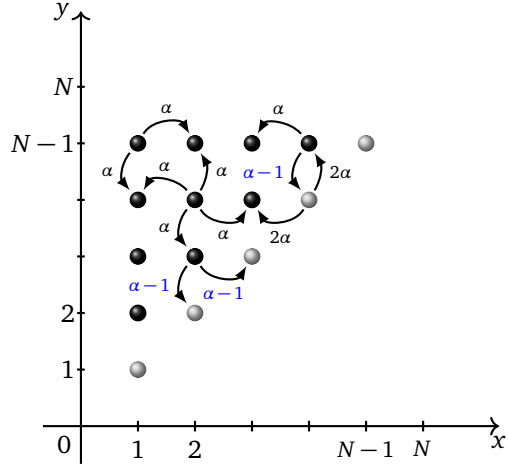
\begin{figure}[h!]
    \centering
\begin{tikzpicture}[thick, scale=1]
\draw[->] (-0.5,0)--(5.5,0) node[anchor=north]{$x$};
\draw[->] (0,-0.5)--(0,5.5) node[anchor=east]{$y$};
\begin{scope}[scale=0.75]
\foreach \x in {1,...,3} 
	\foreach \y in {\x,...,3,4}
		\shade[ball color=black](\x,1+\y) circle (0.15);
  \foreach \x in {1,...,3,4} 
	\foreach \y in {\x}
		\shade[ball color=black](\x,1+\y) circle (0.15);
\foreach \x in {1,...,5} 
	\shade[ball color = gray!50](\x,\x) circle (0.15); 
\end{scope}	
\draw (0,0) node[anchor=north east] {$0$};
\draw (0.75,2pt)--(0.75,-2pt) node[anchor=north] {$1$};
\draw (1.5,2pt)--(1.5,-2pt) node[anchor=north]{$2$};
\draw (2.25,2pt)--(2.25,-2pt);
\draw (3,2pt)--(3,-2pt);
\draw (3.75,2pt)--(3.75,-2pt) node[anchor=north]{\small $N-1$};
\draw (4.5,2pt)--(4.5,-2pt) node[anchor=north]{\small $N$};
\draw (-0.05,.75)--(0.05,.75) node[anchor=east] {$1$};
\draw (-0.05,1.5)--(0.05,1.5) node[anchor=east]{$2$};
\draw (-0.05,2.25)--(0.05,2.25);
\draw (-0.05,3)--(0.05,3);
\draw (-0.05,3.75)--(0.05,3.75) node[anchor= east]{$N-1$};
\draw (-0.05,4.5)--(0.05,4.5) node[anchor=east]{$N$};
\node at (1.5,3) (C) { };
\node at (0.75,3) (Cl) { };
\node at (2.25,3) (Cr) { };
\node at (1.5,3.75) (Cup) { };
\node at (1.5,2.25) (Cdow) { };
\node at (0.75,3.75) (Cor) { };
\node at (0.75,3.75) (B) { };
\node at (0,3.75) (Bl) { };
\node at (0.75,4.5) (Bup) { };
\node at (3,3.75) (D) { };
\node at (2.25,3.75) (Dl) { };
\node at (3,4.5) (Dup) { };
\node at (3,3.75) (Di2) { };
\node at (3,3) (Dow2) { };
\node at (1.5,2.25) (Di) { };
\node at (0.75,2.25) (Dil) { };
\node at (2.25,2.25) (Dir) { };
\node at (1.5,3) (Diup) { };
\node at (1.5,1.5) (Dow) { };
\node at (3,3.75) (Di2) { };
\node at (3,3) (Dow2) { };
\draw[-latex] (Cor) to[out=-130,in=130] node[midway,font=\scriptsize,left] {$\alpha$} (Cl);
\draw[-latex] (Cor) to[out=60,in=1200] node[midway,font=\scriptsize,above] {$\alpha$} (Cup);
\draw[-latex] (Di) to[out=-60,in=-120] node[midway,font=\scriptsize,below,color=blue] {$\alpha-1$} (Dir);
\draw[-latex] (Di) to[out=-130,in=130] node[midway,font=\scriptsize,left,,color=blue] {$\alpha-1$} (Dow);
\draw[-latex] (Di2) to[out=-130,in=130] node[midway,font=\scriptsize,left,color=blue] {$\alpha-1$} (Dow2);
\draw[-latex] (C) to[out=130,in=60] node[midway,font=\scriptsize,above] {$\alpha$} (Cl);
\draw[-latex] (C) to[out=-60,in=-120] node[midway,font=\scriptsize,below] {$\alpha$} (Cr);
\draw[-latex] (C) to[out=60,in=-60] node[midway,font=\scriptsize,right] {$\alpha$} (Cup);
\draw[-latex] (C) to[out=-130,in=130] node[midway,font=\scriptsize,left] {$\alpha$} (Cdow);
\draw[-latex] (D) to[out=130,in=60] node[midway,font=\scriptsize,above] {$\alpha$} (Dl);
\draw[-latex] (Dow2) to[out=60,in=-60] node[midway,font=\scriptsize,right] {$2\alpha$} (Di2);
\draw[-latex] (Dow2) to[out=-120,in=-60] node[midway,font=\scriptsize,below] {$2\alpha$} (Cr);
\end{tikzpicture}
\caption{Illustration of the jump rates of the random walk $\{\widetilde{\mathcal  X}_{t}^i; t \geq 0\}$.}
    \label{fig:triangle}
\end{figure}

\end{enumerate}

For the standard simple symmetric exclusion process, i.e.~the case $\alpha=1$, Proposition \ref{proposition_corr_decay} has been proved in a myriad of articles (see  \cite{LMO,FGN19,GJMN} and references therein). Let us review and adapt this proof. It is not difficult to check that for each $x, y \in \Lambda_N$, the action of the generator $\mcb L_N$ on $\eta(x)\eta(y)$ is given by a linear combination of the functions $(\eta(z) \eta(z'); z,z' \in \Lambda_N)$ - see equation \eqref{eq_action_gen_eta_x_eta_y} of Appendix \ref{equation_correlations_computation}. This means that the correlation function $(\varphi_t^N; t \geq 0)$ satisfies an autonomous, non-homogeneous evolution equation, which involves $(\rho_t^N; t \geq 0)$ as parameters.

For $\alpha = 1$, the correlation function $\varphi^N_t$ is solution to
\begin{align}
 \label{eq_varphi_alpha1}
\partial_t \varphi_t^N(x, y) 
		=  N^2 \Delta^{i}_N \varphi_t^N(x, y) 
			+ g^{N}_t(x,y) \mathbb {1} ((x,y) \in \mathcal{D}_N^{\pm}),
\end{align} 
where $\Delta^{i}_N$ is the operator defined in \eqref{op_A}. Here
\begin{equation*}
    g^{N}_t(x,x+1) = g^{N}_t(x+1,x) = -\left({\tilde{\nabla}}^+_N\rho_t^N(x)\right)^2,
\end{equation*} 
for every $x \in \Lambda_{N-1}$ and $g_t^N(x,y) :=0$ otherwise.
Observe that $\Delta^i_N$ corresponds to the generator of the random walk $\{\mathcal  X_{t}^i; t \geq 0\}$ that moves to nearest-neighbor sites on $V_N$ with annihilation at the boundary and the jumps to the diagonal $\mathcal{D}_N$ are suppressed. {As a} consequence, 
\eqref{eq_varphi_alpha1} does not involve the values of $\varphi_t^N$ at $\mathcal{D}_N$.
By Duhamel's formula, for every $(x,y) \in V_N \setminus \mathcal D_N$, we can represent $\varphi_t^N$ by
{\begin{align} \label{varphi_SEP_1_rep}
    \varphi_t^N(x,y) = \mathbb{E}_{(x,y)} \Big[\varphi_0^N(\mathcal  X_{tN^2}^{i})+ \int_0^t g^{N}_{t-s}( \mathcal  X_{sN^2}^{i})\mathbb {1}( \mathcal  X_{sN^2}^{i} \in \mathcal{D}_N^{\pm}) \,ds\Big],
\end{align}}
where $\mathbb{E}_{(x,y)}$ denotes the expectation of the law of the walk $\{\mathcal  X_{t}^i; t \geq 0\}$ starting from the point $(x,y)$. Now, to obtain the order of decay in $N$ of $\varphi^N_t$ we note that by \eqref{varphi_SEP_1_rep},
\begin{align}
\label{chacabuco}
\max_{\substack{(x,y) \in V_N\\x\neq y }}  |\varphi_t^N(x,y)| 
    		\leq  \max_{\substack{(z,w) \in V_N \\ z \neq w}}  |\varphi_0^N(z,w)| + \sup_{t \geq 0} \max_{z \in \Lambda_{N-1}} |g_t^N(z,z+1)| \max_{\substack{(x,y) \in V_N\\x\neq y }} \mathbb{E}_{(x,y)} \left[\int_0^\infty \mathbb{1}(\mathcal X_{sN^2}^i \in \mathcal{D}_N^+) ds \right].
\end{align}
Observe that
\begin{equation}
\label{chena}
    T_N^i(x,y) := \mathbb{E}_{(x,y)} \left[\int_0^\infty \mathbb{1}(\mathcal X_{tN^2}^i \in \mathcal{D}_N^+) dt \right]
\end{equation}
corresponds to the expected occupation time of the diagonals $\mathcal{D}_N^+$ by the random walk $(\mathcal X_{tN^2}^i; t \geq 0)$. 
By \eqref{chacabuco}, in order to estimate $|\varphi_t^N(x,y)|$, we only need to estimate the simpler quantities
$|\varphi_0^N(z,w)|$ for every $(z,w)\in V_N$ with $z \neq w$,  $|g_t^N(z,z+1)|$ for every $ z \in \Lambda_{N-1}$ and $T_N^i(x,y) $  for every $(x,y) \in V_N \setminus \mathcal D_N $. For details on this, see equations (2.19), (2.20), Lemma 6.2. and Sections 6.1. and 6.2 of \cite{GJMN}.

For $\alpha\geq 2$, we would like to follow a similar strategy to the one outlined above. However, in this case, the Chapman-Kolmogorov equation for $\varphi^N_t$ is more complicated. In the case $\alpha=1$, the relation $\eta(x) = \eta(x)^2$ has as a consequence that no diagonal terms appear in the equation satisfied by $\varphi_t^N$. For $\alpha \geq 2$, this relation is no longer satisfied, and therefore the Chapman-Kolmogorov equation has an additional term - see Appendix \ref{remark_construction_G_N}. At first glance, it would be natural to extend $\varphi_t^N$ to the diagonal $\mathcal D_N$ by taking $\varphi_t^N(x,x)$ as equal to
\begin{equation}\label{eq:corr_1}
\varphi_t^N(x,x) : = \mathbb E_{\mu^N} [(\eta_{tN^2}(x)-\rho_t^N(x))^2].
\end{equation}
However, it turns out that a more convenient definition is to extend $\varphi_t^N$ as
\begin{equation}\label{eq:corr_real}
\varphi_t^N(x,x) := \mathbb E_{\mu^N} \Big[ \frac{\alpha}{\alpha-1} \eta_{tN^2}(x)(\eta_{tN^2}(x)-1) - \rho_t^N(x)^2\Big],
\end{equation}
 and remark here the importance of $\alpha$ being greater or equal to $2$ for this quantity to be well defined. Some motivations and reasons for this choice of defining the function $\varphi_t(x,x)$ are  given in  the Appendix \ref{remark_construction_G_N}. Extending $\varphi_t^N$ in this way, we can verify that $\varphi_t^N$ satisfies the equation
\begin{align}
\label{equations_G_t^N_theta}
\partial_t \varphi_t^N(x, y) 
		=  N^2\Delta^{i}_N \varphi_t^N(x, y) + g^{N}_t(x,x+1)  \mathbb{1}((x,y) \in \mathcal{D}_N^+),
\end{align}
where $\Delta^{i}_N$ is the operator defined in \eqref{op_A}. To simplify,  we will use the same notation as in the case $\alpha=1$ to the occupation time \eqref{chena} for this case, i.e. the case $\alpha\geq 2$. 

Observe that \eqref{equations_G_t^N_theta} generalizes \eqref{eq_varphi_alpha1} in a very convenient way, because the right-hand side is structurally the same; the only difference being the definition of the operator $\Delta_N^i$ which in nothing changes the strategy we followed to bound $\varphi^N_t$ in case $\alpha=1$. In particular, we have the analogous of \eqref{chacabuco} for $\alpha \geq 2$ with the slight  difference that now we need to take into account in the right-hand side of \eqref{chacabuco} the points $(z,w) \in V_N$ with $z=w$.

From here on, we separate the proof of the bounds in \eqref{results_bounds_corr_bulk} and \eqref{results_bounds_corr_bound} in two parts: for Part 1 we treat the case $\theta <2$; and for Part 2 we treat the other case, i.e. $\theta \geq 2$.\\

\textbf{\underline{Part 1: the case $\theta <2$}}
\\

We already saw that 
\begin{align}
\label{chacabuco_2}
\maxx_{\substack{(x,y) \in V_N\\x\neq y }}  |\varphi_t^N(x,y)| 
    		\leq \maxx_{(z,w) \in V_N}  |\varphi_0^N(z,w)| + \sup_{t \geq 0} \maxx_{z \in \Lambda_{N-1}} |g_t^N(z,z+1)| \maxx_{\substack{(x,y) \in V_N\\x\neq y }} T^i_N(x,y),
\end{align}
Using Lemma \ref{estimate_time_oc_absorbe}, the assumptions  \eqref{decay_corr_time0_bulk} and \eqref{decay_corr_time0_bound}, and  Lemma \ref{lemma_bounded_discrete_grad_rho}, we conclude that
\begin{align*}
\sup_{t \geq 0} \maxx_{(x,y) \in V_N} |\varphi_t^N(x,y)| \lsim \begin{cases}
    \frac{1}{N} + \frac{N^\theta}{N}, \textrm{ if } \theta \leq 0,\\
    \frac{1}{N} + \frac{N^\theta}{N^3}, \textrm{ if } \theta > 0,
    \end{cases}
\end{align*} and so, for $\theta < 2$,
\begin{equation}
\sup_{t \geq 0}\maxx_{(x,y) \in V_N} | \varphi_t^N(x,y)| \lsim \frac{1}{N}.
\end{equation}
Moreover, for $x = 1, N-1$,
\begin{align*}
\sup_{t \geq 0} \maxx_{y \in \Lambda_N} | \varphi_t^N(x,y)| \lsim 
    R_N^\theta, \textrm{ if } \theta \leq  2 \;.
\end{align*} 

For  the case $\theta \geq 2$ repeating the previous arguments we get the bound $    \frac{N^\theta}{N^3}$ and this is not enough for our results. For this reason  we need to  consider another random walk.
\\

\textbf{\underline{Part 2: the case $\theta \geq 2$}}
\\

Here we follow a different strategy to improve the bound for $T^i_N$ found previously, following the ideas presented in \cite{GJMN} for the case $\alpha = 1$, and extending the argument for $\alpha \in \mathbb{N}$. We rewrite  \eqref{equations_G_t^N_theta} as

{\begin{equation*}
    \partial_t \varphi_t^N(x, y) =  N^2 \mathfrak C^i_N \varphi_t^N(x, y) + \mathfrak V^i_N(x,y) \varphi_t^N(x, y) + g^N_t(x,x+1)  \mathbb{1}(y = x+1),
\end{equation*}} where $\mathfrak C^i_N$ is, as defined in \eqref{op_C}, the generator of the random walk 
$\{\widetilde{\mathcal{X}}_{t}^i ; t \geq 0\}$ and,
{\begin{equation*}
    \mathfrak V^i_N(x,y) = -\frac{\alpha N^2}{N^\theta}[\lambda^\ell \mathbb{1}(x=1) + \lambda^r \mathbb{1}(y=N-1)].
\end{equation*}}
By Feynman-Kac's formula, we have that 
\begin{align*}
    \varphi_t^N(x,y)= \widetilde{\mathbb{E}}_{(x,y)} \left[\varphi_0^N( \widetilde{\mathcal  X}_{tN^2}^i )e^{\int_0^t \mathfrak V^i_N({\widetilde{\mathcal X}}_{tN^2}^i)ds} + \int_0^t g_{t-s}^N( \widetilde{\mathcal  X}_{sN^2}^i) \mathbb 1( \widetilde{\mathcal  X}_{sN^2}^i\in\mathcal D^\pm)e^{\int_0^s \mathfrak V^i_N(r, \widetilde{\mathcal  X}_{rN^2}^i)dr}ds \right],
\end{align*}
\noindent where  $\widetilde{\mathbb{E}}_{(x,y)}$ denotes the  expectation given that $\widetilde{\mathcal  X}_{sN^2}^i$ starts from the point $(x,y)$. Now, since $\mathfrak V^i_N$ is negative, then
\begin{align} \label{bound1}
    \max_{\substack{(x,y) \in V_N\\ x\neq y}}\Big| \widetilde{\mathbb{E}}_{(x,y)} \left[\varphi_0^N( \widetilde{\mathcal  X}_{tN^2}^i )e^{\int_0^t \mathfrak V^i_N(s, \widetilde{\mathcal  X}_{tN^2}^i)ds} \right] \Big| \lesssim \max_{(z,w) \in V_N} |\varphi_0^N(z,w)|.
\end{align}
For the other term, by changing the integrals using Fubini's theorem and using the fact that $g_t^N$ and $\mathfrak V^i_N$ are both negative, we have that
{\begin{align*}
    \Big| \widetilde{\mathbb{E}}_{(x,y)} \left[\int_0^t g_{t-s}^N( \widetilde{\mathcal  X}_{sN^2}^i) e^{\int_0^s \mathfrak V^i_N( \widetilde{\mathcal  X}_{tN^2}^i)dr}ds \right] \Big| \leq \int_0^t \widetilde{\mathbb{E}}_{(x,y)} \left[-g_{t-s}^N( \widetilde{\mathcal  X}_{sN^2}^i) \right]ds.
\end{align*}} By similar arguments as in the case $\theta < 2$, we obtain that
{\begin{align} \label{bound2}
    \Big| \widetilde{\mathbb{E}}_{(x,y)} \left[\int_0^t g_{t-s}^N(\widetilde{\mathcal  X}_{sN^2}^i ) e^{\int_0^s \mathfrak V^i_N(\widetilde{\mathcal  X}_{tN^2}^i)dr}ds \right] \Big| &\leq \sup_{t \geq 0} \max_{z \in \Lambda_{N-1}} |g_t^N(z,z+1)|  \widetilde T_t^N(x,y),
\end{align}}
where
{\begin{align} \label{time_refl}
   \widetilde T_t^N(x,y):=\int_0^{t} \widetilde{\mathbb{E}}_{(x,y)} \left[ \mathbb{1}(\widetilde{\mathcal X}_{sN^2}^i \in \mathcal{D}_N^+) \right]ds \;.
\end{align}}
Observe that  we did not bound the last integral (from $0$ to $t$) by the integral over the interval from 0 to infinity and the reason is that the bound we will obtain for that time integral depends on $t$ and blows up when $t\to+\infty$. From Lemma \ref{lemma_bound_oc_temp_theta_big}
together with  \eqref{bound1} and \eqref{bound2}, we obtain 
\begin{align*}
    \sup_{t \in [0,T]} \max_{(x,y) \in V_N} |\varphi_t^N(x,y)| \lsim \frac{T + 1}{N},
\end{align*}
and, the same bound holds from $(x,y)\in\partial V_N$.
This concludes the proof.

\subsection{Proof of Proposition \ref{proposition_corr_decay_theta_less_1}} 
\label{proof_proposition_corr_decay_theta_less_1}

Recall that here we will only consider $\theta < 1$.
Since the result of Proposition \ref{proposition_corr_decay_theta_less_1} for $\alpha = 1$ and  $\theta < 0$  was not considered before, we will present a proof that works for every $\alpha \in \mathbb{N}$ and every $\theta < 1$. Let $\epsilon > 0$ and recall from the statement of Proposition \ref{proposition_corr_decay_theta_less_1} that we denote the set $\{1,\dots,\epsilon(N-1)\}$ by $\Lambda_N^{\epsilon,\ell}$. We want to show that, for every $\epsilon > 0$ and every $t \in (0,T]$,
\begin{align*}
    \max_{\substack{(x,y) \in \Lambda_N^{\epsilon,\ell} \times \Lambda_N \\ y \neq x}} |\varphi^N_t(x,y)| \lsim \left( 1 + \frac{1}{\sqrt{t}} \right)
    \frac{\epsilon}{N} + o\left(\frac{1}{N}\right).
\end{align*}

Since $\varphi^N_t$ is the solution to \eqref{eq_varphi_alpha1} then it  admits the representation \eqref{varphi_SEP_1_rep}. {As a consequence,} for every $t \in [0,T]$, we have that
{\begin{align*}
    \max_{\substack{(x,y) \in \Lambda_N^{\epsilon,\ell} \times \Lambda_N \\ y \neq x}} |\varphi^N_t(x,y)|
    &\leq 
    \max_{\substack{(x,y) \in \Lambda_N^{\epsilon,\ell} \times \Lambda_N \\ y > x}}  \left[ |\mathbb{E}_{(x,y)} [\varphi_0^N(\mathcal  X^i_{tN^2})]| + \Big| \mathbb{E}_{(x,y)}\Big[\int_0^t g^{N}_{t-s}( \mathcal  X^i_{sN^2})\mathbb 1(\mathcal  X^i_{sN^2}\in\mathcal D^\pm_N)\,ds\Big] \Big| \right]\\
    &\leq \max_{\substack{(z,w) \in V^\alpha_N} } |\varphi_0^N(z,w)| \max_{\substack{(x,y) \in \Lambda_N^{\epsilon,\ell} \times \Lambda_N \\ y > x}}  \mathcal{P}_{(x,y)}\left[ \mathcal  X^i_{sN^2} \notin \partial V_N \right]\\
    &+ \sup_{r \geq 0} \max_{z \in \Lambda_{N-1}} |g^{N}_r(z,z+1)| \max_{\substack{(x,y) \in \\ \Lambda_N^{\epsilon,\ell} \times \Lambda_N \\ x \neq y}}  T^i_N(x,y),
\end{align*}} where $V^\alpha_N$ was defined in \eqref{new_V}, $\{\mathcal  X^i_{t}\ ; \ t \geq 0\}$ is the bi-dimensional random walk on  $V_N$ with Markov generator $\Delta_N^{i}$ and $\mathcal{P}_{(x,y)}\left[ \mathcal  X^i_{tN^2} \notin \partial V_N \right]$ represents the probability that, starting from $(x,y)$, at time $tN^2$, the random walk $\{\mathcal  X^i_{t}\ ; \ t \geq 0\}$ is still not absorbed at the boundary. Recalling the proof of the estimate of $T^i_N(x,y)$ (see Lemma \ref{estimate_time_oc_absorbe}), one can easily see that
    \begin{equation*}
    \max_{\substack{(x,y) \in \Lambda_N^{\epsilon,\ell} \times \Lambda_N \\ y \neq x}}  T^i_N(x,y) \lsim \frac{\epsilon}{N}+ \frac{N^\theta}{N^3} \mathbb{1}(0 < \theta < 1 ) + \frac{N^\theta}{N} \mathbb{1}(\theta < 0).
    \end{equation*}
    Moreover, by Lemma \ref{lemma_bounded_discrete_grad_rho} and assumption \eqref{decay_corr_time0_bulk}, we have that
\begin{align} \label{bound_strip_triangle_correlation}
    \max_{\substack{(x,y) \in \Lambda_N^{\epsilon,\ell} \times \Lambda_N \\ x \neq y}} |\varphi^N_t(x,y)|\lsim \frac{1}{N} \max_{\substack{(x,y) \in \Lambda_N^{\epsilon,\ell} \times \Lambda_N \\ y >  x}}  \mathcal{P}_{(x,y)}\left[ \mathcal  X^i_{sN^2} \notin \partial V_N \right] + \frac{\epsilon}{N}+ \frac{N^\theta}{N^3} \mathbb{1}(0 < \theta < 1 ) + \frac{N^\theta}{N} \mathbb{1}(\theta < 0).
    \end{align}
  We are only left with estimating $\mathcal{P}_{(x,y)}\left[ \mathcal  X_{sN^2}^i \notin \partial V_N \right]$, when $(x,y) \in \Lambda_N^{\epsilon,\ell} \times \Lambda_N$ and $y > x$. This is the content of the {next} result.

\begin{proposition} \label{prop_prob_stip_triangle}
Let $\alpha \in \mathbb{N}$ and $\Lambda_N^{\epsilon,\ell}$ as defined in Proposition \ref{proposition_corr_decay_theta_less_1}. For every $t \in (0,T]$, there exists $\epsilon_0 > 0$ such that, for every $0 < \epsilon < \epsilon_0$,
    \begin{align} \label{bound_prob_on_band_of_boundary_of_triang}
        \max_{\substack{(x,y) \in \Lambda_N^{\epsilon,\ell} \times \Lambda_N \\ y > x } } \mathcal{P}_{(x,y)}\left[ \mathcal{X}^i_{tN^2} \notin \partial V_N \right] \lsim \frac{\epsilon}{\sqrt{t}},
    \end{align}
where $\mathcal{P}_{(x,y)}\left[ \mathcal  X^i_{tN^2} \notin \partial V_N \right]$ represents the probability that, starting from $(x,y)$, at time $tN^2$, the random walk $\{\mathcal  X^i_{t}\ ; \ t \geq 0\}$ is still not absorbed at the boundary.
\end{proposition}

Using the bound in \eqref{bound_prob_on_band_of_boundary_of_triang} and what we already proved in \eqref{bound_strip_triangle_correlation}, we conclude that
\begin{equation}
\max_{\substack{(x,y) \in \Lambda_N^{\epsilon,\ell} \times \Lambda_N \\ y \neq x}} |\varphi^N_t(x,y)| \lsim \left( 1 + \frac{1}{\sqrt{t}} \right)
    \frac{\epsilon}{N} + o\left(\frac{1}{N}\right),
\end{equation} as we wanted.

\begin{proof}[Proof of Proposition \ref{prop_prob_stip_triangle}]
We divide the proof in two cases: $\alpha = 1$ and $\alpha \geq 2$.
\\

\textbf{\underline{{Part 1: the case $\alpha=1$}}}
\\

For $\alpha=1$ the exclusion rule creates a natural order in the system. Indeed, starting the dynamics from a configuration $\eta$ and enumerating the particles from left to right, such order lasts for every $t \geq 0.$ This implies that, the leftmost particle of $\eta$ will remain the leftmost particle of the system until it is absorbed. This is the main idea behind the next argument.

Given $(x,y) \in \Lambda_N^{\epsilon,\ell}  \times \Lambda_N$ with $x < y$, then $\mathcal{P}_{(x,y)}\left[ \mathcal{X}^i_{tN^2} \notin \partial V_N \right]$ represents the probability that, at time $tN^2$, none of the two particles in the bulk were absorbed, knowing that one started close to the boundary, at the site $x \in \Lambda_N^{\epsilon,\ell}$. Roughly speaking, since $x < y$, if we track the movements, up to time $tN^2$, of the particle that started at $x$, i.e. the leftmost particle in the bulk, then, if it is absorbed with high probability, i.e. of the order  $1 - \frac{\epsilon}{\sqrt{t}}$, then the event $\{\mathcal{X}^i_{tN^2} \notin \partial V_N\}$ has to have a probability at least of order $\frac{\epsilon}{\sqrt{t}}$. The advantage of tracking just the leftmost particle on the bulk relies on the fact that we can compare it with a simple random walk, whose absorption probabilities are known. 

Let us formalize this argument. Recall the definition of $V^\alpha_N$ from \eqref{new_V}.  We also define $\overline{V}^{\alpha}_N = V^\alpha_N \cup \partial V^\alpha_N$ the closure of $V^\alpha_N$. The proof will follow by a sequence of definitions of other processes that can be related with $\{\mathcal{X}^i_{t N^2} \ ; \ t \geq 0\}$. We will divide our strategy in three steps. 
\\

\underline{\textbf{Step 1: Projecting $\{\mathcal{X}^i_{t}\ ; \ t \geq 0\}$ on the line}}
\\

Recall that $\mathcal{X}^i_{\cdot}:V^\alpha_N \to \mathcal{D}([0,T]; \overline{V}^{\ \alpha}_N)
$ is a process evolving on the triangle $\overline{V}^{\ \alpha}_N$. We can now project this process in $\overline{\Lambda}_N$, in the following way: 
 let $\overline{\Omega}_N:= \{ \eta \in \{0,1\}^{\overline{\Lambda}_N} \ | \ \eta(0) = 0, \eta(N) = 0, \textrm{ and } {\sum_{x \in\Lambda_N} \eta(x) = 2\}}$ the set of initial configurations of the process on the line and define $\xi^2_\cdot: \overline{\Omega}_N  \to \mcb{D}([0, T]; \{0,1\}^{\overline{\Lambda}_N})$ to be such that, for every $(x,y) \in V^\alpha_N$ setting $\eta = \eta_{(x,y)} \in \overline{\Omega}_N$ with $\eta(x) = 1$ and $\eta(y) = 1$ (and therefore $\eta(z) = 0$ for every $z \notin \{x,y\}$),
\begin{equation}
    \xi^2_t(\eta_{(x,y)})(z) = \begin{cases}
        0, \textrm{ if } z \neq \Pi_1\mathcal{X}^i_t(x,y) \textrm{ and } z \neq \Pi_2\mathcal{X}^i_t(x,y),\\
        1,\textrm{ if } z = \Pi_1\mathcal{X}^i_t(x,y) \textrm{ or } z = \Pi_2\mathcal{X}^i_t(x,y),
    \end{cases}
\end{equation} where again $\Pi_1$ and $\Pi_2$ are the projection functions on the first and second coordinates, respectively. Since there exists a bijection between $V^\alpha_N$ and $\overline{\Omega}_N$, the previous definition completely defines the process $\xi^2_\cdot$.

\underline{\textbf{Step 2: Construction of a lazy random walk that follows the movements of the leftmost particle}}
\\

To $\xi^2_\cdot$, which can be interpreted as a SEP(1) with only two particles and an absorbing boundary, we will associate another process on the line that will be defined as follows: let $\Tilde{\Omega}_N:= \{ \eta \in \{0,1\}^{\overline{\Lambda}_N} \ | \ \eta(0) = \eta(N) = 0, \textrm{ and } \sum_{x = 2}^{N} \eta(x) = 1\}$ the set of initial configurations on the line with only one particle that starts on the bulk and define $\xi^1_\cdot: \Tilde{\Omega}_N \to \mcb{D}([0, T]; \{0,1\}^{\overline{\Lambda}_N})$ as, for every $(x,y) \in V^\alpha_N$ setting $\eta = \eta_{(x)} \in \Tilde{\Omega}_N$ to be such that $\eta(x) = 1$ (and therefore $\eta(z) = 0$ for every $z \neq x$), then
    \begin{equation}
    \xi^1_t(\eta_{(x)})(z) = \begin{cases}
        0, \textrm{ if } z \neq \Pi_1\mathcal{X}^i_t(x,y), \\
        1,\textrm{ if } z = \Pi_1\mathcal{X}^i_t(x,y),
    \end{cases}
\end{equation} where $\Pi_1$ is the projection function on the first coordinate. Thus, $\xi^1_\cdot$ is the process that follows the left and right movements of $\mathcal{X}^i_\cdot$ in $V^\alpha_N$, i.e. it follows the particle in the system that starts at $x$.

To define $\xi^1_\cdot$ we are using the fact that, as we remarked above, the two particles on the line cannot exchange the order of their positions. We observe that, because of the exclusion rule, if, eventually, the clock of the leftmost particle rings and the jump is suppressed, $\xi^1_\cdot$ remains still until the clock of the leftmost particle rings again for an allowed movement. It is clear that $\xi^1_\cdot \leq \xi^2_\cdot$, in the sense that, for every $z \in \overline{\Lambda}_N$ and every $t \in [0,T]$, $\xi^1_t(z) \leq \xi^2_t(z)$.

    Then, given $(x,y) \in \Lambda_N^{\epsilon,\ell} \times \Lambda_N$ with $x < y$, we see that
    \begin{align*}
        \mathcal{P}_{(x,y)}\left[ \mathcal{X}^i_{tN^2} \notin \partial V_N \right] &\leq \mathcal{P}_{(x,y)} \left[ \textrm{ the leftmost particle of } \mathcal{X}^i_{\cdot} 	\textrm{ was not absorbed until time } tN^2 \right] \\
        &= \mathcal{P}_{\eta_{(x,y)}} \left[ (\xi^2_{tN^2}(\cdot))(0) = 0 \,,\,(\xi^2_{tN^2}(\cdot))(N) = 0 \right] \\
        &\leq \mathcal{P}_{\eta_{(x)}} \left[ (\xi^1_{tN^2}(\cdot))(0) = 0 \right].
    \end{align*}
    
    \underline{\textbf{Step 3: Comparison with a random walk that ignores the exclusion rule of the initial process}}
    \\
    
    Let $\tilde{\xi}^1_{\cdot}$ be the process that follows $\xi^1_\cdot$ up to the first time that a jump is suppressed. Here, the process $\tilde{\xi}^1_\cdot$ realizes the jump and starts following not the leftmost particle but the rightmost particle until a new jump for $\xi^1_\cdot$ was suppressed. Again, $\tilde{\xi}^1_\cdot$ realizes the jump returning to follow the leftmost particle, and so on. This new process $\tilde{\xi}^1_{\cdot}$ also satisfies $\tilde{\xi}^1_\cdot \leq \xi^2_\cdot$ and can be seen as the non-lazy version of $\xi^1_\cdot$ and that describes a continuous time simple symmetric random walk. 
    
    Observe that $$\mathcal{P}_{\eta_{(x)}} \left[ (\xi^1_{tN^2}(\cdot))(0) = 0 \right]\leq \mathcal{P}_{\eta_{(x)}} \left[ (\tilde{\xi}^1_{tN^2}(\cdot))(0) = 0 \right].$$ This is again a consequence of the fact that the two particles on the initial process can not exchange order and so, if the rightmost particle is absorbed at $x=0$ then for sure the leftmost was already absorbed.
    Then, since $0$ and $N$ are absorbing states, if $\xi^1_{tN^2}(\cdot)$ and $\tilde{\xi}^1_{tN^2}(\cdot)$ start with the same configuration, at each time $t$, the point where $\xi^1_{tN^2}(\cdot)$ has a non-zero value is always less or equal to the point where  $\tilde{\xi}^1_{tN^2}(\cdot)$ has a non-zero value. Therefore $\{(\xi^1_{tN^2}(\cdot))(0) = 0\} \subset \{(\tilde{\xi}^1_{tN^2}(\cdot))(0) = 0\}$. This implies that 
    $$\mathcal{P}_{(x,y)}\left[ \mathcal{X}^i_{tN^2} \notin \partial V_N \right] \leq \mathcal{P}_{\eta_{(x)}} \left[ (\xi^1_{tN^2}(\cdot))(0) = 0 \right] \leq \mathcal{P}_{\eta_{(x)}} \left[ (\tilde{\xi}^1_{tN^2}(\cdot))(0) = 0 \right] \leq \mathcal{P}_{\eta_{(x)}} \left[ \tau_1 > tN^2\right],$$ 
    where $\tau_1 = \inf\{t \geq 0 \ | \ (\tilde{\xi}^1_{tN^2}(\cdot))(0) = 1\}$ represents the first time that $\tilde{\xi}^1_{\cdot}$ hits $0$. So, since $x \in \Lambda_N^\epsilon$ and $\tilde{\xi}^1_{tN^2}(\cdot)$ describes a continuous time simple symmetric random walk, we have that, for fixed $t$, there exists $\epsilon_0>0$ such that, for every $0 < \epsilon \leq \epsilon_0$, $\mathcal{P}_{\eta_{(x)}} \left[ \tau_1 > tN^2\right]$ is of order $O(\frac{\epsilon}{\sqrt{t}})$.
\\

\textbf{\underline{{Part 2: the case $\alpha \geq 2$}}}
\\

Clearly in this case the natural ordering is lost, therefore we implement some changes in the previous argument.  Recall that we are working with an absorbing SEP($\alpha$) starting with only two particles, then for every pair $\{x,x+1\}$, for $x \in { \Lambda_{N-1} }$, the jump rates $c_{x,x+1}$ and $c_{x+1,x}$ can take values: $\alpha-1$, $\alpha$ or $2\alpha$ and, as $\alpha$ increases, the jump rates increase. This means that, the time at which a jump will occur, will be as shorter as bigger is the jump rate, namely  the value of $\alpha$. In particular, if $\alpha_1 \geq \alpha_2 > 1$ then the hitting time of the boundary $\partial V_N$ for the SEP($\alpha_2$) is greater or equal to the hitting time for  SEP($\alpha_1$), given that both processes start with the same skeleton of two particles in $\Lambda_N$. This implies that, to complete the proof of \eqref{bound_prob_on_band_of_boundary_of_triang}, it is enough to treat the case $\alpha=2$. 
\\

\textbf{\underline{The case $\alpha=2$}}
\\

Let $\mathcal{Z}_{\cdot}:V_N \to \mcb{D}([0,T]; \overline{V}_N)$ be the representation of  SEP($2$) with only two particles on the system. Fix $(x,y) \in V_N$ which will represent the starting point of $\mathcal{Z}_{\cdot}$, and, to simplify notation, let us denote $\mathcal{Z}_{\cdot}(x,y)$ only by $\mathcal{Z}_{\cdot}$. We remark that $\mathcal{Z}_{\cdot} \in \mcb{D}([0,T]; \overline{V}_N)$, which is a process that takes values on $\overline{V}_N$, can be interpreted as
$\mathcal{Z}_{\cdot} = \xi(\Gamma(\cdot)),
$ where $\xi$ is the skeleton of $\mathcal{Z}_{\cdot}$ and $\Gamma(\cdot)$ represents the Poisson point process associated with the marked Poisson point process $N(\cdot)$ of SEP($2$) given the initial configuration $(x,y)$. I.e., for every $s \in [0,T]$, $\Gamma(sN^2)$ is the number of jumps of the process up to time $sN^2$, which corresponds to counting how many marks, up to time $sN^2$, the marked Poisson point process $N(\cdot)$ had. Observe that, for every $t \in [0,T]$,
\begin{equation*}
    \Gamma_m(t) := \int_0^t dN^m \leq \Gamma(t) \leq \Gamma_M(t) := \int_0^t dN^M,
\end{equation*}
where $N^m$ is a Poisson process with parameter $1$ and $N^M$ is a Poisson process with parameter $4$. The choice of these parameters is due to the fact that, for every $x \in \Lambda_{N-1}$, the jump rates $c^i_{x,x+1}$ and $c^i_{x+1,x}$ only take three possible values: $1$, $2$ or $4$. So we choose the parameter of $N^m$ as the $\min\{1,2,4\}$ and of $N^M$ as the $\max\{1,2,4\}$. Then, denoting
\begin{align*}
    \tau_2 &= \inf\{t \geq 0 \ | \ \mathcal{Z}_{t} \in \partial V_N\},\\
    \tau_m &= \inf\{t \geq 0 \ | \ \xi(\Gamma_m(t)) \in \partial V^\alpha_N\},\\
    \tau_M &= \inf\{t \geq 0 \ | \ \xi(\Gamma_M(t)) \in \partial V^\alpha_N\},
\end{align*} we get that
\begin{align} \label{sandwich}
    \mathcal{P}_{(x,y)} \left[ \tau_m > tN^2\right] \leq \mathcal{P}_{(x,y)} \left[ \tau_2 > tN^2\right] \leq \mathcal{P}_{(x,y)} \left[ \tau_M > tN^2\right].
\end{align}

Since, the processes $\xi(\Gamma_m(\cdot))$ and $\xi(\Gamma_M(\cdot))$ have Poisson clocks with a parameter which is uniform on the triangle $V_N$, they now can be interpret as a continuous time simple symmetric random walk. Tracking the movements of the particle that started at site $x$ and everytime the particles meet and are on top of each other we start moving the particle that jumps from the top of the other, we can deduce that, for fixed $t$, there exists $\epsilon_0>0$ such that, for every $0 < \epsilon \leq \epsilon_0$, $\mathcal{P}_{(x,y)} \left[ \tau_m > tN^2\right]$ and $\mathcal{P}_{(x,y)} \left[ \tau_M > tN^2\right]$ are both of order $\frac{\epsilon}{\sqrt{t}}$. Remark that, since we are taking a bounded interval of time $[0,tN^2]$, the number of meetings between the two particles, when they get on top of each other, is finite, so the number of times that, eventually, we change what is the particle that we will follow next is finite, guarateeing that the process is well defined. From \eqref{sandwich}, we conclude that $\mathcal{P}_{(x,y)} \left[ \tau_2 > tN^2\right]$ is also of order $\frac{\epsilon}{\sqrt{t}}$.
    \end{proof}

\subsection{Proof of Lemma \ref{lemma_bound_corr_bound}}

\label{prob_estimates_Lemma1}
Developing the square in the expectation, using the symmetry of the integrating function on the square and applying Fubini's theorem, we get 
\begin{align} \label{last_bound_correlation}
    \mathbb{E}_{\mu^N} \left[ \left( \int_s^t \bar{\eta}_{sN^2}(x) ds\right)^2\right] =  2 \int_s^t \int_s^r \varphi^N_{v,r}(x,x) dv dr,
\end{align}
where, for $x, y \in \Lambda_N$, 
\begin{equation} \label{2space_2time_correlations}
    \varphi^N_{\nu,r}(x,y) = \mathbb{E}_{\mu^N} [\bar{\eta}_{\nu N^2}(x)\bar{\eta}_{rN^2}(y)] \;.
\end{equation}

Let us fix $v \in [s,t]$ and $x \in \Lambda_N$. For every $r \geq v$ and $y \in \Lambda_N$, a simple computation shows that $\Psi^N_r(y):=\varphi^N_{v,r}(x,y)$ is solution to
\begin{align} \label{system_of_correlation_2times}
    \begin{cases}
        \partial_r \Psi^N_r(y) = N^2 \Delta^i_N \Psi^N_r(y), \textrm{ if } y \in \Lambda_N,\\
        \Psi^N_v(y) = \varphi^N_v(x,y), \textrm{ if } y \in \Lambda_N,\\
        \Psi^N_r(0) = \Psi^N_r(N) = 0,
    \end{cases}
\end{align}
where, for every $f:{\overline{\Lambda}_N} \to \mathbb{R}$ such that $f(0) = f(N) = 0$ 
\begin{align} \label{def_F_operatorTheta}
N^2 \Delta^i_N f(y) = \begin{cases} \alpha N^2[f(y+1)+f(y-1)-2f(y)], \textrm{ if } y \notin \{1, N-1\},\\
\frac{\alpha \lambda^\ell}{N^\theta} N^2[f(0)-f(1)] + \alpha N^2[f(2)-f(1)],\textrm{ if } y = 1,\\
\frac{\alpha \lambda^r}{N^\theta} N^2[f(N)-f(N-1)] + \alpha N^2[f(N-2)-f(N-1)],\textrm{ if } y = N-1.
\end{cases}
\end{align}
Then the solution of the previous equation can be written in terms of the fundamental solution $P^{N,\theta}_{r}(x,y)$ of the initial value problem \eqref{system_of_correlation_2times} as: 
\begin{align} \label{2space_time_corr_representation}
    \Psi^N_r(y) = \sum_{\substack{z=1}}^{N-1} P^{N,\theta}_{r-v}(y,z) \mathbb{E}_{\mu^N}[\bar{\eta}_{v N^2}(y)\bar{\eta}_{v N^2}(z)] \;.
\end{align}
Plugging last identity in \eqref{last_bound_correlation} and  using  \eqref{results_bounds_corr_bound} and the fact that the occupation variables are bounded, we obtain
\begin{align} \nonumber
  \mathbb{E}_{\mu^N} \left[ \left( \int_s^t \bar\eta_{rN^2}(x) dr\right)^2\right]\lesssim
& \int_s^t \int_s^r \Big\{ P_{r-v}^{N,\theta}(x,x) + \sum_{\substack{z=1\\ z \neq x }}^{N-1} P^{N,\theta}_{r-v}(x,z)R_N^\theta\Big\} dv dr \\ \label{bound1_lemma_d_theta}
&\lsim \int_s^t \int_s^r \Big\{P_{r-v}^{N,\theta}(x,x) + R_N^\theta \Big\}dv dr,
\end{align}
where above we used the fact that $\displaystyle \sum_{\substack{z \in \Lambda_N \\ z \neq x }}P^{N,\theta}_{r}(x,z)$ is (uniformly in time) bounded by one. To finish the proof  we just need to estimate $\int_s^t \int_s^r P^{N,\theta}_{r-v}(x,x) dv dr$ for $x \in \{1,N-1\}$. 

Let us define $\tilde{P}_{r-v}^{N,\theta}(x,y)$ the fundamental solution of  {\eqref{system_of_correlation_2times}}  when $\lambda^\ell = \lambda^r = 1$. Remark that
\begin{equation}
f_{r-v}^{N,\theta}(x,y) := P_{r-v}^{N,\theta}(x,y) - \tilde{P}_{r-v}^{N,\theta}(x,y)
\end{equation} is the fundamental solution to
\begin{align} \label{system_of_correlation_2times_n}
    \begin{cases}
        \partial_r g_{r,v}^{N,\theta}(x,y)  = N^2 \Delta^i_N g_{r,v}^{N,\theta}(x,y) + N^2 K^{N,\theta} \tilde{P}_{r-v}^{N,\theta}(x,y), \textrm{ if } y \in \Lambda_N,\\
        g_{r,v}^{N,\theta}(x,y) = 0, \textrm{ if } y \in \Lambda_N,\\
        g_{r,v}^{N,\theta}(x,0) = g_{r,v}^{N,\theta}(x,N) = 0,
    \end{cases}
\end{align} where
\begin{align*}
K^{N,\theta} \tilde{P}_{r-v}^{N,\theta}(x,y) := - \frac{\alpha(1 - \lambda^\ell)}{N^\theta} \tilde{P}_{r-v}^{N,\theta}(1,y)\mathbb{1}(x=1) - \frac{\alpha(1 - \lambda^r)}{N^\theta} \tilde{P}_{r-v}^{N,\theta}(x,N-1)\mathbb{1}(y = N-1).
\end{align*} Thus, $\tilde{P}_{r-v}^{N,\theta}(x,y)$ is a probability and since  $\lambda^\ell, \lambda^r \leq 1$, then $K^{N,\theta} \tilde{P}_{r-v}^{N,\theta}(x,y) \leq 0$, and so, by the Maximum Principle,  Theorem \ref{th_maximum_principle_para}, we obtain 
\begin{equation} \label{bound_P_without_bc_with_bc}
f_{r-v}^{N,\theta}(x,y) \leq 0 \quad \Longleftrightarrow \quad P_{r-v}^{N,\theta}(x,y) \leq \tilde{P}_{r-v}^{N,\theta}(x,y).
\end{equation} Using Proposition \ref{proposition_bound_transition_prob} presented in the next section we have that,
for every $t \in [0,T]$ and $x \in \Lambda_N$
\begin{align*}
    P^{N,\theta}_{t}(x,x) \leq \tilde{P}^{N,0}_{t}(x,x) + \left( \frac{N^\theta}{\lambda^\ell}-1 \right) \tilde{P}^{N,0}_{t}(1,x) + \left( \frac{N^\theta}{\lambda^r}-1 \right) \tilde{P}^{N,0}_{t}(N-1,x), \quad \text{if } \theta \geq 0 
\end{align*}
and
\begin{align*}
    P^{N,\theta}_{t}(x,x) \leq \tilde{P}^{N,0}_{t}(x,x), \quad \text{if } \theta < 0 \;.
\end{align*}

Moreover, a simple computation similar to Lemma 4.3 of \cite{baldasso2017exclusion}, relying in a comparison to the case $\theta = 0$ and $\lambda^l = \lambda^r = \alpha$, shows that for  $x \in \{1,N-1\}$
$$\int_s^t \int_s^r \tilde{P}^{N,0}_{r-v}(x,1) dv dr\lesssim \frac{|t-s|}{N^2}$$ 
and by symmetry the same is true for 
$\int_s^t \int_s^r \tilde{P}^{N,0}_{r-v}(x,N-1) dv dr.$ From this we get that
\begin{align*} \mathbb{E}_{\mu^N} \left[ \left( \int_s^t \bar\eta_{rN^2}(x) dr\right)^2\right]\lsim  \frac{N^\theta}{N^2} |t-s|+ (t-s)^2R_N^\theta \;.
\end{align*}
From  the definitions of $R_N^\theta$ in \eqref{results_bounds_corr_bound} the proof of \eqref{bound_to_use_on_CLM} ends. To conclude \eqref{limit_to_use_on_CLM} we only have to observe that, by the definition of $d_N^\theta$, \eqref{bound_to_use_on_CLM} implies that
\begin{equation}
\mathbb{E}_{\mu^N} \left[ \left( \int_s^t d^\theta_N \bar\eta_{rN^2}(x) dr\right)^2\right]\lsim  |t-s| \begin{cases} N^{\theta - 1}\textrm{ if } \theta < 1 \\ 
N^{1-\theta} \textrm{ if } \theta > 1\end{cases} + \ \ (t-s)^2 (d^\theta_N)^2 R_N^\theta \;.
\end{equation} Since $(d^\theta_N)^2 R_N^\theta = N^{2(1-\theta)} \mathbb{1}(1 < \theta) + \frac{N^\theta}{N}\mathbb{1}(0 \leq \theta \leq 1) + N^{\theta} \mathbb{1}(-1 < \theta < 0) + \frac{1}{N} \mathbb{1}(\theta \leq -1)$, \eqref{limit_to_use_on_CLM} follows.

On the other hand, \eqref{bound_to_use_on_KC} follows once we prove that 
$$\int_s^t \int_s^r (d^\theta_N)^2 N^\theta\tilde{P}^{N,0}_{r-v}(x,1) dv dr\lesssim |t-s|^{1 + \delta_\theta}, $$ 
where  $\delta_\theta$ is the same as in the statement of the lemma. To obtain this, namely the analogous of equation (5.4) of \cite{GJMN}, we can simply repeat the argument used in Section 5.2 of \cite{GJMN}. To this aim we remark that 
$\tilde{P}^{N,0}_{r-v}(x,1) = P^{1,N,0}_{\alpha (r-v)}(x,1)$, where $P^{1,N,0}_{s}(x,y)$ is the unique solution of the initial value problem (5.4) of \cite{GJMN} taking $\theta=0$, i.e. fixed $x \in \Lambda_N$, we have
\begin{align*}
\begin{cases}\partial_t P^{1,N,0}_{t}(x,y) = N^2 \Delta^{1,i}_N P^{1,N,0}_{t}(x,y), \quad y \in \Lambda_N, t > 0,\\
 P^{1,N,0}_{t}(x,0) = P^{1,N,0}_{t}(x,N) = 0, \quad t > 0,\\
 P^{1,N,0}_0(x,y) = \delta_0(x-y), \quad y \in \Lambda_N,
 \end{cases}
\end{align*} where $\Delta^{1,i}_N$ coincide with the operator $\Delta^{i}_N$ when taking $\alpha = 1 = \lambda^\ell = \lambda^r$ and $\delta_0(x) = 1$ if $x=0$, otherwise it is equal to zero. The equality follows simply because they solve the same initial value problem, whose solution is unique.
%

\subsection{Proof of Lemma \ref{statementlemma2}}
\label{prob_estimates_Lemma2}

Recall that for $u \in[0,1]$ we defined $\iota^0_\epsilon(u):=\epsilon^{-1}\mathbb {1}_{(0,\epsilon]}(u)$ and $\iota^1_\epsilon(u):=\epsilon^{-1}\mathbb {1}_{[1-\epsilon,1)}(u)$. Here we will only give the details for the case $j=0$ since, for $j=1$, the proof is analogous. 
By expanding the square, using Fubini's Theorem and the definition of the density field $Y^N_s$, we obtain
\begin{align*}
    \mathbb{E}_{\mu_N}\left[ \left(\int_0^t Y^N_s(\iota_\epsilon^0) ds \right)^2  \right] &=  \frac{2}{\epsilon^2 N}  \sum_{x,y \in \Lambda_N^{\epsilon, \ell} } \int_0^t \int_0^s \varphi^N_{v,s}(x,y) dv ds,
\end{align*}
where $\varphi^N_{v,s}(x,y)$ was defined in \eqref{2space_2time_correlations}. Using the identity  \eqref{2space_time_corr_representation}, last display is equal to 
\begin{align}
  &\frac{2}{\epsilon^2 N} \sum_{x \in \Lambda_N^{\epsilon, \ell}} \int_0^t \int_0^s \varphi^N_{v,s}(x,x) dv ds + \frac{2}{\epsilon^2 N} \sum_{\substack{x,y \in \Lambda_N^{\epsilon, \ell} \\ y \neq x}} \int_0^t \int_0^s \varphi^N_{v,s}(x,y) dv ds\nonumber  \\ \label{eq:useful1}
    = &\frac{2}{\epsilon^2 N} \sum_{x \in \Lambda_N^{\epsilon, \ell}} \int_0^t \int_0^s  P^{N,\theta}_{s-v}(x,x) \mathbb{E}_{\mu^N}[(\bar{\eta}_{vN^2}(x))^2] dv ds +\frac{2}{\epsilon^2 N}  \sum_{x \in \Lambda_N^{\epsilon, \ell}} \int_0^t \int_0^s \sum_{\substack{z \in \Lambda_N \\ z \neq x}} P^{N,\theta}_{s-v}(x,z) \varphi^N_{v}(z,x) dv ds\\ \label{eq:useful2}
    + &\frac{2}{\epsilon^2 N}  \sum_{\substack{x,y \in \Lambda_N^{\epsilon, \ell} \\y \neq x}} \int_0^t \int_0^s P^{N,\theta}_{s-v}(x,y) \mathbb{E}_{\mu^N}[(\bar{\eta}_{vN^2}(y))^2] dv ds + \frac{2}{\epsilon^2 N} \sum_{\substack{x,y \in \Lambda_N^{\epsilon, \ell} \\y \neq x}} \int_0^t \int_0^s \sum_{\substack{z \in \Lambda_N \\ z \neq y}} P^{N,\theta}_{s-v}(x,z) \varphi^N_{v}(z,y) dv ds.
\end{align}

We remark that, for every $x \in \Lambda_N$, $\displaystyle\sum_{\substack{z \in \Lambda_N \\ z \neq x}} P^{N,\theta}_{s-v}(x,z) \leq 1$. Using \eqref{results_bounds_corr_bulk}, we can bound the rightmost term in \eqref{eq:useful1} by
\begin{align*}
    \frac{2}{N \epsilon^2}\left| \sum_{x \in \Lambda_N^{\epsilon, \ell}} \int_0^t \int_0^s \sum_{\substack{z \in \Lambda_N \\ z \neq x}} P^{N,\theta}_{s-v}(x,z) \varphi^N_{v}(z,x) dv ds \right| \leq \frac{2 t^2}{\epsilon} \sup_{v \in [0,T]} \max_{\substack{(x, z) \in V_N \\ z \neq x}} |\varphi^N_v(x,z)| \lsim \frac{t^2}{\epsilon N},
\end{align*}
which goes to zero when taking $N$ to infinity. Moreover, using  \eqref{results_bounds_Gt_bulk_close_bound}, we can bound the rightmost  term of \eqref{eq:useful2} by
\begin{align} \nonumber
    \frac{2}{N \epsilon^2} \left| \sum_{\substack{x,y \in \Lambda_N^{\epsilon, \ell} \\y \neq x}} \int_0^t \int_0^s \sum_{\substack{z \in \Lambda_N \\ z \neq y}} P^{N,\theta}_{s-v}(x,z) \varphi^N_{v}(z,y) dv ds \right| &\lsim N \int_0^t \int_0^s  \max_{\substack{(z,y) \in \Lambda_N \times\Lambda_N^{\epsilon,\ell} \\ z \neq y}} |\varphi^N_v(z,y)| dv ds \\ \label{bound_imp}
    &\lsim \epsilon \int_0^t \int_0^s \left(1 + \frac{1}{\sqrt{v}} \right) dv ds + o\left( \frac{1}{N}\right) \\ \nonumber
    &\lsim C_t \epsilon + o\left( \frac{1}{N}\right),
\end{align}
where $C_t$ is a constant that depends on $t$.
Since, in the last bound, the first term is uniformly bounded in $N$, this term will only go to zero when taking $\epsilon$ to zero.

For the remaining  terms, since the occupation variables are bounded  for every $x \in \Lambda_N$, we can bound the first term in \eqref{eq:useful1} and \eqref{eq:useful2} by
\begin{align} \label{term1}
    \frac{2}{N \epsilon^2} \left| \int_0^t \int_0^s \sum_{x \in \Lambda_N^{\epsilon, \ell}} P^{N,\theta}_{s-v}(x,x)\mathbb{E}_{\mu^N}[(\bar{\eta}_{vN^2}(x))^2]dv ds \right| \lsim \frac{1}{N\epsilon^2}\sum_{x \in \Lambda_N^{\epsilon, \ell}} \int_0^t \int_0^s P^{N,\theta}_{s-v}(x,x) dv ds
\end{align}
and
\begin{align} \label{term3}
    \frac{2}{N \epsilon^2} \left|\int_0^t \int_0^s \sum_{\substack{x,y \in \Lambda_N^{\epsilon, \ell} \\y \neq x}} P^{N,\theta}_{s-v}(x,y) \mathbb{E}_{\mu^N}[(\bar{\eta}_{vN^2}(y))^2] dv ds \right| \lsim \frac{1}{N \epsilon^2} \sum_{\substack{x,y \in \Lambda_N^{\epsilon, \ell} \\y \neq x}} \int_0^t \int_0^s P^{N,\theta}_{s-v}(x,y) dv ds,
\end{align} respectively. The idea now is to estimate $P^{N,\theta}_{t}(x,y)$ using $\Tilde{P}^{N,0}_{t}(x,y)$, where $\Tilde{P}^{N,0}_{t}(x,y)$ represents $\mathbb{P}[\mathcal{X}^i_{tN^2} = y | \mathcal{X}^i_{0} = x]$, where $\mathcal{X}^i_{tN^2}$ is the random walk defined in point 1. in the begining of Section \ref{_correlation_decay_proof} in the case we choose $\theta = 0$ and $\lambda^\ell = \lambda^r = \alpha$. To do this, we will use the maximum principles of Appendix \ref{appendix_maximum_principles}. Inspired by the bound for $P^{N,\theta}_{t}(x,y)$ proved for $\theta \geq 0$ in Lemma 4.2 of \cite{GJMN}, we will show the following estimates.

\begin{proposition} \label{proposition_bound_transition_prob}
Let $\{ \mathcal{X}^i_{tN^2} \ ; \ t \geq 0\}$  be the random walk on $\Lambda_N$ with infinitesimal generator $N^2\Delta^i_N$ which was defined in \eqref{def_F_operatorTheta} and let $P^{N,\theta}_{t}(x,y)$ be the transition probability for this random walk, i.e. for every $(x,y) \in \overline{V}_N $,
\begin{equation*}
    P^{N,\theta}_{t}(x,y) = \mathbb{P}_x[\mathcal{X}^i_{tN^2} = y] = \mathbb{P}[\mathcal{X}^i_{tN^2} = y | \mathcal{X}^i_{0} = x],
\end{equation*}
which coincides with the fundamental solution of \eqref{system_of_correlation_2times}. Denote by $\tilde{P}^{N,0}_{t}$ the transition probability of the random walk  $\{ \mathcal{X}^i_{tN^2} \ | \ t \geq 0\}$ when we take $\theta=0$ and $\lambda^\ell = \lambda^r = \alpha$. Then, for every $t \in [0,T]$ and $(x,y) \in V_N$, for $\theta \geq 0$,
\begin{align*}
    P^{N,\theta}_{t}(x,y) \leq \tilde{P}^{N,0}_{t}(x,y) + \left( \frac{N^\theta}{\lambda^\ell}-1 \right) \tilde{P}^{N,0}_{t}(1,y) + \left( \frac{N^\theta}{\lambda^r}-1 \right) \tilde{P}^{N,0}_{t}(N-1,y)],
\end{align*}
and, for $\theta < 0$,
\begin{align*}
    P^{N,\theta}_{t}(x,y) \leq \tilde{P}^{N,0}_{t}(x,y).
\end{align*}

\end{proposition}

Remark that Proposition \ref{proposition_bound_transition_prob} is valid for every $\alpha \in \mathbb{N}$, extending what was known for the case $\alpha = 1$ and $\theta \geq 0$ to the case $\alpha \geq 2$ with $\theta \geq 0$ as well as the case $\theta < 0$ for all $\alpha \in \mathbb{N}$.

\begin{proof}[Proof of Proposition \ref{proposition_bound_transition_prob}]


Let $\theta \in \mathbb{R}$ and fix $t_0 \in [0,T]$ and $y_0 \in \Lambda_N$. Define the function $h^{N,\theta}_{t_0,y_0} : \overline{V}_N \to \mathbb{R}$ to be such that, for $x \in \Lambda_N$,
\begin{align*}
    h^{N,\theta}_{t_0,y_0}(x) = P^{N,\theta}_{t_0}(x,y_0) - \tilde{P}^{N,0}_{t_0}(x,y_0),
\end{align*} and at the boundary we define it as
\begin{equation*}
\begin{cases}
h^{N,\theta}_{t_0,y_0}(0)= \left(\frac{N^\theta}{\lambda^\ell} - 1 \right) \tilde{P}^{N,0}_{t_0} (1,y_0) \textrm{ and } h^{N,\theta}_{t_0,y_0}(N)= \left(\frac{N^\theta}{\lambda^r} - 1 \right) \tilde{P}^{N,0}_{t_0}(N-1,y_0) \textrm{ if } \theta \geq 0\\
h^{N,\theta}_{t_0,y_0}(0)=h^{N,\theta}_{t_0,y_0}(N) = 0 \textrm{ if } \theta < 0.
\end{cases}
\end{equation*}
Using the fact that, for every $t \in [0,T]$ and $(x,y) \in V_N$, $P^{N,\theta}_{t}(x,y)$ and $\tilde{P}^{N,0}_{t}(x,y)$ are fundamental solutions of \eqref{system_of_correlation_2times} for $\theta \in \mathbb{R}$ and for $\theta = 0$ and $\lambda^\ell = \lambda^r=\alpha$, respectively, we get 
\begin{align*}
    0=\partial_t h^{N,\theta}_{t_0,y_0}(x) &= N^2 \Delta^i_N h^{N,\theta}_{t_0,y_0}(x)\\&+ \frac{\alpha N^2}{N^\theta} \left[(N^\theta-\lambda^\ell) \tilde{P}^{N,0}_{t}(1,y_0)\mathbb{1}(x=1) + (N^\theta-\lambda^r) \tilde{P}^{N,0}_{t}(N-1,y_0) \mathbb{1}(x=N-1) \right]  \mathbb{1}(\theta < 0).
\end{align*}
Since, for $\theta < 0$, $\frac{N^2(N^\theta - \lambda^j)}{N^\theta} \leq 0$ where $j \in \{\ell,r\}$, then,
by the maximum principle, Theorem \ref{th_max_min_principle_disc_elliptic_operators}, if $\theta \geq 0$ and Theorem \ref{th_maximum_principle} if $\theta < 0$, for every $x \in \overline{V}_N$, we have that, for every $\theta \in \mathbb{R}$,
$$h^{N,\theta}_{t_0,y_0}(x) \leq \max\{h^{N,\theta}_{t_0,y_0}(0),h^{N,\theta}_{t_0,y_0}(N)\}.$$
This then implies that, for every $t \in [0,T]$ and $(x,y) \in V_N$,
\begin{equation*}
P^{N,\theta}_{t}(x,y) \leq \tilde{P}^{N,0}_{t}(x,y) + \left(\frac{N^\theta}{\lambda^\ell} - 1 \right) \tilde{P}^{N}_{t} (1,y_0) + \left(\frac{N^\theta}{\lambda^r} - 1 \right) \tilde{P}^{N,0}_{t}(N-1,y), \quad \text{for } \theta \geq 0 
\end{equation*}
and
\begin{equation*}
P^{N,\theta}_{t}(x,y) \leq \tilde{P}^{N,0}_{t}(x,y), \quad \text{for } \theta < 0 \;,
\end{equation*} as we wanted to show.
\end{proof}

We conclude this section with the following auxiliary results.
\begin{lemma} \label{bounds_prod_sum_int_box}
Let $\theta < 1$ then the following holds:
\begin{enumerate}
\item [i.] For every $\epsilon > 0$ and $t \in [0,T]$ 
\begin{align} \label{proof_inequality_sum_int_probXX}
    \limsup_{N \to +\infty} \sum_{\substack{x \in \Lambda_N^{\epsilon, \ell}  }}  \int_0^t \int_0^s \tilde{P}^{N,0}_{s-v}(x,x) dv ds \lsim t \epsilon \;.
    \end{align}

\item [ii.]  For every $\epsilon > 0$ and $t \in [0,T]$ 
  \begin{align}   \label{proof_inequality_sum_int_2prob}
    \limsup_{N \to +\infty} \sum_{\substack{x \in \Lambda_N^{\epsilon, \ell} }} \int_0^t \int_0^s [ \tilde{P}^{N,0}_{s-v}(x,1) + \tilde{P}^{N,0}_{s-v}(x,N-1) ]dv ds \lsim t \epsilon.
\end{align}

\item [iii.]  For every $p \geq 1$ and $t \in [0,T]$ 
\begin{align} \label{bdp}
     \lim_{\epsilon \to 0}\limsup_{N \to +\infty} \frac{1}{\epsilon^p N} \sum_{x \in \Lambda_N^{\epsilon, \ell} } \int_0^t \int_0^s P^{N,\theta}_{s-v}(x,x) dv ds = 0.
\end{align}
\item [iv.] For any $t \in [0,T]$ 
\begin{align} \label{second_bound_lemma4_3}
    \lim_{\epsilon \to 0} \limsup_{N \to +\infty} \frac{1}{\epsilon^2 N} \sum_{\substack{x,y  \in \Lambda_N^{\epsilon, \ell}  \\ y \neq x}} \int_0^t \int_0^s  \tilde{P}^{N,0}_{s-v}(x,y) dv ds = 0 .
\end{align}
\item [v.] For any $t \in [0,T]$ 
\begin{align} \label{second_new_bound_lemma4_3}
    \lim_{\epsilon \to 0} \limsup_{N \to +\infty} \frac{N^\theta}{\epsilon} \sum_{x \in \Lambda_N^{\epsilon, \ell} }  \int_0^t \int_0^s \tilde{P}^{N,0}_{s-v}(x,1) + \tilde{P}^{N,0}_{s-v}(x,N-1) dv ds = 0 .
\end{align}
\end{enumerate}
We also note that the  same results hold by replacing $\Lambda_N^{\epsilon, \ell} $ by $\Lambda_N^{\epsilon, r}$.
\end{lemma}
\noindent
Combining Proposition \ref{proposition_bound_transition_prob} and Lemma \ref{bounds_prod_sum_int_box}, for any $t \in [0,T]$ we have that
\begin{align*}
    \lim_{\epsilon \to 0} \limsup_{N \to +\infty} \frac{1}{N \epsilon^2} & \sum_{\substack{x,y \in \Lambda_N^{\epsilon, \ell}  \\y \neq x}} \int_0^t \int_0^s P^{N,\theta}_{s-v}(x,y) dv ds \lesssim \lim_{\epsilon \to 0} \limsup_{N \to +\infty} \frac{1}{N \epsilon^2} \sum_{\substack{x,y  \in \Lambda_N^{\epsilon, \ell}  \\y \neq x}}\int_0^t \int_0^s \tilde{P}^{N,0}_{s-v}(x,y) dv ds \\
   + \frac{N^\theta}{\epsilon}   & \sum_{y \in \Lambda_N^{\epsilon, \ell} } \int_0^t \int_0^s [\tilde{P}^{N,0}_{s-v}(1,y) + \tilde{P}^{N,0}_{s-v}(N-1,y)] dv ds \mathbb{1}(0 \leq \theta < 1)  = 0
\end{align*} 
and the same holds for 
 $\Lambda_N^{\epsilon, r}$.
With this we  complete the proof of Lemma \ref{statementlemma2}. 
Indeed, the previous observation together with equation \eqref{bdp} imply that the terms on the right-hand side of \eqref{term3} and \eqref{term1}, respectively, also go to zero when taking the limit as $N \to +\infty$ and then as $\epsilon \to 0$, from which the proof is complete.

\begin{proof}[Proof of Lemma \ref{bounds_prod_sum_int_box}]
To show all the estimates above recall that for every $t \in [0,T]$ and $x,y \in \Lambda_N$ we can explicitly write $\tilde{P}^{N,0}_{t}(x,y)$ via the eigenvalues and eigenfunctions of the operator $N^2 \Delta^i_N$ , see also Lemma 4.3. of \cite{GJMN}. Indeed, 
\begin{equation} \label{probability_in_terms_of_eigenvalues_and_functions}
    \tilde{P}^{N,0}_{t}(x,y) = \sum_{l \in \Lambda_N}e^{-\alpha \lambda^N_l t} v^N_l(x) v^N_l(y),
\end{equation}
where for every $x, l \in \Lambda_N $, $v^N_l(x) = \sqrt{\frac{2}{N}} \sin\left(\frac{\pi l x}{N}\right)$ and $\lambda^N_l = 4N^2 \sin^2\left(\frac{\pi l}{2N}\right)$ are respectively the eigenfunctions and eigenvalues of $N^2\Delta^{1,i}_N$. \\
We start with item \textit{i.} 
For $x=y$ after two times integration of \eqref{probability_in_terms_of_eigenvalues_and_functions} we get
\begin{align*}
    \sum_{\substack{x \in \Lambda_N^{\epsilon, \ell} }}  \int_0^t \int_0^s \tilde{P}^{N,0}_{s-v}(x,x) dv ds = \sum_{l \in \Lambda_N} t^2 \psi( \alpha \lambda^N_l t)  \sum_{\substack{x \in \Lambda_N^{\epsilon, \ell} }}  \frac{2}{N} \sin^2\left(\frac{\pi l x}{N}\right),
\end{align*}
where $\psi(u) := \frac{e^{-u} - 1 + u}{u^2}$. We observe that, for every $u \geq 0$, $|\psi(u)| \leq \min\{1,\frac{1}{u}\}$, then
\begin{align*}
    \sum_{l \in \Lambda_N} t^2 \psi( {\alpha} \lambda^N_l t)  \sum_{\substack{x \in \Lambda_N^{\epsilon, \ell}}}  \frac{2}{N} \sin^2\left(\frac{\pi l x}{N}\right) & \lsim \sum_{l \in \Lambda_N} \dfrac{2t \epsilon }{\pi^2 \alpha l^2} \frac{ \frac{\pi^2 l^2}{4N^2}}{ \sin^2\left(\frac{\pi l}{2N}\right) } \;.
\end{align*}
Noticing that $\frac{x^2}{sin^2(x)}$ is bounded  for $0 \leq x \leq 2$  we finally have that 
\begin{equation*}
   \limsup_{N \to +\infty} \sum_{\substack{x \in \Lambda_N^{\epsilon, \ell}  }}  \int_0^t \int_0^s \tilde{P}^{N,0}_{s-v}(x,x) dv ds  \lsim    \limsup_{N \to +\infty}\sum_{l \in \Lambda_N} \dfrac{2t \epsilon }{\pi^2 \alpha l^2} \lsim t \epsilon \;.
\end{equation*}
Now we prove item \textit{ii.} Again we start with the expression \eqref{probability_in_terms_of_eigenvalues_and_functions} for $y=1$ and $y=N-1$. We observe that, for every $t \in [0,T]$ and $x \in \Lambda_N$, since $\sin\left( \frac{\pi l (N-1)}{N} \right) = -\cos\left( \pi l \right) \sin\left( \frac{\pi l}{N} \right)$, then
\begin{align*}
    \tilde{P}^{N,0}_{t}(x,1) + \tilde{P}^{N,0}_{t}(x,N-1) = \sum_{l \in \Lambda_N} \frac{2[1-\cos\left( \pi l \right)]}{N} e^{- \alpha \lambda^N_l t} \sin\left(\frac{\pi l x}{N}\right) \sin\left(\frac{\pi l}{N}\right).
\end{align*} 
Thus, integrating twice in time both sides above, we get 
\begin{align*}
    \sum_{\substack{x \in  \Lambda_N^{\epsilon, \ell} }} \int_0^t \int_0^s \tilde{P}^{N,0}_{s-v}(x,1) + \tilde{P}^{N,0}_{s-v}(x,N-1) dv ds = \sum_{l \in \Lambda_N}   t^2 \psi(\alpha \lambda^N_l t)  2[1-\cos\left( \pi l \right)] \sin\left(\frac{\pi l}{N}\right) \sum_{\substack{x \in  \Lambda_N^{\epsilon, \ell}}} \dfrac{1}{N} \sin\left(\frac{\pi l x}{N}\right).
\end{align*}
As before, using the expression of $\lambda^N_l  $ we can bound the left-hand side of the last display by
\begin{align*}
 \sum_{l \in \Lambda_N} \dfrac{4t \epsilon}{\alpha \pi^2 l^2} \dfrac{\frac{\pi^2 l}{4 N^2}}{\sin^{2}\left(\frac{\pi l}{2N} \right) }.
\end{align*}
Using again that  $\frac{x^2}{sin^2(x)}$ is bounded for $0 \leq x \leq 2$ we conclude that
\begin{equation*}
    \limsup_{N \to +\infty} \sum_{\substack{x \in \Lambda_N^{\epsilon, \ell} }} \int_0^t \int_0^s [ \tilde{P}^{N,0}_{s-v}(x,1) + \tilde{P}^{N,0}_{s-v}(x,N-1) ]dv ds \lsim \lim_{N \to + \infty} \sum_{l \in \Lambda_N} \dfrac{4t \epsilon }{\alpha \pi^2 l^2}  \lsim  t \epsilon \;.
\end{equation*}
Now we prove item \textit{iii.} It simply follows from  equations \eqref{proof_inequality_sum_int_probXX}, \eqref{proof_inequality_sum_int_2prob} and Proposition \ref{proposition_bound_transition_prob}.  For every $p \geq 1$ we can conclude that
\begin{align}
    \lim_{\epsilon \to 0} \lim_{N \to +\infty} \frac{1}{\epsilon^p N} \sum_{\substack{x \in  \Lambda_N^{\epsilon, \ell} } } \int_0^t \int_0^s P^{N,\theta}_{s-v}(x,x) dv ds \lsim \begin{cases} \displaystyle
    \lim_{\epsilon \to 0} \lim_{N \to +\infty} \frac{t [1+N^\theta]}{\epsilon^{p-1} N} \ \ \ \ \  \textrm{ if } 0 \leq \theta < 1 \\ \displaystyle
    \lim_{\epsilon \to 0} \lim_{N \to +\infty} \frac{t }{ \epsilon^{p-1} N} \ \ \ \ \ \ \ \ \ \ \  \textrm{ if } \theta < 0
     \end{cases} = 0 \;.
\end{align}
Now we prove item \textit{iv}. A simple computation shows that
\begin{align*}
    \frac{1}{N \epsilon^2} \sum_{\substack{x,y = 1 \\ y \neq x}}^{\epsilon(N-1)}  \int_0^t \int_0^s \tilde{P}^{N,0}_{s-v}(x,y) dv ds =
    \sum_{l =1}^{N-1} \alpha^2 t^2 \psi( \alpha \lambda^N_l t)  \sum_{\substack{x, y = 1 \\ y \neq x }}^{\epsilon(N-1)}  \frac{2}{N^2 \epsilon^2} \sin \left(\frac{\pi l x}{N}\right) \sin \left(\frac{\pi l y}{N}\right).
\end{align*}
Trying to recover a Riemann sum from the right-hand side of the last identity, we can write
\begin{align} \label{integral_terms_double_epsilon}
    \sum_{\substack{x, y = 1 \\ y \neq x }}^{\epsilon(N-1)}  \frac{2}{N^2} \sin \left(\frac{\pi l x}{N}\right) \sin \left(\frac{\pi l y}{N}\right) = &\int_0^\epsilon \int_0^\epsilon  \sin \left(\pi l z\right) \sin \left( \pi l w\right) dz dw\\
    &\sum_{\substack{x, y = 1 \\ y \neq x }}^{\epsilon(N-1)}  \frac{2}{N^2} \sin \left(\frac{\pi l x}{N}\right) \sin \left(\frac{\pi l y}{N}\right) - \int_0^\epsilon \int_0^\epsilon  \sin \left(\pi l z\right) \sin \left( \pi l w\right) dz dw,
\end{align}
and we remark that
\begin{align*}
    \frac{1}{\epsilon^2} \int_0^\epsilon \int_0^\epsilon  \sin \left(\pi l z\right) \sin \left( \pi l w\right) dz dw = \frac{1}{\epsilon^2} \left( \int_0^\epsilon  \sin \left(\pi l z\right) dz \right)^2 = \left( \frac{1 - \cos\left(\pi l \epsilon \right)}{\pi l \epsilon}\right)^2.
\end{align*}
Therefore,
\begin{align} \nonumber
    &\frac{1}{\epsilon^2} \sum_{l =1}^{N-1} t^2 \psi(\lambda^N_l t)  \int_0^\epsilon \int_0^\epsilon  \sin \left(\pi l z\right) \sin \left( \pi l w\right) dz dw\\ \label{splitting_terms_with_min}
    = &\sum_{l =1}^{\min\{N-1, (\epsilon \pi)^{-1} \}} t^2 \psi(\lambda^N_l t)  \left( \frac{1 - \cos\left(\pi l \epsilon \right)}{\pi l \epsilon}\right)^2 + \sum_{\substack{l =\min\{N-1, (\epsilon \pi)^{-1} \} \\ l \in \mathbb{N} }}^{N-1} t^2 \psi(\lambda^N_l t)  \left( \frac{1 - \cos\left(\pi l \epsilon \right)}{\pi l \epsilon}\right)^2.
\end{align}
For the leftmost term of \eqref{splitting_terms_with_min}: by a third order Taylor expansion of $\cos\left(\pi l \epsilon \right)$ around zero, the fact that $x^p \leq \sqrt{x}$, for every $p \geq 1$ and $x \in [0,1]$, also that $l \leq (\epsilon \pi)^{-1}$, i.e. $\pi l \epsilon \leq 1$ and that $\psi(u) \leq 1/u$, then, for each $l$ in the above conditions, there exists $\xi_l \in (0,\pi l \epsilon)$, such that
\begin{align*}
    \sum_{l =1}^{\min\{N-1, (\epsilon \pi)^{-1} \}} t^2 \psi(\lambda^N_l t)  \left( \frac{1 - \cos\left(\pi l \epsilon \right)}{\pi l \epsilon}\right)^2 &= \sum_{l =1}^{\min\{N-1, (\epsilon \pi)^{-1} \}} t^2 \psi(\lambda^N_l t) \left( \frac{\pi l \epsilon}{2} - \cos (\xi_l)\frac{(\pi l \epsilon)^2}{3!}\right)^2 \\
    &\lesssim \sum_{l =1}^{\min\{N-1, (\epsilon \pi)^{-1} \}} \frac{t}{\lambda^N_l} \sqrt{\pi l \epsilon} \\
    &\lesssim \sqrt{\epsilon} \sum_{l =1}^{N-1} \frac{ t}{(\pi l)^{3/2}} \frac{\pi^2 l^2}{4N^2 \sin^2 \left(\frac{\pi l}{2N}\right)} \lesssim t \sqrt{\epsilon}.
\end{align*}

For the rightmost term of \eqref{splitting_terms_with_min}, for $\epsilon > 0$ and close to zero, for $N\in \mathbb{N}$ sufficiently large, we have that $\min\{N-1, (\epsilon \pi)^{-1} \} = (\epsilon \pi)^{-1}$ and therefore
\begin{align*}
    \sum_{\substack{l =(\epsilon \pi)^{-1} \\ l \in \mathbb{N} }}^{N-1} t^2 \psi(\lambda^N_l t)  \left( \frac{1 - 1cos\left(\pi l \epsilon \right)}{\pi l \epsilon}\right)^2 &\lesssim \sum_{\substack{l =(\epsilon \pi)^{-1} \\ l \in \mathbb{N} }}^{N-1} \frac{t}{\lambda^N_l} \left( \frac{1}{\pi l \epsilon}\right)^2 = \sum_{\substack{l =(\epsilon \pi)^{-1} \\ l \in \mathbb{N} }}^{N-1} \frac{t}{\pi^4 l^4 \epsilon^2} \underbrace{\frac{\pi^2 l^2}{4N^2 \sin^2 \left( \frac{\pi l}{2N}\right)}}_{\leq 5} \\
    &\lesssim \sum_{\substack{l =(\epsilon \pi)^{-1} \\ l \in \mathbb{N} }}^{N-1} \frac{t}{\pi^4 l^4 \epsilon^2} \lesssim (\epsilon \pi)^{3-\delta}\sum_{\substack{l =(\epsilon \pi)^{-1} \\ l \in \mathbb{N} }}^{N-1} \frac{t}{\pi^4 l^{1+\delta} \epsilon^2} \lesssim t \epsilon^{1-\delta},
\end{align*} where $0< \delta < 1$. Putting these estimates together in \eqref{splitting_terms_with_min}, we finally obtain that 
\begin{align*}
    \frac{2}{\epsilon^2} \sum_{l =1}^{N-1} t^2 \psi(\lambda^N_l t)  \int_0^\epsilon \int_0^\epsilon  \sin \left(\pi l z\right) \sin \left( \pi l w\right) dz dw \lesssim t \max\{\sqrt{\epsilon},\epsilon^{1-\delta}\} \longrightarrow 0 \textrm{ as } N \to +\infty \textrm{ and then } \epsilon \to 0.
\end{align*}
Finally,
\begin{align*}
    &\frac{1}{\epsilon^2} \sum_{l =1}^{N-1} t^2 \psi(\lambda^N_l t) \left[ \sum_{\substack{x, y = 1 \\ y \neq x }}^{\epsilon(N-1)}  \frac{1}{N^2} \sin \left(\frac{\pi l x}{N}\right) \sin \left(\frac{\pi l y}{N}\right) - \int_0^\epsilon \int_0^\epsilon  \sin \left(\pi l z\right) \sin \left( \pi l w\right) dz dw \right] \\
    &\leq \frac{1}{\epsilon^2} \sum_{l =1}^{N-1} \frac{t}{(\pi l)^2} \underbrace{\frac{\pi^2 l^2}{4 N^2 sin^2\left( \frac{\pi l}{2 N} \right)}}_{\leq 5}\left| \sum_{\substack{x, y = 1 \\ y \neq x }}^{\epsilon(N-1)}  \frac{1}{N^2} \sin \left(\frac{\pi l x}{N}\right) \sin \left(\frac{\pi l y}{N}\right) - \int_0^\epsilon \int_0^\epsilon  \sin \left(\pi l z\right) \sin \left( \pi l w\right) dz dw \right|\\
    &\leq \frac{5 t}{\pi^2} \sum_{l =1}^{N-1} \frac{1}{l^2} \left| \sum_{\substack{x, y = 1 \\ y \neq x }}^{\epsilon(N-1)}  \frac{1}{\epsilon ^2 N^2} \sin \left(\frac{\pi l x}{N}\right) \sin \left(\frac{\pi l y}{N}\right) - \frac{1}{\epsilon^2}\int_0^\epsilon \int_0^\epsilon  \sin \left(\pi l z\right) \sin \left( \pi l w\right) dz dw \right|.
\end{align*}
To finish the argument, it is enough to show that
\begin{align} \label{to_finish_the_argument}
    \lim_{\epsilon \downarrow 0} \limsup_{N \to +\infty} \sum_{l =1}^{N-1} \frac{1}{l^2} \left| \sum_{\substack{x, y = 1 \\ y \neq x }}^{\epsilon(N-1)}  \frac{1}{\epsilon ^2 N^2} \sin \left(\frac{\pi l x}{N}\right) \sin \left(\frac{\pi l y}{N}\right) - \frac{1}{\epsilon^2}\int_0^\epsilon \int_0^\epsilon  \sin \left(\pi l z\right) \sin \left( \pi l w\right) dz dw \right| = 
    0.
\end{align}

A simple computation shows that since
\begin{align*}
    \frac{1}{\epsilon^2}\int_0^\epsilon \int_0^\epsilon  \sin \left(\pi l z\right) \sin \left( \pi l w\right) dz dw = \frac{1}{\epsilon^2}\sum_{\substack{x, y = 0 }}^{\epsilon(N-1)} \int_{\frac{x}{N}}^{\frac{x+1}{N}} \int_{\frac{y}{N}}^{\frac{y+1}{N}}  \sin \left(\pi l z\right) \sin \left( \pi l w\right) dz dw,
\end{align*}
and $\sin(x)$ is a Lipschitz continuous function, then, for every $l \in \Lambda_N$,

\begin{align*}
    &\sum_{l =1}^{N-1} \frac{1}{l^2} \left|\sum_{\substack{x, y = 1 \\ y \neq x }}^{\epsilon(N-1)}  \frac{1}{\epsilon ^2 N^2} \sin \left(\frac{\pi l x}{N}\right) \sin \left(\frac{\pi l y}{N}\right) - \frac{1}{\epsilon^2}\int_0^\epsilon \int_0^\epsilon  \sin \left(\pi l z\right) \sin \left( \pi l w\right) dz dw \right|\\
    = &\sum_{l =1}^{N-1} \frac{1}{l^2} \left|\frac{1}{\epsilon ^2} \sum_{\substack{x, y = 0 }}^{\epsilon(N-1)} \int_{\frac{x}{N}}^{\frac{x+1}{N}} \int_{\frac{y}{N}}^{\frac{y+1}{N}} \left[\sin \left(\frac{\pi l x}{N}\right) \sin \left(\frac{\pi l y}{N}\right) -  \sin \left(\pi l z\right) \sin \left( \pi l w\right) \right] dz dw - \sum_{\substack{x = 1}}^{\epsilon(N-1)}  \frac{1}{\epsilon ^2 N^2} \sin^2 \left(\frac{\pi l x}{N}\right) \right|\\
    \leq &\sum_{l = 1}^{N-1}\frac{2}{l^2 \epsilon } \sum_{\substack{x = 0 }}^{\epsilon(N-1)} \int_{\frac{x}{N}}^{\frac{x+1}{N}} \left|\sin \left(\frac{\pi l x}{N}\right) - \sin \left(\pi l z\right) \right| dz + \frac{\pi^2 }{6 \epsilon  N}\\
    \leq &\sum_{l = 1}^{N-1} \frac{2\pi }{l \epsilon } \sum_{\substack{x = 0 }}^{\epsilon(N-1)} \int_{\frac{x}{N}}^{\frac{x+1}{N}} \left(z - \frac{ x}{N} \right) dz + \frac{\pi^2 }{6 \epsilon  N} \lesssim \frac{\log(N)}{N} +  \frac{1}{\epsilon N} \longrightarrow 0 \textrm{ as } N \to +\infty,
\end{align*}
which proves \eqref{to_finish_the_argument}.
\newline
Item \textit{v.} For the final estimate we observe that the result immediately follows from \eqref{proof_inequality_sum_int_2prob} when $\theta < 0$.
For $0\leq \theta <1$ the idea is to improve the estimates done in \eqref{proof_inequality_sum_int_2prob}. Indeed we can write
\begin{align*}
\frac{N^\theta}{\epsilon}  \sum_{\substack{x \in \Lambda_N^{\epsilon, \ell} }} \int_0^t \int_0^s \tilde{P}^{N,0}_{s-v}(x,1) + \tilde{P}^{N,0}_{s-v}(x,N-1) dv ds   \lesssim \dfrac{1}{N^{1- \theta}} \sum_{l \in \Lambda_N} \dfrac{t}{\pi \alpha l} \lesssim  \dfrac{t}{ N^{(1 - \theta)/2}} \sum_{l \in \Lambda_N} \dfrac{1}{\pi \alpha l^{1+ (1 - \theta)/2}}
\end{align*}
 where in the first bound we used the same reasoning of item \textit{ii.} and that $\sin(2x) = 2 \sin(x) \cos(x)$ while for the last one we used that $l < N $. The result follows again by considering the limit as $N \to \infty$, since $1+ (1 - \theta)/2$ is bigger than one the series converges.
\end{proof}

\section{Results on occupation times}
\label{sec_results_occupation_times}

In this section we collect some of the results that were necessary regarding occupation times of all the random walks we used in the article.  The proof of  our results  uses an artefact that consists in comparing our random walk with another one for which explicit results are known. 
To that end, in the first subsection below we make a comparison with an absorbed random walk and in the following subsection we make a comparison with a reflected random walk.

\subsection{Comparison with an absorbed random walk}

\begin{lemma} \label{estimate_time_oc_absorbe}
Recall the function $T^{i}_N$ defined in \eqref{chena}. Then, for every $(x,y) \in V_N$
\begin{align*}
    T^i_N(x,y) 
    &\lesssim \begin{cases}
   \frac{1}{N} \mathbb{1}((x,y) \notin U_N) + \frac{1}{N^2} \mathbb{1}((x,y)\in U_N) +\frac{N^\theta}{N}, \textrm{ if } \theta \leq 0,\\
   \frac{1}{N} \mathbb{1}((x,y) \notin U_N) + \frac{1}{N^2} \mathbb{1}((x,y)\in U_N) + \frac{N^\theta}{N^3}, \textrm{ if } \theta > 0,
    \end{cases}
    \end{align*} where $U_N = \{(x,y) \in V_N \ | \ x = 1 \textrm{ or } y =N-1\}$.
\end{lemma}

\begin{proof}[Proof of Lemma \ref{estimate_time_oc_absorbe}]
To prove the result we will use the  random walk  $(\mathcal X_{tN^2}^i; t \geq 0)$ generated by the operator \eqref{op_A} with the choice $\lambda^\ell = \lambda^r = \alpha$ and $\theta = 0$.  Denote by $\mathcal T^{\textrm{a}}_N$ the expected occupation time of the diagonals $\mathcal{D}_N^+$ by that random walk.
A simple computation shows that $\mathcal T^{\textrm{a}}_N(x,y)$ is the solution of 
\begin{align*}
\begin{cases}
N^2 \Delta_N^{0,i} \mathcal T^{\textrm{a}}_N(x,y)=-\delta_{y=x+1}, \textrm{ if } (x,y) \in V_N\\
\mathcal T^{\textrm{a}}_N(x,y)= 0, \textrm{ if } (x,y) \in \partial V_N.
\end{cases}
\end{align*}
where $\Delta_N^{0,i}$ is the operator defined in  \eqref{op_A} with the choice $\lambda^\ell = \lambda^r = \alpha$ and $\theta=0$. Solving explicitly the previous system of linear equations, we obtain
\begin{align} \label{solution_T_N}
   \mathcal T^{\textrm{a}}_N(x,y) = \frac{(N-y)x}{N^2(\alpha N-1)}-\frac{1}{2N(\alpha N-1)}\mathbb{1}(y=x),
\end{align}
and therefore
\begin{equation} \label{bounds_T_0}
\max_{(x,y) \in V_N} \mathcal T^{\textrm{a}}_N(x,y) \lsim \frac{1}{N}, \quad \max_{x \in \Lambda_N} \mathcal T^{\textrm{a}}_N(x,N-1) \lsim \frac{1}{N^2} \quad \textrm{ and } \quad \max_{y \in \Lambda_N} \mathcal T^{\textrm{a}}_N(1,y) \lsim \frac{1}{N^2}.
\end{equation}
Now, let us consider the function
\begin{equation*}
    W^i_N(x,y) := T^i_N(x,y) - \mathcal T^{\textrm{a}}_N(x,y) + C^i_N(x,y),
\end{equation*}
where $C^i_N$ is given on $(x,y)\in \overline{V}_N$ by
\begin{equation*}
    C^i_N(x,y) = \Big(\frac{N^\theta}{\lambda^\ell+\lambda^r} - 1\Big) \min\limits_{(z,w) \in V_N^\alpha} \textrm{sgn} (\theta) T^{0,i}_N(z,w) \mathbb{1}( (x,y) \in V_N).
\end{equation*} 
Recall the expression of $\mathcal T^{\textrm{a}}_N$ given in \eqref{solution_T_N}. A simple computation shows that 
\begin{align*}
    \max_{(x,y) \in V_N} \mathcal T^{\textrm{a}}_N(x,y) = \begin{cases} \frac{2(N-\lfloor N/2 \rfloor)(\lfloor N/2 \rfloor)-N}{2N^2 (\alpha N - 1)}, \textrm{ if } N/2-\lfloor N/2 \rfloor < \lceil N/2 \rceil - N/2, \quad (\textrm{chosing in } \eqref{solution_T_N} \ x=y=\lfloor N/2 \rfloor),\\
    \frac{2(N-\lceil N/2 \rceil)(\lceil N/2 \rceil)-N}{2N^2 (\alpha N - 1)}, \textrm{ if } N/2-\lfloor N/2 \rfloor \geq \lceil N/2 \rceil - N/2, \quad (\textrm{chosing in } \eqref{solution_T_N} \ x=y=\lceil N/2 \rceil).
    \end{cases}
    \end{align*}
and
\begin{align*}
\min_{(x,y) \in V_N} \mathcal T^{\textrm{a}}_N(x,y) = \frac{1}{N^2(\alpha N-1)} \quad (\textrm{chosing in } \eqref{solution_T_N} \ x=1, y=N-1).
\end{align*}
Recall that $T^i_N$ is the solution of 
\begin{align*}
\begin{cases}
N^2 \Delta_N^{i} T^{i}_N(x,y)=-\delta_{y=x+1}, \textrm{ if } (x,y) \in V_N,\\
 T^{i}_N(x,y)= 0, \textrm{ if } (x,y) \in \partial V_N.
\end{cases}
\end{align*}
Then, a simple computation shows that $W^i_N$ is solution to
\begin{align}
    \begin{cases}
    N^2 \Delta_N^i W^i_N(x,y) +  \Big(N^2 \Delta_N^i-  N^2 \Delta_N^{0,i}\Big) \mathcal T^{\textrm{a}}_N+  N^2 \Delta_N^i C^i_N(x,y) = 0, \textrm{ if } (x,y) \in V_N,\\
    W^i_N(x,y) = 0, \textrm{ if } (x,y) \in \partial V_N,
    \end{cases}
\end{align}
A simple computation shows that for every $(x,y) \in V_N$
\begin{align*}
    \Big(N^2 \Delta_N^i-&  N^2 \Delta_N^{0,i}\Big) \mathcal T^{\textrm{a}}_N(x,y)+  N^2 \Delta_N^i C^i_N(x,y) \\&=(1+\mathbb{1}(y=x))N^2 \Big( \alpha\Big[1 - \frac{\lambda^\ell}{N^\theta} \Big] \mathcal T^{\textrm{a}}_N(1,y) \mathbb{1}(x=1) + \alpha \Big[1 - \frac{\lambda^r}{N^\theta}\Big] \mathcal T^{\textrm{a}}_N(x,N-1) \mathbb{1}(y=N-1) \Big)\\
    &+ (1+\mathbb{1}(y=x)) N^2 \Big( \frac{\lambda^\ell \alpha}{N^\theta}C^i_N(1,y) \mathbb{1}(x=1) + \frac{\lambda^r \alpha}{N^\theta}C^i_N(x,N-1) \mathbb{1}(y=N-1) \Big).
\end{align*}

Observe that the unique solution $f$ of 
\begin{equation*}
    \begin{cases}
    N^2 \Delta^i_N f(x,y) = 0, \textrm{ if } (x,y) \in V_N,\\
    f(x,y) = 0, \textrm{ if } (x,y) \in \partial V_N,
    \end{cases}
\end{equation*}
is $f(x,y) = 0$ for all $(x,y) \in \overline{V}_N$. Moreover, from the definition of $C^i_N$, for every $(x,y) \in V_N$, it holds
\begin{equation*}
   \Big(N^2 \Delta_N^i-  N^2 \Delta_N^{0,i}\Big) \mathcal T^{\textrm{a}}_N(x,y)+  N^2 \Delta_N^i C^i_N(x,y) \leq 0.
\end{equation*}
Therefore, $W^i_N$ is the solution of the initial value problem given by
\begin{align*}
    \begin{cases}N^2 \Delta^i_N W^i_N(x,y) \geq 0, \textrm{ if } (x,y) \in V_N,\\
    W^i_N(x,y) = 0, \textrm{ if } (x,y) \in \partial V_N.
    \end{cases}
\end{align*}
Applying a version of the maximum principle for discrete elliptic operators that are Markov generators, i.e. Theorem \ref{th_maximum_principle} below, we get for every $(x,y) \in \overline{V}_N$  that $ W^i_N(x,y) \leq 0$, i.e. 
\begin{align*}
    T^i_N(x,y) &\leq \mathcal T^{\textrm{a}}_N(x,y) - C^i_N(x,y) \\
    &\lesssim \begin{cases}
   \frac{1}{N} \mathbb{1}((x,y) \notin U_N) + \frac{1}{N^2} \mathbb{1}((x,y)\in U_N) +\frac{N^\theta}{N}, \textrm{ if } \theta \leq 0,\\
   \frac{1}{N} \mathbb{1}((x,y) \notin U_N) + \frac{1}{N^2} \mathbb{1}((x,y)\in U_N) + \frac{N^\theta}{N^3}, \textrm{ if } \theta > 0,
    \end{cases}
    \end{align*} where $U_N = \{(x,y) \in {V_N} \ | \ x = 1 \textrm{ or } y =N-1\}$.
This ends the proof. 
\end{proof}

\subsection{Comparison with a reflected random walk}

\begin{lemma} \label{lemma_bound_oc_temp_theta_big}
Recall \eqref{time_refl}. Then, for every $t \in [0,T]$,
\begin{equation*}
    \max_{\substack{(x,y) \in V_N \setminus \mathcal{D}_N}} \widetilde{T}^{N}_t(x,y)\lesssim \frac{t+1}{N} .
\end{equation*}
\end{lemma}
\begin{proof}[Proof of Lemma \ref{lemma_bound_oc_temp_theta_big}]
Recall that $\{\tilde{\mathcal{X}}_{tN^2} \ ; \ t \geq 0 \}$  represents a two-dimensional random walk  on  $V_N$ that jumps to every nearest-neighbor site at rate $\alpha$, except at the diagonal $\mathcal D_N^+$ where it jumps left/up at rate $\alpha$ and right/down at rate $\alpha-1$ and moreover, it is reflected at $\partial V_N$. Let $\tilde{\mathbb{E}}_{(x,y)}$ denote the  expectation given that $\tilde{\mathcal  X}_{tN^2}$ starts from the point $(x,y)$. From Dynkin’s formula, for every function $f: V_N \to \mathbb{R}$ and for every $(x,y) \in V_N \setminus \mathcal{D}_N$,
\begin{align}~\label{eq:imp}
    &0 = \tilde{\mathbb{E}}_{(x,y)} \left[M^N_t(f) \right] =  \tilde{\mathbb{E}}_{(x,y)} \left[ f(\tilde{\mathcal X}_{tN^2}) - f(\tilde{\mathcal X}_0) - \int_0^t N^2\mathfrak C_N f(\tilde{\mathcal X}_{sN^2}) ds \right].
    \end{align}
where $\mathfrak C^i_N$ is, as defined in \eqref{op_C}. From \eqref{eq:imp} we get 
    \begin{align*}
\tilde{\mathbb{E}}_{(x,y)} \left[\int_0^t N^2\mathfrak C_N f(\tilde{\mathcal X}_{sN^2}) ds \right] \leq \max_{z,w \in V_N} \{ f(z) - f(w) \}.
\end{align*}
For  the choice $f(x,y) = -(x-\frac 12)^2-(y-(N-\frac 12))^2$, a long but  elementary computation shows that for every $(x,y) \in V_N$:
\begin{align*} 
    N^2\mathfrak C_N f(x,y) = \begin{cases}
        -4\alpha N^2, \textrm{ if } |x-y| \geq 2 \textrm{ but } (x,y) \neq (1,N-1),\\
        -2 \alpha N^2, \textrm{ if } (x,y) = (1,N-1),\\
         N^2(2N - 4 \alpha -2) , \textrm{ if } |x-y| = 1,\\
        2\alpha N^2(2N - 1), \textrm{ if } y=x \textrm{ and } y,x\neq 1, N-1,\\
        2\alpha N^2(2N - 7), \textrm{ if } y=x=1 \textrm{ or } y=x=N-1.
    \end{cases}
\end{align*}

From last display, we conclude that
\begin{align*}
    \tilde{\mathbb{E}}_{(x,y)} \left[\int_0^t N^2\mathfrak C_N f(\tilde{\mathcal X}_{sN^2}) ds \right] 
    &= N^2(2N-4\alpha - 2) \int_0^t \tilde{\mathbb{E}}_{(x,y)} \left[ \mathbb{1}(\tilde{\mathcal X}_{sN^2} \in \mathcal{D}_N^+ \right] ds \\
    &+ 2\alpha N^2(2N-1) \int_0^t \tilde{\mathbb{E}}_{(x,y)} \left[ \mathbb{1}(\tilde{\mathcal X}_{sN^2} \in \mathcal{D}_N\setminus \{(1,1),(N-1,N-1)\}) \right] ds \\
    &+ 2\alpha N^2(2N-7) \int_0^t \left( \tilde{\mathbb{E}}_{(x,y)} \left[ \mathbb{1}(\tilde{\mathcal X}_{sN^2} = (1,1))\right] + \tilde{\mathbb{E}}_{(x,y)} \left[ \mathbb{1}(\tilde{\mathcal X}_{sN^2} = (N-1,N-1)) \right] \right) ds\\
    &- 2 \alpha N^2 \int_0^t \tilde{\mathbb{E}}_{(x,y)} \left[ \mathbb{1}(\tilde{\mathcal X}_{sN^2} = (1,N-1)) \right]ds - 4\alpha N^2 \int_0^t \tilde{\mathbb{E}}_{(x,y)} \left[ \mathbb{1}(\tilde{\mathcal X}_{sN^2} \in \mathcal{C} \right] ds,
\end{align*}
where $\mathcal{C} = \{(x,y) \in V_N \ | \ |x-y| \geq 2 \textrm{ and } (x,y) \neq (1,N-1)\}$. By noting that the time integral of the rightmost term in the first line of last display is equal to  $T^{i,N}_t(x,y)$, we conclude that 
\begin{align*}
    T^{i,N}_t(x,y)  &\leq -\frac{2\alpha N^2(2N-1)}{N^2(2N-4\alpha - 2)} \int_0^t \tilde{\mathbb{E}}_{(x,y)} \left[ \mathbb{1}(\tilde{\mathcal X}_{sN^2} \in \mathcal{D}_N\setminus \{(1,1),(N-1,N-1)\}) \right] ds \\
    &- \frac{2\alpha N^2(2N-7)}{N^2(2N-4\alpha - 2)}  \int_0^t \tilde{\mathbb{E}}_{(x,y)} \left[ \mathbb{1}(\tilde{\mathcal X}_{sN^2} \in \{(1,1),(N-1,N-1)\})\right] ds\\
    &+ \frac{ 2 \alpha N^2}{N^2(2N-4\alpha - 2)} \int_0^t \left(\tilde{\mathbb{E}}_{(x,y)} \left[ \mathbb{1}(\tilde{\mathcal X}_{sN^2} = (1,N-1))\right] + 2 \tilde{\mathbb{E}}_{(x,y)} \left[ \mathbb{1}(\tilde{\mathcal X}_{sN^2}\in \mathcal{C}\right] \right)ds \\
    &+\frac{1}{N^2(2N-4\alpha-2)}\max_{z,w \in V_N} \{ f(z) - f(w) \}.
\end{align*}
For $N \geq 2 \alpha +1$, since the first two terms of the last bound for $T^{i,N}_t(x,y)$ are negative, we have that
\begin{align*}
     \max_{\substack{(x,y) \in V_N}} T^{i,N}_t(x,y) &\lsim \frac{t}{N} + \frac{\displaystyle\max_{(x,y), (z,w) \in V_N} \left\{ \left(z-\frac{1}{2} \right)^2 + \left(w - N + \frac{1}{2} \right)^2 - \left(x-\frac{1}{2} \right)^2 - \left(y - N + \frac{1}{2} \right)^2 \right\}}{N^2(2N-4\alpha-2)}\\
    &\lsim \frac{t + 1}{N}.
\end{align*}
\end{proof}

\appendix

\section{Maximum principles} \label{appendix_maximum_principles}

\begin{theorem} \label{th_maximum_principle}
Let $\mathcal{E}$ be the Markov generator of the continuous time Markov chain $\{X_t\}_{t \geq 0}$ and denote by $\mathcal{D}(\mathcal{E})$ its domain. Let $\Omega$ be a discrete set with a non-empty $\partial \Omega$. If $f \in \mathcal{D}(\mathcal{E})$ with domain $\Omega$ is solution to
\begin{align*}
\begin{cases}
    \mathcal{E}f \geq 0 \textrm{ in } \Omega,\\
    f(x) = 0 \textrm{ in } \partial \Omega,
\end{cases}
\end{align*}
then
$
    f \leq 0 \textrm{ in } \Omega.
$
\end{theorem}

\begin{proof}
Let $f$ be the solution of 
\begin{align*}
\begin{cases}
    \mathcal{E}f = h \textrm{ in } \Omega,\\
    f(x) = 0 \textrm{ in } \partial \Omega,
\end{cases}
\end{align*}
with $h \geq 0$ in $\Omega$. Then, given the stopping time $\tau_{\partial \Omega} = \inf\{ t \geq 0 \ | \ X_t \in \partial \Omega \}$, $f$ can be represented, for every $x \in \Omega \cup \partial \Omega$ by
\begin{equation*}
    f(x) = -\mathbb{E}_x\Big[ \int_0^{\tau_{\partial \Omega}} h(X_t) dt \Big].
\end{equation*}
Since $h \geq 0$ in $\Omega$ by assumption, the result is a simple consequence of the previous formula.
\end{proof}

\begin{theorem} \label{th_max_min_principle_disc_elliptic_operators}
Let $A$ be a finite set. Define $\mathcal{F}(A)$ as the set of functions $f:A\to \mathbb{R}$. Consider a connected graph $(A,E)$ and define the non-empty subset of $A$, that we denote by $\partial A$, that is the set of vertices with degree one. Let $\mathcal{E}: \mathcal{F}(A) \to \mathcal{F}(A)$ be an operador of the form
\begin{align*}
\mathcal{E} f (\eta) = \sum_{\{\xi,\eta\} \in E} c(\eta,\xi) [f(\xi) - f(\eta)],
\end{align*} where $c(\cdot,\cdot)$ is a positive function.  If there exists $f \in \mathcal{F}(A)$ solution to
$
\mathcal{E}f = 0 \textrm{ in } A \setminus \partial A,
$
then
\begin{equation*}
    \max_{x \in A} f(x) \leq \max_{w \in \partial A} f(w) \quad \textrm{ and } \quad 
    \min_{x \in A} f(x) \geq \min_{w \in \partial A} f(w).
\end{equation*}
\end{theorem}

\begin{proof}
We prove the maximum case, since, to obtain the minimum, we only have to take $g=-f$ and the result follows.

If $f$ is constant, there is nothing to prove. So, assume this is not the case and let us proceed by contradiction. Since $A$ is finite, if $f$ was such that $\max_{x \in A \setminus \partial A} f(x) > \max_{w \in \partial A} f(w)$, then there would exist $y \in A \setminus \partial A$ such that $f(y) = \max_{x \in A} f(x)$ and $f(y) > f(w)$ for all $w \in \partial A$. Then
\begin{equation} \label{identity_to_see_then_average}
0 = \mathcal{E} f (y) = \sum_{\{\xi,y\} \in E} c(y,\xi) [f(\xi) - f(y)],
\end{equation} and, because $c > 0$, \eqref{identity_to_see_then_average} imply that 
\begin{equation}
\max_{x \in A} f(x) = \frac{1}{a_y}\sum_{\{\xi,y\} \in E} c(y,\xi) f(\xi),
\end{equation} where $a_y := \sum_{\{\xi,y\} \in E} c(y,\xi)$. Since the left-hand-side of last display is a weighted average, in order to an average to attain the maximum of a function, then all the points have to be equal to the maximum value. This means that, for all the vertices that are connected to $y$ by an edges, the maximum of $f$ is also attained there. Repeatting the argument now for this vertices, we obtain that, for all the vertices that are connected to them by an edge, the maximum of $f$ is also attained there, and so on. Because $G$ is connected, we know that for every two points of the graph there must exists a path that connect them. Therefore, by the previous reasoning, we showed that the maximum of the function has to be attained in $\partial A$, which is a contradiction.
\end{proof}

\begin{theorem} \label{th_maximum_principle_para}
Let $\mathcal{E}$ be the Markov generator of the continuous time Markov chain $\{X_t\}_{t \geq 0}$ and denote by $\mathcal{D}(\mathcal{E})$ its domain. Let $\overline{\Omega}$ be a discrete set and $\partial \Omega$ a non-empty subset of $\overline{\Omega}$. Let $\Omega = \overline{\Omega} \setminus \partial \Omega$. If $f:[0,T] \times \overline{\Omega} \to \mathbb{R}$ is a function that it is differentiable in time and that is solution to
\begin{align*}
\begin{cases}
    \partial_t f \leq \mathcal{E}f \textrm{ in } (0,T) \times \Omega,\\
    f(t,x) = 0 \textrm{ in } [0,T] \times \partial \Omega,\\
    f(0,x) = f_0(x), \textrm{ in } \Omega,
\end{cases}
\end{align*}
then
$f(y) \leq \max_{x \in \Omega}\{0,f_0(x)\}$, for every  $y \in [0,T] \times \overline{\Omega}.
$
\end{theorem}

The proof of the previous theorem can be obtained by adapting the proof of \eqref{th_maximum_principle} for the time-dependent case. It is a simple combination of Feynman-Kac's representation of the solution to the problem
\begin{align*}
\begin{cases}
    \partial_t f = \mathcal{E}f + h \textrm{ in } (0,T) \times \Omega,\\
    f(t,x) = 0 \textrm{ in } [0,T] \times \partial \Omega,\\
    f(0,x) = f_0(x), \textrm{ in } \Omega,
\end{cases}
\end{align*} where the function $h$ is non-positive.

\section{Details on the Chapman-Kolmogorov equation of $\varphi^N_t$, when $\alpha \geq 2$} \label{equation_correlations_computation}
For completeness we perform here some standard computations regarding one and two-point correlations, used in the proof of Proposition \ref{proposition_corr_decay}. For every $(x,y) \in V_N$, we have
\begin{align*}
\partial_t \varphi^N_t(x,y) &= \mathbb{E}_{\mu^N}[N^2\mcb{L}_N (\bar{\eta}_{tN^2}(x)\bar{\eta}_{tN^2}(y))]\\
&= \mathbb{E}_{\mu^N}[N^2\mcb{L}_N (\eta_{tN^2}(x)\eta_{tN^2}(y))] - \rho^N_t(y)\mathbb{E}_{\mu^N}[N^2\mcb{L}_N \eta_{tN^2}(x)] - \rho^N_t(x)\mathbb{E}_{\mu^N}[N^2\mcb{L}_N \eta_{tN^2}(y)],
\end{align*} by the forward Kolmogorov equation and the linearity of $\mcb{L}_N$.
It is worthy to compute $\mathbb{E}_{\mu^N}[\mcb{L}_N (\eta(x)\eta(y))]$ and $\mathbb{E}_{\mu^N}[\mcb{L}_N \eta(x)]$ (resp. $\mathbb{E}_{\mu^N}[\mcb{L}_N \eta(y)]$). We start with the latter.
\noindent
The action of the SEP($\alpha$) generator $\mcb{L}_N$ on the one-point correlation function is
\begin{align*}
\mcb{L}_N \eta(x) &= \alpha [\eta(x-1) - \eta(x)] \mathbb{1} (x \neq 1) + \alpha [\eta(x+1) - \eta(x)] \mathbb{1} (x \neq N-1) \\
&+ \frac{\alpha \lambda^\ell }{N^\theta}\left[\rho^\ell - \eta(1)\right] \mathbb{1} (x = 1) + \frac{\alpha \lambda^r }{N^\theta}\left[\rho^r - \eta(N-1)\right]\mathbb{1} (x = N-1) \;, 
\end{align*} 
for $x \in \Lambda_N$.
Similarly for $x,y \in \Lambda_N$ the action on the two-point correlation function can be conveniently written as
\begin{equation} \label{eq_action_gen_eta_x_eta_y}
\mcb L_N ( \eta(x) \eta(y) ) = \eta(x) \mcb L_N \eta(y) + \eta(y) \mcb L_N \eta(x) + \Gamma \left( \eta\left( x\right), \eta\left( y\right) \right)   \;,
\end{equation}
where 
\begin{align*}
 \Gamma \left( \eta\left( x\right), \eta\left( y\right) \right) = \begin{cases} 
 \frac{ \lambda^\ell \rho^\ell}{N^\theta} [\alpha -\eta(1)] +  \frac{ \lambda^\ell \eta(1)}{N^\theta} [\alpha - \rho^\ell ] + \alpha[\eta(1) + \eta(2)] - 2\eta(1)\eta(2) &\textrm{for } x=y=1, \\
 \alpha [\eta(x-1) + 2 \eta(x) + \eta(x+1)] - 2 \eta(x) [\eta(x-1) + \eta(x+1)] &\textrm{for } y = x \neq 1,N-1,  \\
 2 \eta(x)\eta(y) - \alpha [\eta(x) \eta(y)] &\textrm{for } y = x+1, \\ \frac{ \lambda^r \rho^r}{N^\theta} [\alpha -\eta(N-1)] +  \frac{ \lambda^r \eta(N-1)}{N^\theta} [\alpha - \rho^r] + \\ \alpha[\eta(N-1) + \eta(N-2)] - 2\eta(N-1)\eta(N-2) &\textrm{for } x=y=N-1, \\{0}&\textrm{otherwise}.
 \end{cases}
\end{align*}

\section{Extention of $\varphi^N_t$ to the diagonal} \label{remark_construction_G_N} 

The role of this section is to give two different approaches in order to extend the value of the correlation function to the diagonal $  \mathcal{D}_N$. We first start with an approach based on stochastic duality, while for the second one we use an analytic approach based on degree two functions.

\subsection{Stochastic Duality} \label{remark_duality}
Based on properties of duality (see \cite{carinci2013duality} for a survey on duality results for several boundary driven interacting systems), we show how to extend $\varphi^N_t$ to the diagonal $\mathcal{D}_N$. It is well known that the SEP($\alpha$) with open boundary has SEP($\alpha$) with absorbing boundary as its dual process with duality function $D:\Omega_N \times \Omega^{dual}_N \to \mathbb{R}$ given by
\begin{align} \label{duality_function_expression}
D(\eta, \hat{\eta}) = \left[\rho^\ell\right]^{\hat{\eta}(0)} \prod\limits_{x = 1}^{N-1} \frac{\eta(x)!(\alpha-\hat{\eta}(x))!}{[\eta(x)-\hat{\eta}(x)]!\alpha!} \mathbb{1}(\eta(x) \geq \hat{\eta}(x)) \left[\rho^r\right]^{\hat{\eta}(N)},
\end{align}
for every $(\eta,\hat{\eta}) \in \Omega_N \times \Omega^{dual}_N $, where $\Omega_N^{dual} = \mathbb{N} \times \{0, \ldots, \alpha \}^{\Lambda_N} \times \mathbb{N}$ is the state space of the absorbing dual process. If we now take $\hat{\eta} = \delta_x + \delta_y$ in \eqref{duality_function_expression}, we have that
\begin{align} \label{truncatedcorr}
    \mathbb{E}_{\mu^N}[D(\cdot, \delta_x + \delta_y)] = \begin{cases} \mathbb{E}_{\mu^N}\left[\frac{\eta(x)\eta(y)}{\alpha^2} \right], \textrm{ if } y \neq x\vspace{0,1cm}\\\vspace{0,1cm}
    \mathbb{E}_{\mu^N}\left[\frac{\eta(x)(\eta(x)-1)}{\alpha(\alpha-1)} \right], \textrm{ if } y = x.
    \end{cases}
\end{align}
A simple computation shows that in fact $ \varphi^N_t(x,y)$ as defined in  \eqref{time_dependent_correlation_def} for $x \neq y$ and in \eqref{eq:corr_real} for $x=y$ satisfies
\begin{align} \label{G_t_using_duality}
\varphi^N_t(x,y)=\alpha^2(\mathbb{E}_{\mu^N}[D(\eta_{tN^2}, \delta_x + \delta_y)] - \mathbb{E}_{\mu^N}[D(\eta_{tN^2}, \delta_x)]\mathbb{E}_{\mu^N}[D(\cdot, \delta_y)]) \;.
\end{align}
In other words, the function $ \varphi^N_t(x,y)$ can be written in a natural way in terms of the duality function \eqref{duality_function_expression} without distinguishing the case $x=y$. 

\subsection{Degree Two Functions}

Now we show an analytic argument to choose the extension of $\varphi^N_t$ to $\mathcal{D}_N$  as in \eqref{eq:corr_real}. In this subsection, for simplicity of the presentation, we neglect the boundary dynamics of the process and we explain the argument for the bulk dynamics. The general case, follows from adapting the ideas we present here.

Let us call $\tilde{\varphi}^N_t$ the extension of $\varphi^N_t$ to $\mathcal{D}_N$ as $\mathbb{E}_{\mu^N}[(\bar{\eta}(x))^2]$, i.e. for every $(x,y) \in V_N$
\begin{equation}
\tilde{\varphi}^N_t(x,y) = \begin{cases}
\varphi^N_t(x,y), \textrm{ if } y \neq x,\\
\mathbb{E}_{\mu^N}[(\bar{\eta}(x))^2], \textrm{ if } x = y.
\end{cases}
\end{equation} 
For $\alpha=1$ and since $\eta(x)\in\{0,1\}$ then there is no need to extend the correlation function to the diagonal $\mathcal{D}_N$.  Moreover, the Chapman-Kolmogorov equation for  $\varphi^N_t$ is very simple as we saw in \eqref{eq_varphi_alpha1}. Nevertheless, if $\alpha \geq 2$, the Chapman-Kolmogorov equation for $\tilde{\varphi}^N_t$ is not as simple. In fact, $\tilde{\varphi}^N_t$ is solution, for every $(x,y) \in V_N$ to
\begin{align} \label{system_correlation_alpha_geq_2}
\partial_t \tilde{\varphi}^N_t(x,y) &= N^2 \mathcal{H}_N \tilde{\varphi}^N_t(x,y) \\ \label{term_on_diag_y_xplus1}
&+ N^2 \left\{ 2 \tilde{\varphi}^N_t(x,x+1) - \Tilde{\chi}^{N,t}_\alpha(x,x+1) \right\}\mathbb{1}(y = x+1) \\ \label{term_on_diag_y_x}
&- N^2 \left\{ 4\tilde{\varphi}^N_t(x,x) - [   \Tilde{\chi}^{N,t}_\alpha(x,x+1) + \Tilde{\chi}^{N,t}_\alpha(x,x-1) ]\right\} \mathbb{1}(y = x),
\end{align}
where the operator $\mathcal{H}^i_N$ is the generator of a two dimensional random walk that jumps to each neighbor at rate $\alpha$, apart when it is on the diagonal $\mathcal{D}_N$ that jumps at rate $\alpha-1$ to each one of its neighbors, i.e. for every function $f: \overline{V}_N \to \mathbb{R}$ such that $f(x,y) = 0$ if $(x,y) \in \partial V_N$, and for every $(x,y) \in V_N$, 
\begin{align*}
    &\mathcal{H}_N f (x,y) = \begin{cases}
        \alpha[ f(x-1,y) + f(x+1,y) + f(x,y-1) +  f(x,y+1) -4f(x,y)], \textrm{ if } |x-y|\geq 1,\\
        2(\alpha - 1)[f(x-1,x) + f(x,x+1) - 2f(x,x)] , \textrm{ if } y = x,
    \end{cases}
\end{align*}
and, for every $(x,y) \in V_N$ such that $y \neq x$, $$\Tilde{\chi}^{N,t}_\alpha (x,y) = \rho^N_t(x)[\alpha - \rho^N_t(y)]+\rho^N_t(y)[\alpha - \rho^N_t(x)].$$

Since we have different signs for the extra terms that appear on the upper diagonal $\mathcal{D}^+_N$ and main diagonal $\mathcal{D}_N$, i.e. \eqref{term_on_diag_y_xplus1} and \eqref{term_on_diag_y_x}, and also they are not uniformly bounded in $N$, we observe that the argument used for the case $\alpha = 1$ can not be applied directly here. This motivates us to redefine the function on the diagonal values in such a way that it becomes the solution of an equation with a similar structure to  \eqref{eq_varphi_alpha1}. As we will see below, that function is exactly the function $\varphi^N_t$ defined in \eqref{eq:corr_real}.


We now observe that we can rewrite  \eqref{system_correlation_alpha_geq_2}, \eqref{term_on_diag_y_xplus1} and \eqref{term_on_diag_y_x} as:
\begin{align} \nonumber
\partial_t \tilde{\varphi}^N_t(x,y) &= N^2 \widetilde{\mathcal{H}}_N \tilde{\varphi}^N_t(x,y) \\ \label{extra_term_new_diag_y_xplus1}
&- N^2 \left\{ 2 \tilde{\varphi}^N_t(x,x) + 2 \tilde{\varphi}^N_t(x+1,x+1) + \Tilde{\chi}^{N,t}_\alpha(x,x+1) \right\}\mathbb{1}(y = x+1) \\ \nonumber
&- N^2 \left\{ 4\tilde{\varphi}^N_t(x,x) - [  \Tilde{\chi}^{N,t}_\alpha(x,x+1) +  \Tilde{\chi}^{N,t}_\alpha(x,x-1) ]\right\} \mathbb{1}(y = x),
\end{align}
where $\widetilde{ \mathcal{H}}_N$ is the operator given, for every $f:  {\overline{V}_N} \to \mathbb{R}$ and $(x,y) \in V_N$, by 
\begin{align*}
   \widetilde{\mathcal{H}}_N f (x,y) = \begin{cases}
        \alpha[  f(x-1,y) + f(x+1,y) + f(x,y-1) +  f(x,y+1) - 4f(x,y)], \textrm{ if } |x-y|\geq 2,\\
        \alpha[  f(x-1,x+1) +  f(x,x+2) - 2f(x,y)] \\
        \quad + (\alpha - 1)[f(x+1,x+1) + f(x,x) - 2f(x,x+1)], \textrm{ if } y = x+1,\\
        2(\alpha - 1)[ f(x-1,x) +f(x,x+1) - 2f(x,x)] , \textrm{ if } y = x.
    \end{cases}
\end{align*}
 With this new choice of operator acting on $\tilde{\varphi}^N_t$, we have corrected the sign problem but now the equation on the diagonal $\mathcal{D}^+_N$ is no longer closed for $\tilde{\varphi}^N_t$, i.e. the extra terms we get in \eqref{extra_term_new_diag_y_xplus1} also depend on $\tilde{\varphi}^N_t(x,x)$ and $\tilde{\varphi}^N_t(x+1,x+1)$, which are terms of the diagonal $\mathcal{D}_N $. Also, even though on the main diagonal we have a closed equation for $\tilde{\varphi}^N_t$, the extra terms that do not depend on $\tilde{\varphi}^N_t$ are non-negative. 

This motivates us to take as a candidate for $\varphi^N_t$ a function of the following form: $$\varphi^N_t(x,y)= C \tilde{\varphi}^N_t(x,y) + \mathbb{E}_{\mu^N}[f^N_t(x)] \mathbb{1}(y = x),$$ with $f^N_t(x) := A \eta(x)^2 + B \eta(x) + D$ (where $A$, $B$, $C$ and $D$ will be chosen later on). This choice is due to the fact that, since $\mcb{L}_N$ preserves the degree of  functions and we want to obtain a system of equations for degree two functions, then the function $f^N_t(x)$ should be of degree two. 

With this new definition we have 
\begin{align*}
    \partial_t \varphi^N_t(x,y) &= N^2 \Delta^i_N \varphi^N_t(x,y) + h_t(x,y),
\end{align*}
where $N^2 \Delta^i_N$ is the operator defined in \eqref{op_A} considered without the part that involve boundary terms and 
\begin{align} \nonumber
    h_t(x,y) &= - C [\tilde{\nabla}^+_N\rho^N_t(x)]^2 \mathbb{1}(y = x+1)\\ \nonumber
    &- N^2 [(\alpha-1)A - C][\mathbb{E}_{\mu_N}[\eta(x)^2] + \mathbb{E}_{\mu_N}[\eta(x+1)^2]] \mathbb{1}(y = x+1)\\ \nonumber
    &- N^2\{[(\alpha-1)B + \alpha C] [\rho^N_t(x) + \rho^N_t(x+1)] + 2(\alpha-1) D \}\mathbb{1}(y = x+1)\\ \nonumber
    &+ 2N^2 [(\alpha-1)A - C][\varphi^N_t(x,x+1) + \varphi^N_t(x-1,x) + (\rho^N_t(x+1) + \rho^N_t(x-1))\rho^N_t(x)]\mathbb{1}(y = x) \\ \nonumber
    &+ \alpha N^2 \{  [B + C + A][\rho^N_t(x-1) + \rho^N_t(x+1) + 2\rho^N_t(x)] + 4 D \}\mathbb{1}(y = x).
\end{align}
We observe that, since we want $h_t$ to not depend on $\varphi^N_t$, then it cannot  depend on $\mathbb{E}_{\mu_N}[\eta(x)^2]$ nor on $\mathbb{E}_{\mu_N}[\eta(x+1)^2]$, meaning that the second and fourth lines of last display have to be equal to zero. Then $(\alpha-1)A - C = 0$, i.e. $A = \tfrac{C}{\alpha-1}$. We can then simplify $h_t$ to
\begin{align} \nonumber
    h_t(x,y) &= - C [\tilde{\nabla}^+_N \rho^N_t(x)]^2 \mathbb{1}(y = x+1)\\ \label{choose_B}
    &- N^2\{[(\alpha-1)B + \alpha C] [\rho^N_t(x) + \rho^N_t(x+1)] + 2(\alpha-1) D \}\mathbb{1}(y = x+1)\\ \label{choose_B_main_diag}
    &+ \tfrac{\alpha}{\alpha-1}N^2 \{  [(\alpha-1)B + \alpha C][\rho^N_t(x-1) + \rho^N_t(x+1) + 2\rho^N_t(x)] + 4 D \}\mathbb{1}(y = x).
\end{align}
Now, by the fact that we want $h_t$ to be uniformly (in $N$) bounded, from \eqref{choose_B_main_diag} we need $D \leq 0$ and $(\alpha-1)B + \alpha C \leq 0$, but from \eqref{choose_B} we also need $D \geq 0$ and $(\alpha-1)B + \alpha C \geq 0$. To make these two requirements compatible, we finally obtain that $D = (\alpha-1)B + \alpha C = 0$, i.e. $D = 0$ and $B = -\tfrac{\alpha C}{\alpha-1}$. This implies that $h_t(x,y) = - C [\tilde{\nabla}^+_N \rho^N_t(x)]^2 \mathbb{1}(y = x+1)$. We impose that $C \geq 0$. 
For simplicity, we will take $C = 1$, and this coincides with the definition of  $\varphi^N_t$ from \eqref{eq:corr_real}.

\section{Proof of Lemma \ref{lemma_bounded_discrete_grad_rho}} \label{edp_result_bound_discrete_derivative}

The proof of last lemma follows exactly the same steps as in the proof of Lemma 6.2 of \cite{GJMN}, which was done for the case $\theta\geq 0$. For completeness and convenience of the reader we decided to present it here with the necessary adaptations to accommodate the case $\theta<0$.  In fact the proof we present below works for any $\theta<1$ and we note that the proof for $\theta>1$ follows exactly the same steps as the proof of Lemma 6.2 of \cite{GJMN}. Assume now that $\theta<1$. 
The idea of the proof  is to  find a sequence of functions  $\{\phi_N\}_N$, such that $\phi_N(t,\tfrac xN)$ is close to $\rho_t^N(x)$ with an error of order $O(N^{-1})$.  Therefore, we  consider a sequence of functions  of class  $C^4$ in space and for that we need to restrict to  initial profiles $\rho_0$ of class $C^6$. To this end let $\{\phi_N(t,u)\}_{N\geq 1}$ be the solution of 
 \begin{equation}\label{eq:robin}
\begin{cases}
\partial_t \phi_N(t,u)\;=\; \alpha \partial_u^2 \phi_N(t,u)\,, & \textrm{ for } t>0\,,\, u\in (0,1)\,,\\
\partial_u \phi_N(t,0^+) \;=\;\mu^\ell_N( \phi_N(t,0^+)-\rho^\ell)\,, & \textrm{ for } t>0\,,\\
\partial_u \phi_N(t,1^-) \;=\;\mu_N^r( \rho^r-\phi_N(t,1^-))\,, & \textrm{ for } t>0\,,\\
\phi_N(t,0)\;=\;\rho^\ell\,,\;\; \phi_N(t,1)=\rho^r\,,& \textrm{ for } t>0\,,\\
\phi_N(0,u)\;=\;g_N(u)\,,& u\in [0,1]\,,
\end{cases}
\end{equation}
where, for $j\in\{\ell, r\}$, we define $\mu_N^j=\tfrac {N\lambda^j}{N^\theta-\lambda^j}$, and $g_N$ is a function of class $C^6$ and that satisfies \eqref{H1_6} and \eqref{H2}. Repeating the proof of Section
6.4 of \cite{GJMN}, we see that $\phi_N \in C^{1,4}$, which is a consequence of the fact that the initial condition of the equation above is of class $C^6$ and $\phi_N$ satisfies \eqref{eq:robin}.

For $x\in \overline{\Lambda}_N$, let  $\gamma^N_t(x):=\rho^N_t(x)-\phi_N(t,\tfrac xN)\,$. A simple computation shows that $\gamma^N_t$ is solution of 
\begin{equation}\label{eq28}
\left\{
\begin{array}{ll}
 \partial_t\gamma_t^N(x)\,=\,(N^2\Delta^i_N \gamma_t^N)(x)+F_t^N(x)\,, \;\; x\in\Lambda_N\,,\;\;t \geq 0\,,\\
 \gamma_t^N(0)=0\,, \quad \gamma^N_t(N)=0\,, \;\;t \geq 0\,,\\
\end{array}
\right.
\end{equation}
where $\Delta^i_N$ was defined in \eqref{laplaciannn} and  $F_t^N(x)=(N^2\Delta^i_N- \alpha\partial_u^2)\phi_N(t,\tfrac xN)$. Since 
 $\phi_N(t,\cdot)$ is sufficiently regular, we are done if we show that 
$\big|\gamma^N_t(x)\big|\lesssim \tfrac 1N$.  From Duhamel's formula, we have 
\begin{equation*}
\gamma_t^N(x)\;=\; \mathbb E_x\Big[\gamma^N_0( X_{tN^2}^i)+\int_{0}^t F_{t-s}^N(  X_{sN^2}^i)\,ds\Big],
\end{equation*}
where   $\{ X_s^i,\, s\geq 0\}$ is  the random walk on $\overline{V}_N$ with generator $\Delta^i_N$, absorbed at the boundary $\{0,N\}$ and  $\mathbb E_x$ denotes the  expectation with respect to the probability induced by the generator $\Delta^i_N$ and the initial position $x$.   Therefore,
\begin{equation}\label{eq:estimate_lemma6}
\sup_{t\geq 0}\max_{x\in\Lambda_N}|\gamma_t^N(z)|\;\leq\; \max_{x\in \Lambda_N}|\gamma^N_0(x)|\;+\;\sup_{t\geq 0}\max_{x\in \Lambda_N}\Big|\mathbb E_x\Big[\int_{0}^t F_{t-s}^N( X_{sN^2}^i)\,ds\Big]\Big|.
\end{equation}

From \eqref{H1_6}, we have that
\begin{equation*}
\max_{x\in \Lambda_N}|\gamma^N_0(x)| = \max_{x\in \Lambda_N}|\rho^N_0(x)-g_N(\tfrac xN)| \lsim \frac{1}{N}.
\end{equation*} Then, it remains to analyse the rightmost term in last display. Note that 
\begin{equation}\label{eq:imp_1}
\Big|\mathbb E_x\Big[\int_{0}^t F_{t-s}^N( X_{sN^2}^i)\,ds\Big]\Big|\leq \int_{0}^{t}\sum_{z\in\Lambda_N}\mathbb P_x\Big[ X_{sN^2}^i=z\Big]| F_{t-s}^N(z)|\,ds.
\end{equation}
Since  $\phi_N\in C^4$,  then  $F_t^N(x)\lesssim{1/N^2}$ for $x\in\{2,\ldots,N-2\}$ and for any $t\geq 0$ and last display is bounded by 
\begin{equation}\label{713}
\frac{C}{N}+\sum_{k\in \{1,N-1\}} \mathbb E_x\Big[\int_{0}^{\infty}\textbf{1}_{\{ X_{sN^2}^i=k\}}\, ds\Big]\cdot|F_t^N(k)|.
\end{equation}
Last   expectation is  the average time spent by the random walk at the site $k$ until its absorption. This is the  solution of the elliptic equation  
\begin{equation*}
-N^2\Delta^i_N T^N(x)=\delta_{x=k},\forall x\in\Lambda_N
\end{equation*}
with null Dirichlet conditions $ T^N(0)=0$ and $T^N(N)=0.$
A simple computation shows that 
$$T^N(x)=\frac{N^\theta}{N^2}\Big[-A_N^i x+ B_N^i\Big]
$$
where $$A_N^i:=\frac{\lambda^r}{\lambda^\ell \lambda^r (N-2)+\alpha N^\theta(\lambda^\ell + \lambda^r))}\quad \textrm{and} \quad B_N^i:=\frac{1}{\lambda^\ell} \Big(1-\Big(\alpha-\frac{\lambda^\ell}{N^\theta}\Big)A_N^iN^\theta\Big).$$
From this it follows that
$\max_{x\in\Lambda_N}|T^N(x)|\lesssim \tfrac {N^\theta}{N^2}$. Now we analyse 
$\max_{k\in \{1,N-1\}}|F^N_t(k)|.$ We do the proof for the case $k=1$ and we leave the case $k=N-1$ to the interested reader. 
Note that 
\begin{equation*}
\begin{split}
F_t^N(1)&=(N^2\Delta^i_N- \alpha\partial_u^2)\phi_N(t,\tfrac 1N)\\
&=\alpha N^2 (\phi_N(t,\tfrac 2N)-\phi_N(t,\tfrac 1N))+\alpha N^{2-\theta}\lambda^\ell(\phi_N(t, 0)-\phi_N(t,\tfrac 1N))-\alpha \partial_u^2\phi_N(t,\tfrac 1N).
\end{split}
\end{equation*}
Now we use the regularity of $\phi_N$ and make a Taylor expansion to get
\begin{equation*}
\begin{split}
F_t^N(1)&=\alpha N \partial_u \phi_N(t,0^+) +O(1) + \alpha N^{2-\theta}\lambda^\ell\Big(\phi_N(t, 0)-\phi_N(t,0^+)-\frac{1}{N}\partial_u\phi_N(t,0^+)\Big)+O(N^{-\theta}).
\end{split}
\end{equation*}
If we now use the condition 
\begin{equation*}
\begin{split}
\alpha N (1-\frac{\lambda^\ell}{N^\theta}) \partial_u \phi_N(t,0^+)= \alpha N^{2-\theta}\lambda^\ell\Big(\phi_N(t,0^+)-\phi_N(t, 0)\Big),
\end{split}
\end{equation*}
which (by noting that $\phi_N(t,0)=\rho^\ell$) coincides with   $\partial_u \phi_N(t,0^+)=\mu^\ell_N( \phi_N(t,0^+)-\rho^\ell)$,
then we obtain 
\begin{equation*}
\begin{split}
\sup_{t\geq 0} |F_t^N(1)|&\lesssim 1 + N^{-\theta}.
\end{split}
\end{equation*}
Putting all the estimates together we find the bound for  \eqref{eq:estimate_lemma6}  given by 
\begin{equation*}
\sup_{t\geq 0}\max_{x\in\Lambda_n}|\gamma_t^N(x)|\lesssim\frac{1}{N}+\frac{N^\theta}{N^2}+\frac{1}{N^2}
\end{equation*}
from where the proof ends, since $\theta<1$. 

\begin{remark}
We observe that, for each $N \in \mathbb{N}$, the stationary solution of \eqref{eq:robin}, that we denote by $\bar{\rho}_{\mu^j_N}$, under the assumption that $\lambda^\ell = \lambda^r := \lambda$, is given by
\begin{equation} \label{expression_rho_stationary_robin}
\bar{\rho}_{\mu^j_N}(u):= \frac{\rho^r + \rho^l(1+\mu^j_N)}{2+\mu^j_N} + \frac{\mu^j_N(\rho^r - \rho^l) u}{2+\mu^j_N}.
\end{equation}
So, taking $g_N = \bar{\rho}_{\mu^j_N} + f \in C^6$ where $f$ is a $C_c^\infty[0,1]$ function, we have that $g_N$ satisfies \eqref{H1_6}. Indeed, using \eqref{expression_rho_stationary_robin} and the definition of $\mu^j_N$, we get that
\begin{align*}
\bar{\rho}_{\mu^i_N}(u)&= \frac{(N^\theta - \lambda)(\rho^r + \rho^l) + N \lambda \rho^l}{2(N^\theta - \lambda) + N \lambda} + \frac{N \lambda(\rho^r-\rho^l)u}{2(N^\theta - \lambda) + N \lambda}= N a_N u + b_N.
\end{align*} Therefore, because $f$ has compact support, we have that
\begin{equation*}
\partial^k_u g_N(u) = \partial^k_u \bar{\rho}_{\mu^i_N}(u),
\end{equation*} for $u \in \{0,1\}$ and $k=0,1,2,3$. Moreover, if we restrict $\rho_0^N$ to be such that $\rho_0^N(x) = g_N \left(\tfrac{x}{N}\right)$, then \eqref{H2} is trivially satisfied and we can find $\gamma \in C^6$ which satisfies \eqref{H1_3}. Indeed,
\begin{equation*}
\bar{\rho}_{\mu^j_N}(u) \xrightarrow[N \to +\infty]{} \bar{\rho}(u):= \begin{cases}
\rho^l + (\rho^r-\rho^l)u, \textrm{ if } \theta < 1,\\
\frac{\rho^r + (1+\lambda)\rho^l}{2 + \lambda} + \frac{\lambda(\rho^r-\rho^l)u}{2 + \lambda}, \textrm{ if } \theta = 1,\\
\frac{\rho^r + \rho^l}{2}, \textrm{ if } \theta > 1.
\end{cases}
\end{equation*} where the limit is taken uniformly in $u$. Taking $\gamma = \bar{\rho} + f$ we have that 
\begin{align*}
\frac{1}{N} \sum_{x \in \Lambda_N} \Big|\rho^N_0(x) - \gamma\left(\frac{x}{N} \right) \Big| = \frac{1}{N} \sum_{x \in \Lambda_N} \Big| \bar{\rho}_{\mu^i_N}\left(\frac{x}{N}\right) - \bar{\rho}\left(\frac{x}{N}\right)\Big| \leq \sup_{u \in [0,1]}|\bar{\rho}_{\mu^i_N}(u) - \bar{\rho}(u)| \xrightarrow[N \to +\infty]{} 0
\end{align*}
and so \eqref{H1_3} is satisfied.
\end{remark}

\section{Replacement Lemma} \label{sec_RL}

For a configuration $\eta\in\Omega_N$ and $x\in\Lambda_N$ we define the translation by $x$ of $\eta$ as  $(\tau_x\eta)(y)=\eta(x+y)$. Recall  \eqref{eq:averages}.

\begin{lemma}[Replacement Lemma]\label{RL_lemma_bulk}
Recall from Proposition \ref{proposition_corr_decay_theta_less_1} the definition of $\Lambda_N^{\epsilon,\ell},\Lambda_N^{\epsilon,r} $.  Fix  $x \notin \Lambda_N^{\epsilon,r}$ and let $\varphi:\Omega_N\to\mathbb R$ be a function whose support does not intersects the set of points in $\{x+1,\cdots, x+\epsilon N\}$. Then for any $\theta\in\mathbb R$ and for any  $t \in [0,T]$, it holds
\begin{equation} \label{RL_bulk_equation}
\lim_{\epsilon\to 0}\lim_{N\to+\infty}\mathbb{E}_{\mu^N} \left[\Big | \int_0^t \varphi(\tau_x\eta)\Big(\eta_{sN^2}(x) - \overrightarrow{\eta}^{\lfloor \epsilon N\rfloor}_{sN^2}(x)\Big) ds\Big| \right] =0.
\end{equation}
If  $x \notin \Lambda_N^{\epsilon,\ell}$ and  for $\varphi:\Omega_N\to\mathbb R$  a function whose support does not intersects the set of points in $\{x-\epsilon N, \cdots, x-1\}$, the same statement holds replacing $\overrightarrow{\eta}^{\lfloor \epsilon N\rfloor}_{sN^2}(x)$ by $\overleftarrow{\eta}^{\lfloor \epsilon N\rfloor}_{sN^2}(x)$. 
\end{lemma}

In the case $\varphi\equiv 1$, the last result was proved in Lemma 4.3 of \cite{FGS2022} but for sake of completeness we give here a sketch of the proof  of the more general result stated above, by following the strategy of the proof of the Lemma 4.3 of  \cite{FGS2022}.

\begin{proof}
Our starting point is to change the measure $\mu_N$ to a reference measure, which in fact should be the invariant state of the system that we do not know, but instead we consider another suitable measure that we define as follows. To this end, 
let $\varrho: [0,1] \to (0,1)$ be a Lipschitz function, bounded away from zero and one, and let
\begin{equation} \label{bernoulli_prod_measure_Lipschitz}
\nu^N_{\varrho(\cdot)}(\eta) := \prod_{x = 1}^{N-1} {\alpha \choose \eta(x)} \left(\varrho(\tfrac xN)\right)^{\eta(x)} \left(1-\varrho(\tfrac xN)\right)^{\alpha - \eta(x)}\end{equation}
be the inhomogeneous Binomial product measure of parameter $\varrho(\cdot)$.

From  the entropy  and Jensen's inequalities, the fact that
$e^{|x|} \leq e^{x}+e^{-x}$ and that for sequences of positive real numbers $(a_N)_N, (b_N)_N$ it holds
\begin{equation*}
\limsup_{N \rightarrow \infty} \dfrac{1}{N} \log (a_N + b_N ) = \max \left\lbrace \limsup_{N \rightarrow \infty} \dfrac{1}{N} \log (a_N  ), \; \limsup_{N \rightarrow \infty} \dfrac{1}{N} \log ( b_N ) \right\rbrace,
\end{equation*} together with Feynman-Kac’s formula, the expectation in \eqref{RL_bulk_equation} is bounded from above by 
\begin{align*}
& \frac{H(\mu^N | \nu^N_{\varrho(\cdot)})}{BN}+ t  \sup_{f \textrm{density}} \Big\{ \pm \langle \varphi(\tau_x\eta)(\eta(x) - \overrightarrow{\eta}^{\lfloor \epsilon N\rfloor}(x) ), f \rangle_{\nu^N_{\varrho(\cdot)}} + \tfrac{N}{B} \langle \mcb{L}_N \sqrt f, \sqrt f \rangle_{\nu^N_{\varrho(\cdot)}} \Big\},
\end{align*} 
where $B>0$.

Now we note that a bound on the entropy can be obtained as  $H(\mu^N | \nu^N_{\varrho(\cdot)}) \lesssim N$, see for example beginning of Section 4 of \cite{FGS2022}).
Moreover, we can use the estimate $N^2 \langle \mcb{L}_N \sqrt f, \sqrt f \rangle_{\nu^N_{\varrho(\cdot)}}$ given in Lemma 4.1 of \cite{FGS2022} (where the parameters $\epsilon,\gamma,\delta,\beta$ there have the correspondence given in \eqref{eq:choiceofparameters}). Putting this all together, we get that the expectation in the statement of the lemma is bounded from above by 
a constant times\begin{align*}
&\frac{1}{B}+ t\sup_{f \textrm{density}} \Big\{ \pm \langle \varphi(\tau_x\eta)(\eta(x) - \overrightarrow{\eta}^{\lfloor \epsilon N\rfloor}(x) ), f \rangle_{\nu^N_{\varrho(\cdot)}} - \frac{N}{B} D_{\nu^N_{\varrho(\cdot)}}(\sqrt{f}) \Big\} + \frac{1}{BN},
\end{align*} 
 where 
\begin{align*}
D_{\nu^N_{\varrho(\cdot)}}(\sqrt{f}) := D_{\nu^N_{\varrho(\cdot)}}^\ell(\sqrt{f}) + D^{bulk}_{\nu^N_{\varrho(\cdot)}}(\sqrt{f}) + D_{\nu^N_{\varrho(\cdot)}}^r(\sqrt{f})
\end{align*} with
\begin{align*}
D_{\nu^N_{\varrho(\cdot)}}^\ell(\sqrt{f}) &:= \int_{\Omega_N} \!\left[ \frac{\lambda^\ell \varrho^\ell \eta(1)}{N^\theta}\Big\{\!\sqrt{f}(\eta^{1,0}) - \!\sqrt{f}(\eta)\!\Big\}^2 + \frac{\lambda^\ell [\alpha-\varrho^\ell][\alpha- \eta(1)]}{N^\theta}\Big\{\sqrt{f}(\eta^{0,1}) - \sqrt{f}(\eta)\Big\}^2 \right]d\nu^N_{\varrho(\cdot)} \\
        D^{bulk}_{\nu^N_{\varrho(\cdot)}}(\sqrt{f}) &:=\sum_{x = 1}^{N-2} D_ {\nu^N_{\varrho(\cdot)}}^{x,x+1}(\sqrt{f}) +D_{\nu^N_{\varrho(\cdot)} }^{x+1,x}(\sqrt{f} ) \\&= \sum_{x = 1}^{N-2} \int_{\Omega_N} \!\eta(x)[\alpha-\eta(x+1)] \Big\{ \!\sqrt{f}(\eta^{x,x+1}) - \! \sqrt{f}(\eta) \! \Big\}^2 d\nu^N_{\varrho(\cdot)}\\
        &+\sum_{x = 1}^{N-2} \int_{\Omega_N} \! \eta(x+1)[\alpha-\eta(x)] \Big\{ \! \sqrt{f}(\eta^{x+1,x}) - \! \sqrt{f}(\eta) \! \Big\}^2 d\nu^N_{\varrho(\cdot)}
\end{align*} and the definition of  $D_{\nu^N_{\varrho(\cdot)}}^r(\sqrt{f})$  is analogous to the one of  $ D_{\nu^N_{\varrho(\cdot)}}^\ell(\sqrt{f}) $ by replacing $0$ and $1$ by $N$ and $N-1$, respectively, and also $\lambda^\ell$ and $\varrho^\ell$ by $\lambda^r$ and $\varrho^r$, respectively. We are now left with estimating
\begin{align*}
\langle \varphi(\tau_x\eta)(\eta(x) - \overrightarrow{\eta}^{\lfloor \epsilon N\rfloor}(x) ), f \rangle_{\nu^N_{\varrho(\cdot)}}
\end{align*}
for every $f$ density with respect to $\nu^N_{\varrho(\cdot)}$.
Note that 
\begin{align*}
\langle \varphi(\tau_x\eta)(\eta(x) - \overrightarrow{\eta}^{\lfloor \epsilon N\rfloor}(x) ), f \rangle_{\nu^N_{\varrho(\cdot)}} = \frac{1}{\lfloor \epsilon N\rfloor}\sum_{y=x+1}^{x+\lfloor \epsilon N\rfloor}\sum_{w=x+1}^{y-1} \langle \left[\eta(w)-\eta(w+1) \right]\varphi(\tau_x\eta), f \rangle_{\nu^N_{\varrho(\cdot)}}.
\end{align*} 
Since 
\begin{align} \nonumber
&\langle \left[\eta(w)-\eta(w+1) \right]\varphi(\tau_x\eta), f \rangle_{\nu^N_{\varrho(\cdot)}} \\ \label{bound_using_Young}
&= \frac{1}{2} \int_{\Omega_N} \left[\eta(w)-\eta(w+1) \right]\varphi(\tau_x\eta)[f(\eta) - f(\eta^{w,w+1})] d\nu^N_{\varrho(\cdot)} \\ 
\label{change_varibles_here}
&+ \frac{1}{2} \int_{\Omega_N} \left[\eta(w)-\eta(w+1) \right]\varphi(\tau_x\eta)[f(\eta) + f(\eta^{w,w+1})] d\nu^N_{\varrho(\cdot)},
\end{align} making a change of variables $\eta \mapsto \xi = \eta^{w,w+1}$ in \eqref{change_varibles_here} (and noting that the support of $\varphi$ does not overlap with the set of points where this change is done) and  splitting the state space $\Omega_N$ as is done in Lemma 4.3 of \cite{FGS2022}, we get
\begin{align*}
\eqref{change_varibles_here}=\frac{1}{2} \int_{\Omega_N} \left[\eta(w)-\eta(w+1) \right]\varphi(\tau_x\eta)\left( 1 - \frac{\varrho\left(\frac{w}{N}\right)[1-\varrho\left(\frac{w+1}{N}\right)]}{\varrho\left(\frac{w+1}{N}\right)[1 - \varrho\left(\frac{w}{N}\right)]}\right) f(\eta) d\nu^N_{\varrho(\cdot)}.
\end{align*} 
Since $\varrho(\cdot)$ is  Lipschitz  and bounded away from zero and one; the  occupation variables are bounded and $f$ is a density, the last display is bounded from above by a constant times $\Big| \varrho\left(\frac{w+1}{N}\right) - \varrho\left(\frac{w}{N}\right) \Big|$. 

Since $\eta(w)-\eta(w+1)= \frac{1}{\alpha}\left(\eta(w)[\alpha-\eta(w+1)]- \eta(w+1)[\alpha-\eta(w)]\right)$ and $x^2-y^2 =(x-y)(x+y)$, we get  that \eqref{bound_using_Young} is equal to 
\begin{equation} \begin{split}\label{display1}
&\frac{1}{2\alpha} \int_{\Omega_N} \eta(w)[\alpha-\eta(w+1)]\varphi(\tau_x\eta)[\sqrt{f}(\eta) - \sqrt{f}(\eta^{w,w+1})][\sqrt{f}(\eta) + \sqrt{f}(\eta^{w,w+1})] d\nu^N_{\varrho(\cdot)} \\ 
- &\frac{1}{2\alpha} \int_{\Omega_N} \eta(w+1)[\alpha-\eta(w)]\varphi(\tau_x\eta)[\sqrt{f}(\eta^{w+1,w}) - \sqrt{f}(\eta)][\sqrt{f}(\eta^{w+1,w}) + \sqrt{f}(\eta)] a_w d\nu^N_{\varrho(\cdot)}
\end{split}\end{equation} Using Young's inequality and then that $(x+y)^2 \leq 2(x^2+y^2)$, we can bound  \eqref{display1}  by
\begin{align*}
&\frac{1}{4\alpha A} \int_{\Omega_N} \eta(w)[\alpha-\eta(w+1)][\sqrt{f}(\eta) - \sqrt{f}(\eta^{w,w+1})]^2 d\nu^N_{\varrho(\cdot)} \\
&+ \frac{A}{2\alpha } \int_{\Omega_N} \eta(w)[\alpha-\eta(w+1)](\varphi(\tau_x\eta))^2[f(\eta) + f(\eta^{w,w+1})] d\nu^N_{\varrho(\cdot)} \\
&+ \frac{1}{4\alpha A} \int_{\Omega_N} \eta(w+1)[\alpha-\eta(w)][\sqrt{f}(\eta^{w+1,w}) - \sqrt{f}(\eta)]^2 d\nu^N_{\varrho(\cdot)} \\
&+ \frac{A}{2\alpha } \int_{\Omega_N} \eta(w+1)[\alpha-\eta(w)](\varphi(\tau_x\eta))^2[f(\eta^{w+1,w}) + f(\eta)] (a_w)^2 d\nu^N_{\varrho(\cdot)}.
\end{align*} where $A >0$ will be chosen later.

Putting together the previous  bounds, we get that  \eqref{bound_using_Young} and \eqref{change_varibles_here} are bounded from above by
\begin{align} \label{usefull_bound1_RL}
\langle \varphi(\tau_x\eta)\left[\eta(w)-\eta(w+1) \right], f \rangle_{\nu^N_{\varrho(\cdot)}} \lesssim \frac{1}{A} \left[D^{w,w+1}_{\nu^N_{\varrho(\cdot)}}(\sqrt{f}) + D^{w+1,w}_{\nu^N_{\varrho(\cdot)}}(\sqrt{f})\right] + A + \Big| \varrho\left(\frac{w+1}{N}\right) - \varrho\left(\frac{w}{N}\right) \Big|.
\end{align} 

From this it follows that 
\begin{align*}
& \pm \frac{1}{\lfloor \epsilon N\rfloor}\sum_{y=x-\lfloor \epsilon N\rfloor}^{x-1}\sum_{w=x}^{y-1} \langle \left[\eta(w)-\eta(w+1) \right][\alpha-\eta(x+1)], f \rangle_{\nu^N_{\varrho(\cdot)}} - \frac{N}{B} D_{\nu^N_{\varrho(\cdot)}}(\sqrt{f}) \\
&\lesssim \frac{1}{\lfloor \epsilon N\rfloor}\sum_{y=x-\lfloor \epsilon N\rfloor}^{x-1}\sum_{w=x}^{y-1} \left[\frac{1}{4A} \left[D^{w,w+1}_{\nu^N_{\varrho(\cdot)}}(\sqrt{f}) + D^{w+1,w}_{\nu^N_{\varrho(\cdot)}}(\sqrt{f})\right] - \frac{N}{B} D_{\nu^N_{\varrho(\cdot)}}(\sqrt{f}) \right.\\
&+ A \epsilon N + \frac{1}{\lfloor \epsilon N\rfloor}\sum_{y=x-\lfloor \epsilon N\rfloor}^{x-1}\sum_{w=x}^{y-1} \Big| \varrho\left(\frac{w+1}{N}\right) - \varrho\left(\frac{w}{N}\right) \Big| \\
&\lesssim \frac{1}{4A} - \frac{N}{B}  D_{\nu^N_{\varrho(\cdot)}}(\sqrt{f}) +  A \epsilon N + \frac{1}{\lfloor \epsilon N\rfloor}\sum_{y=x-\lfloor \epsilon N\rfloor}^{x-1}\sum_{w=x}^{y-1} \Big| \varrho\left(\frac{w+1}{N}\right) - \varrho\left(\frac{w}{N}\right) \Big| .
\end{align*} Choosing $A = \frac{B}{4N}$ and using the fact that $\varrho(\cdot)$ is Lipschitz, then 
\begin{align*}
& \limsup_{N \rightarrow \infty} \mathbb{E}_{\mu^N} \left[\Big | \int_0^t \varphi(\tau_x\eta)\Big(\eta_{sN^2}(x) - \overrightarrow{\eta}^{\lfloor \epsilon N\rfloor}_{sN^2}(x)\Big) ds\Big| \right] \lesssim \frac{1}{B} + \left[\frac{B \epsilon}{4} + \epsilon \right].
\end{align*} Finally, taking the limit $\epsilon \to 0 $ and then  $B \to \infty$, we are done. The proof of the other average to the left is completely analogous and we leave it to the reader. 
\end{proof}

\begin{acknowledgements}
B. S. thanks FCT/Portugal for the financial support through the PhD scholarship with reference 2022.13270.BD. P.G. thanks  Funda\c c\~ao para a Ci\^encia e Tecnologia FCT/Portugal for financial support through the
projects UIDB/04459/2020 and UIDP/04459/2020. M.J.~has been funded by  CNPq grant 201384/2020-5 and FAPERJ grant E-26/201.031/2022.  
This project has received funding from the European Research Council (ERC) under the European Union’s Horizon 2020 research and innovative programme (grant agreement No. 715734).
\end{acknowledgements}

\end{document}